\documentclass[11pt,a4paper,twoside]{elsarticle}
\usepackage[utf8]{inputenc}
\usepackage[T1]{fontenc}
\usepackage[english]{babel}
\usepackage{enumerate}
\usepackage{graphicx} 
\usepackage{booktabs}
\usepackage[left=2.5cm,right=2.5cm,top=3.5cm,bottom=3.5cm]{geometry}
\usepackage[font=small,labelfont=bf]{caption} 
\usepackage[lofdepth,lotdepth]{subfig}
\usepackage{mathtools}
\usepackage{amsmath,bm}
\usepackage{amssymb}
\usepackage{fancyhdr}
\usepackage{fourier} 
\usepackage{pgfplots} 
\usepackage{cancel}
\usepackage{graphicx} 

\pgfplotsset{
  log x ticks with fixed point/.style={
      xticklabel={
        \pgfkeys{/pgf/fpu=true}
        \pgfmathparse{exp(\tick)}
        \pgfmathprintnumber[fixed relative, precision=2]{\pgfmathresult}
        \pgfkeys{/pgf/fpu=false}
      }
  },
  log y ticks with fixed point/.style={
      yticklabel={
        \pgfkeys{/pgf/fpu=true}
        \pgfmathparse{exp(\tick)}
        \pgfmathprintnumber[fixed relative, precision=3]{\pgfmathresult}
        \pgfkeys{/pgf/fpu=false}
      }
  }
}
\usepackage{tikz}  
\usetikzlibrary{positioning}
\usepgfplotslibrary{external}
\tikzsetexternalprefix{External-TikZ/}
\tikzset{math3d/.style=
    {y= {(0.353cm,0.353cm)}, z={(0cm,1cm)},x={(1cm,0cm)}}}
\usepackage[colorlinks = true,
            linkcolor = blue,
            urlcolor  = blue,
            citecolor = black,
            anchorcolor = blue]{hyperref} 
            
\usepackage{titlesec}  
\titleformat*{\section}{\bf}
\titleformat*{\subsection}{\bf}

\newtheorem{remark}{Remark}[section]

\newcommand{\vect}[1]{\mathbf{#1}}
\newcommand{\tdomain}{\Omega_1(t) \cup \Omega_2(t)}

\newcommand{\bfx}{\mathbf{x}}
\newcommand{\bfv}{\mathbf{v}}
\newcommand{\bfu}{\mathbf{u}}

\newcommand{\bfV}{\mathbf{V}}
\newcommand{\bfw}{\mathbf{w}}

\newcommand{\bfn}{\mathbf{n}}
\newcommand{\bfg}{\mathbf{g}}
\newcommand{\bff}{\mathbf{f}}
\newcommand{\bfI}{\mathbf{I}}
\newcommand{\bfeps}{\bm{\epsilon}}

\newcommand{\mcN}{\mathcal{N}}
\newcommand{\lsf}{\phi}

\newcommand{\vel}{\mathbf{u}} 

\newcommand{\R}{\mathbb{R}}

\newcommand{\enstq}[2]{\left\{#1\mathrel{}\middle|\mathrel{}#2\right\}}

\newcommand{\norme}[1]{\left\Vert #1\right\Vert}
\newcommand{\norm}[1]{\left\Vert #1\right\Vert}

\renewcommand{\div}[1]{\nabla \cdot #1}

\newcommand{\grad}[1]{\nabla #1}
\newcommand{\gradS}[1]{\nabla_\Gamma #1}
\newcommand{\divS}[1]{ \nabla_\Gamma \cdot #1}
\newcommand{\diff}{\mathop{ }\mathopen{ }\mathrm{d}}
\newcommand{\restreinta}[1]{\mathclose{}|\mathopen{}_{#1}}
\newcommand{\prodscal}[2]{\left( #1\ ,\ #2\right)}
\newcommand{\jump}[1]{\llbracket #1 \rrbracket}
\newcommand{\dOmega}{\partial \Omega}
\newcommand{\dw}{\partial \Omega}

\newcommand{\w}{k} 
\newcommand{\mcF}{\mathcal{F}}

\newcommand{\meancurv}{\mathbf{H}}
\newcommand{\diffSurf}{D_\Gamma}

\newcommand{\Bspace}{\mathcal{V}_{h}}

\newcommand{\mesh}{\mathcal{T}_h}
\newcommand{\mcT}{\mathcal{T}}
\newcommand{\nablas}{\nabla_{\Gamma}}
\newcommand{\spacetimeS}{\mathcal{G}_{I,\Gamma}}
\begin{document}

\begin{frontmatter}
\title{A Cut Finite Element Method for two-phase flows with insoluble surfactants}
\author[1]{Thomas Frachon\corref{cor1}%
}
\ead{frachon@kth.se}
\address[1]{Department of Mathematics, KTH Royal Institute of Technology,
SE-100 44 Stockholm, Sweden}

\author[1]{Sara Zahedi}
\ead{sara.zahedi@math.kth.se}

\cortext[cor1]{Corresponding author}

\begin{abstract}
We propose a new unfitted finite element method for simulation of two-phase flows in presence of insoluble surfactant. The key features of the method are 1) discrete conservation of surfactant mass; 2) the possibility of having meshes that do not conform to the evolving interface separating the immiscible fluids; 3) accurate approximation of quantities with weak or strong discontinuities across evolving geometries such as the velocity field and the pressure. The new discretization of the incompressible Navier--Stokes equations coupled to the convection-diffusion equation modeling the surfactant transport on evolving surfaces is based on a space-time cut finite element formulation with quadrature in time and a stabilization term in the weak formulation that provides function extension. The proposed strategy utilize the same computational mesh for the discretization of the surface Partial Differential Equation (PDE) and the bulk PDEs and can be combined with different techniques for representing and evolving the interface, here the level set method is used.  Numerical simulations in both two and three space dimensions are presented including simulations showing the role of surfactant in the interaction between two drops.  
\end{abstract}


\end{frontmatter}

\section{Introduction}
Surfactants are present in many forms and play an important role in many applications, for example in drug delivery for their ability of solubilizing drugs and encapsulating microbubbles with a shell or in petroleum industries for promoting oil flow\cite{SchStaMar03}. Computer simulation is an important tool in such applications and is used to study the effects of surfactants on drop deformation, breakup or coalescence in multiphase flow systems.

There are several challenges that computational methods need to handle both with respect to the representation and evolution of the interface separating the immiscible fluids and the approximation of quantities such as the surfactant concentration, the fluid velocity, and the pressure. Here we assume that we know the PDEs that model two-phase flow systems in the presence of insoluble surfactant and that we also have a numerical representation technique for representing and evolving the interface. Given that, we propose a strategy based on finite element methods to accurately discretize the given PDEs both in two and three space dimensions so that the total surfactant mass is conserved and with accurate approximation of discontinuities across evolving interfaces without conforming the mesh to the time dependent domains.

Many computational methods have been developed for simulations of two-phase flow systems with surfactants, see for example \cite{MurTry08, TeiSonLowVoi11, KhaTor11}, where discretizations based on finite differences are proposed, \cite{BazAndMei06, SorTor18, PalSieTor19,HsuChuLaiTsa19} for strategies based on boundary integral methods, \cite{JaLo04} for a conservative discretization based on finite volume methods, and \cite{GaTo09, BaGaNu15} for different strategies using finite element methods.  

Our method falls into the class of so called unfitted finite element methods. We use a regular mesh that covers the computational domain but does not need to conform to outer or internal boundaries (in this case the interface). We refer to this mesh as the background mesh. On the background mesh we start with appropriate standard finite element spaces. We then define 1) active meshes corresponding to the subdomains occupied by the immiscible fluids and the interface; 2) active finite element spaces; 3) a weak formulation so that 
the proposed method is accurate and robust independent of the interface position relative the background mesh. 
Active spaces are defined as the restriction of the standard finite element spaces defined on the background mesh to the active mesh. From the active spaces we build pressure and velocity spaces with functions that are double valued on elements in a region around the interface. Hence functions in our function spaces can be discontinuous across interfaces separating the two fluids. Therefore, we don't need to regularize discontinuities by introducing regularized delta functions to approximate the surface tension force as for example in \cite{LaiTseHua08, MurTry08, KhaTor11} or/and regularized Heaviside functions to approximate densities and viscosities across the interface as for example in \cite{XuYaLo2012}. In the weak formulation of the Navier--Stokes equations coupled to the surfactant transport equation, we handle space and time similarly, physical interface conditions are imposed weakly, and we add stabilization terms in the weak form that are vital for the proposed method and its straightforward implementation.

This work is a further development of the space-time Cut Finite Element Method (CutFEM) proposed in \cite{FrZa19} to two-phase Navier--Stokes flows where also surfactant is present. We note that a CutFEM for the two-phase Navier--Stokes equations without surfactant has also been proposed in \cite{ClaKer19}. However, the discretization in \cite{ClaKer19} is based on implicit Euler, linear elements for both velocity and pressure, and a first order approximation of the mean curvature proposed in \cite{HanLarZah15} and is therefore limited to first order accuracy. For the discretization of the incompressible Navier--Stokes equations we here start with a regular mesh and the inf-sup stable Taylor-Hood element P2-P1 pair for velocity and pressure. See \cite{GuOl18} and \cite{KiGr16} for stability results. We also use an integration by part result in the presence of surfactant and avoid explicit computation of the mean curvature of the interface.  The proposed space-time strategy for the surfactant transport equation is similar to \cite{HanLarZah16} but here we rewrite the weak formulation so that the surfactant mass is discretely conserved without enforcing the mass with a Lagrange multiplier. As in \cite{HanLarZah16} and \cite{FrZa19} the space-time integrals in the weak formulation are approximated using quadrature rules, utilizing that the stabilization term defines the extension we need to have a robust and stable method. The proposed method decouples the representation and evolution of the interface from the discretization of the PDEs and although we in this work use the standard level set method \cite{OshFed01, Set01} we emphasize that other interface representation techniques can also be used. See for example \cite{Zah18} where a space-time cut finite element discretization is used together with an explicit representation of the interface using cubic splines.     

The paper is organized as follows. We start with some basic notation in Section \ref{sec:notation}. In Section \ref{sec:mathmod} we state the mathematical model and a weak formulation.  In Section \ref{sec:mesh_space} we introduce  the mesh and the finite element spaces we use in our finite element method. In Section \ref{sec:surfactant} we propose a space-time cut finite element method for a surface PDE modeling the surfactant transport on an evolving interface but here we assume the velocity field is known and explicitly given. We show numerical examples and convergence studies. We then propose a space-time CutFEM  for the coupled problem where the velocity field is given by the Navier--Stokes equations in Section \ref{sec: weakform_nummeth} and we show simulations in both two and three space dimensions in Section \ref{sec:numexpcoup}. Finally, we conclude in Section \ref{sec:con}.

\section{Notation}\label{sec:notation}
\begin{figure}[ht]
\begin{center}
\includegraphics[width=0.35\textwidth]{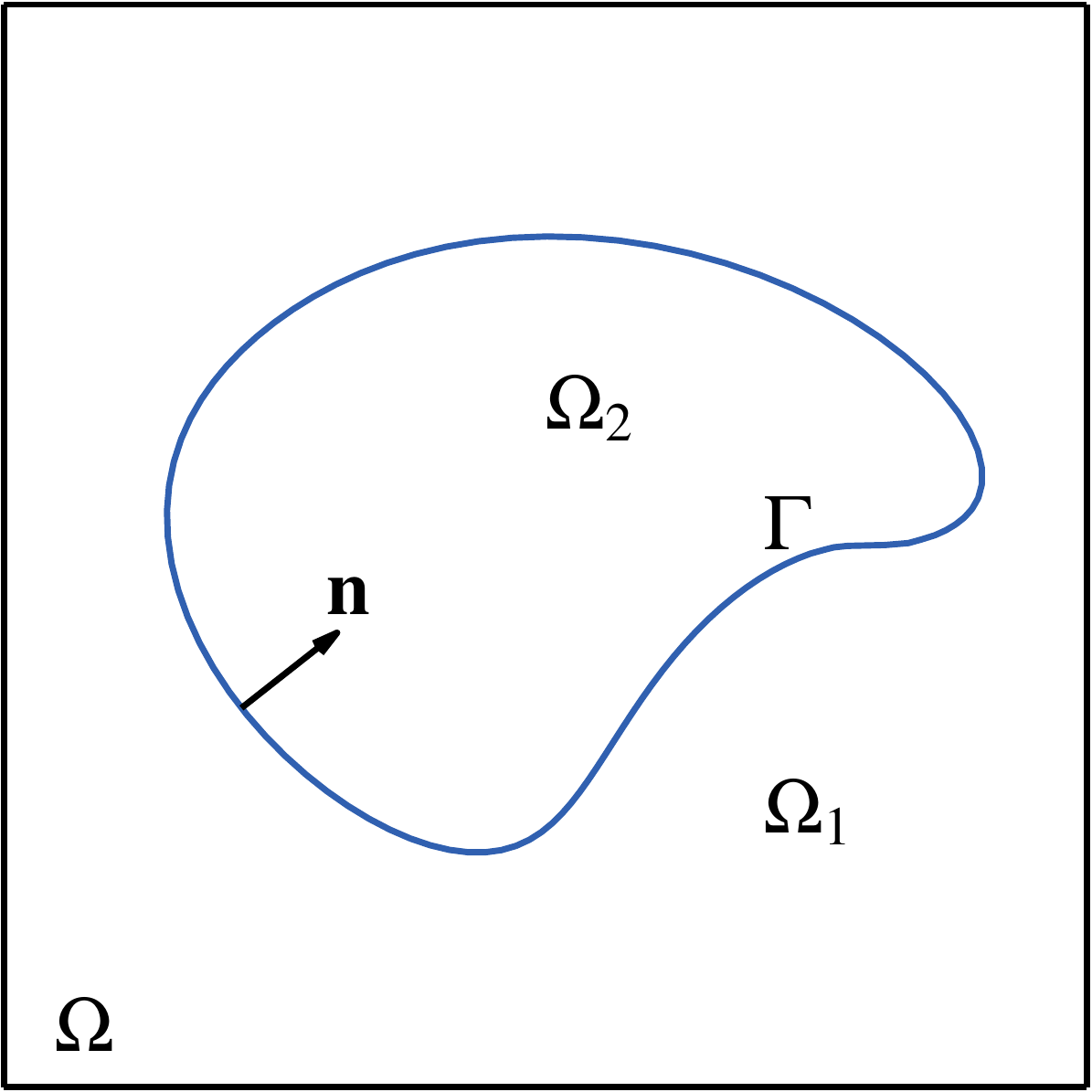}  \hspace{1cm}
\caption{Illustration of the domain $\Omega \in \R^2$ and the two subdomains $\Omega_1(t)$ and $\Omega_2(t)$ separated by the interface $\Gamma(t)$ at a fixed time instance $t \in [0,T]$.   \label{fig:domain}}
\end{center}
\end{figure}

\noindent Consider a bounded convex domain $\Omega$ in $\R^d$, $d=2,3$, with polygonal boundary $\dw$. We assume that at any time instance $t$ in the time interval $I= [0,T]$ ($T\in(0,\infty)$) the domain $\Omega$ contains two subdomains $\Omega_i(t) \subset \Omega$, $i=1,2$ that are separated by a sufficiently smooth closed $d-1$ dimensional internal interface $\Gamma(t) = \dw_1(t) \cap \dw_2(t)$. The time dependent interface is transported with a sufficiently smooth velocity field $\vel: I \times \Omega \rightarrow \R^d$. We assume that for all $t \in I$, the subdomain $\Omega_2(t)$ is enclosed by $\Gamma(t)$ and $\bar{\Omega} = \bar{\Omega}_1(t) \cup \bar{\Omega}_2(t)$ so $\dw \subset \dw_1$. Let $\bfn$ be the unit normal to $\partial \Omega_1$ which is outward directed with respect to $\Omega_1$. Thus, at the interface $\Gamma$, the normal $\bfn$ is pointing from $\Omega_1$ into $\Omega_2$. See Figure~\ref{fig:domain} for an illustration in two space dimensions ($d=2$).

For a scalar function $v:I \times \Omega \rightarrow \R$ we use the notation $v_i = v|_{\Omega_i(t)}$, $i=1,2$ and define the jump across the interface $\Gamma$ by $$\jump{v} = (v_1 - v_2)\restreinta{\Gamma}.$$ 
We denote by $\jump{v}_t$ the jump of a function $v$ in time at time instance $t$, i.e.,
\begin{equation}
\jump{v}_t = v^+(t,\bfx) - v^-(t,\bfx), \quad v^{\pm}(t,\bfx)=\lim_{\varepsilon \to 0^+}  v(t \pm \varepsilon,\bfx).
\end{equation} 
We do similarly for a vector field (i.e. $\bfv:I \times \Omega \rightarrow \R^d$) or a matrix-valued function but we denote those functions in bold.

The usual $\R^d$ gradient is denoted by $\nabla$ and $\nablas$ denotes the surface gradient on $\Gamma$ and is defined by $\nablas=P_{\Gamma} \nabla$, where $P_{\Gamma} = \bfI - \bfn \otimes \bfn$, $\mathbf{I} \in \R^{d\times d}$ is the identity matrix, and $\bfn$ is the unit normal to the surface $\Gamma$. The surface divergence of a vector field $\bfv(\bfx)$, is defined by 
\begin{equation}
\gradS \cdot \bfv=\text{tr}(\bfv \otimes \nabla_\Gamma) 
=\text{tr}(\bfv \otimes \nabla) - \bfn \cdot (\bfv \otimes \nabla )\cdot \bfn=\nabla \cdot \bfv- \bfn \cdot \nabla \bfv \cdot \bfn, \quad \bfx \in \Gamma.
\end{equation}
Note that for tangent vector fields $\bfv$ on $\Gamma$ we have the identity $\gradS \cdot \bfv= \nabla \cdot \bfv$, also note that $\bfv \otimes \nabla=(\nabla \otimes \bfv)^T$. The Laplace--Beltrami operator is denoted by $\Delta_\Gamma = \nabla_\Gamma \cdot \nabla_\Gamma$.

Recall the Sobolev spaces
\begin{equation}
L^2(U)=\enstq{v}{\int_{U} |v|^2 \diff \bfx=\| v \|_{L^2(U)}^2 < \infty},
\end{equation}
and
\begin{equation}
H^1(U)=\enstq{v}{v \in L^2(U), \grad v \in L^2(U)},
\end{equation}
where $U \subset \R^d$ and is in this paper $\Omega$, $\Omega_i(t)$, or $\Gamma(t)$ for $t\in I$. We will use the notation $\prodscal{v}{w}_U=\int_U v(\bfx) w(\bfx) \diff \bfx \ $ for the $L^2$-inner product on $U$ (similarly for inner products in $[L^2(U)]^d$). 
For $t\in I$ we write $v(t) \in H^1(\Omega_1(t) \cup \Omega_2(t))$ and we mean $v(t)=v(t,\cdot)=(v_1,v_2)$ where $v_i = v|_{\Omega_i(t)} \in H^1(\Omega_i(t))$. We use the notation
\begin{equation}
\prodscal{v}{w}_{\Omega_1 \cup \Omega_2}=\sum_{i=1}^2\prodscal{v_i}{w_i}_{\Omega_i}.
\end{equation}

Given the time interval $I$ we use the notation $\spacetimeS$ to denote the following space-time surface
\begin{equation}
  \spacetimeS=
  \bigcup_{t\in I} \{t\}\times \Gamma(t).
\end{equation}

\section{The Mathematical model}\label{sec:mathmod}
In this paper we consider two-phase flow systems with insoluble surfactants present.  The subdomains $\Omega_1$ and $\Omega_2$ are occupied by incompressible and immiscible fluids with densities $\rho_i \in \R_{>0}$ and viscosities $\mu_i \in \R_{>0}$, $i=1,2$. 
Let $\vel : I \times \Omega \rightarrow \R^d$ denote the fluid velocity, $p : I \times \Omega \rightarrow \R$ denote the pressure, and $w : \spacetimeS \rightarrow [0, w_\infty)$ denote the surfactant concentration where $w_\infty$ is the highest surfactant concentration on $\Gamma$. We assume the dynamics is governed by the incompressible Navier--Stokes equations coupled to a time-dependent convection-diffusion equation on the interface $\Gamma(t)$ that separates the two immiscible fluids:
\begin{alignat}{2}
\rho  \left( \partial_t \vel + (\vel \cdot \grad){\vel} \right)  
	- \div{\left( 2\mu \bfeps(\vel) - p \bfI \right)}  &= \bff  & &\quad \text{in} \ \tdomain,  \quad t \in (0,T] \label{eq : Navier-Stokes 1} \\
\div{\vel} & =  0   & &\quad \text{in} \ \tdomain, \label{eq : Incompressibility}\\
\jump{\vel} & =  0 & &\quad \text{on} \ \Gamma(t), \label{eq : Continuity velocity}\\
\jump{(2\mu \bfeps(\vel) - p \mathbf{I}) \vect{n}} & =  \sigma(w) \meancurv  + \gradS{\sigma(w)} & & \quad \text{on} \ \Gamma(t) , \label{eq : Surface tension 1}\\
\vel(t, \bfx) & =  \bfg(t,\bfx) &  &\quad \text{on} \ \partial \Omega \\
\vel(0, \bfx) & =  \vel^0(\bfx) &  &\quad \text{in} \  \Omega_1(0) \cup \Omega_2(0)  \\
\partial_t w + \vel \cdot \grad{w}  + (\divS{\vel}) w 
- \divS{(\diffSurf \gradS{w})} & = 0 & & \quad  \text{on}  \ \Gamma(t) \label{eq : Surfactant} \\
w(0,\bfx) & = w^0(\bfx) & & \quad \bfx \in \Gamma(0) \label{eq : Initial condition surfactant}
\end{alignat}
Here, $\partial_t = \frac{\partial}{\partial t}$, $\bfeps(\vel) = (\grad{\vel} + (\grad{\vel})^T)/2$ is the strain rate tensor, $\sigma$ is the surface tension which depends on the surfactant concentration, $\meancurv$ is the mean curvature vector defined by $\meancurv= \Delta_\Gamma \bfx_\Gamma $ with $\bfx_\Gamma: \Gamma \ni \bfx \rightarrow \bfx \in {\R}^d$ being the coordinate map or embedding of $\Gamma$ into $\R^d$,  $\diffSurf >0$ is the diffusion coefficient and is assumed to be constant, $\bff: I \times \Omega \rightarrow \R^d$ is a given external force (e.g. the gravitational force), and $\bfg$ is a given function such that $\int_{\partial \Omega} \bfg \cdot \bfn \diff s=0$.  We can also have other type of suitable boundary conditions on $\partial \Omega$.  The initial configuration, that is the interface $\Gamma(0)$ and hence the subdomains $\Omega_i(0)$, $i=1,2$, the initial data $\vel^0:\Omega \rightarrow \R^d$, and $w^0:\Gamma(0) \rightarrow [0,w_\infty)$ are also given. 
In this work we do not consider high Reynolds number flows. 

The evolution of the interface $\Gamma(t)$ is given by the following advection equation 
\begin{align}
\partial_t \lsf + \vel \cdot \grad{\lsf} = 0  \quad \textrm{in } \Omega, t \in I. \label{eq:Level set}
\end{align}
Here $\lsf(t,\bfx) : I \times \Omega \rightarrow \R$ is the signed distance function with positive sign in $\Omega_2(t)$, the subdomain enclosed by the interface, and the zero level set of this function represents the interface $\Gamma(t)$, i.e.,
\begin{equation}
  \Gamma(t)=\enstq{\bfx \in \Omega}{\lsf(t,\bfx)=0}.
\end{equation}
The initial condition $\lsf(0, \bfx) = \lsf_0(\bfx)$, $\bfx \in \Omega$ is given by the configuration of $\Gamma(0)$.

Finally, we use the following linear equation of state: 
\begin{equation}\label{eq:lineareqstate}
\sigma(w) = \sigma_0(1 - \beta w ),
\end{equation} 
where $\sigma_0 \in \R_{>0}$ and $\beta \in \R_{>0}$ are given parameters. Note that $\beta=0$ corresponds to the surfactant-free (clean) interface with $\sigma=\sigma_0$ being constant. Other models then this linear model describing the relation between the surface tension $\sigma:\R_{>0} \rightarrow \R_{>0}$ and the surfactant concentration $w$ can also be used, see~\cite{RaFeLi00}.  The numerical method we will present is not restricted to this linear relation and we discuss this in Section \ref{sec: implem_coupledp}.

\subsection{A weak formulation}\label{sec:weakform}
Let
\begin{align*}
\bfV = & \left[H^1(\tdomain) \right]^d, \\
Q =& \enstq{q \in L^2(\tdomain) }{\prodscal{\mu^{-1}q}{1}_{\tdomain} = 0},\\
W = & H^1(\spacetimeS).
\end{align*}
A weak formulation of the two-phase flow model~\eqref{eq : Navier-Stokes 1}-\eqref{eq : Initial condition surfactant} is: Find $\vel(t)=\vel(t,\cdot) \in \bfV$, $p(t)=p(t,\cdot) \in Q$, $w(t)=w(t,\cdot)  \in W$ such that for almost all $t\in I$ 
\begin{align}\label{eq:weakformfullp}
\prodscal{\rho (\partial_t \vel+(\vel \cdot \grad){\vel})}{\bfv}_{\tdomain } 
+a(t, \vel,\bfv)-b(t,\bfv,p)+&b(t,\vel,q)
+\prodscal{\sigma(w)\nabla_{\Gamma} \bfx_{\Gamma}}{\gradS \langle \bfv \rangle}_{\Gamma(t)}
\nonumber \\
&=l(t,v,q), \quad \forall v \in \bfV, q\in Q, \nonumber \\
\frac{\diff}{\diff t} \prodscal{w}{r}_{\Gamma(t)} 
 - \prodscal{w}{\partial_t r+\vel \cdot \grad{r}}_{\Gamma(t)} 
 +\prodscal{\diffSurf \gradS{w}}{\gradS{r}}_{\Gamma(t)}
 &=0,
  \quad  \forall r\in W, 
\end{align}
and $\vel(0,\bfx)  =  \vel^0(\bfx)$ in $\Omega_1(0) \cup \Omega_2(0)$ and $w(0,\bfx) = w^0(\bfx)$ on $\Gamma(0)$.
Here, the time derivative $\partial_t \vel$ is defined in a suitable weak sense and 
\begin{align}\label{eq:forma}
	a(t,\vel,\bfv) & = 
        \prodscal{2\mu\varepsilon(\vel)}{\varepsilon(\bfv)}_{\tdomain}  \nonumber \\
&	- \prodscal{\lbrace 2\mu \varepsilon(\vel)  \bfn \rbrace}{ \jump{\bfv}}_{\Gamma(t)}
	-  \prodscal{ \jump{\vel}}{\lbrace 2\mu \varepsilon(\bfv)  \bfn \rbrace}_{\Gamma(t)}
	+ \prodscal{\lambda_{\Gamma} \jump{\vel}}{\jump{\bfv}}_{\Gamma(t)}
	\nonumber\\
&	- \prodscal{2\mu \varepsilon(\vel)  \bfn}{\bfv}_{\dOmega}
	- \prodscal{\vel}{2\mu \varepsilon(\bfv)  \bfn}_{\dOmega}	
	+ \prodscal{\lambda_{\dOmega}  \vel}{\bfv}_{\dOmega},
\end{align}
\begin{align}\label{eq:formb}
	b(t,\bfv,q) & = \prodscal{\div \bfv}{q}_{\tdomain}  
	- \prodscal{\jump{\bfv \cdot \bfn}}{ \{ q \}}_{\Gamma(t)}	
	- \prodscal{\bfv \cdot \bfn}{q}_{\partial \Omega}, 
\end{align}
\begin{align}\label{eq:forml}
	l(t,\bfv,q) & = \prodscal{\mathbf{f}}{\bfv}_{\tdomain} 
	- \prodscal{\bfg }{ 2\mu \varepsilon(\bfv)  \bfn}_{\dOmega} 
        +\prodscal{\lambda_{\dOmega} \vect{g} }{ \bfv}_{\dOmega} 
        - \prodscal{\bfg \cdot \bfn }{q}_{\dOmega}, 
\end{align}
with $\lambda_{\Gamma} \in L^{\infty}(\Gamma)$, $\lambda_{\partial \Omega} \in L^{\infty}(\partial \Omega)$,  and 
\begin{equation}\label{eq:averagop}
\{f \}=\w_1f_1+\w_2f_2, \qquad  \langle f \rangle=\w_2f_1+\w_1f_2,
\end{equation}
where $f_i=f\restreinta{\Omega_i}$ and the weights $\w_1$ and $\w_2$ are positive real numbers satisfying $\w_1+\w_2=1$.

To derive the weak form we use the following variant of Reynolds' transport theorem: For $w,r \in W$ we have 
\begin{align}\label{eq:Reynoldwr}
\frac{\diff}{\diff t} \int_{\Gamma(t)} wr  \diff s
	& = \int_{\Gamma(t)} \partial_t(wr) + \vel \cdot \grad(wr)+wr\divS{\vel} \ \diff s \nonumber \\
& = \int_{\Gamma(t)}  (\partial_tw  + \vel \cdot \grad w+ w\divS{\vel })r+
w (\partial_t r + \vel \cdot \grad r) \ \diff s.
\end{align}
See Lemma 5.2 of \cite{DzEl13}. Assume that $(\vel,p,w)$ is a sufficiently smooth solution to~\eqref{eq : Navier-Stokes 1}-\eqref{eq : Initial condition surfactant} and so that Reynolds' transport theorem~\eqref{eq:Reynoldwr} holds. Then, it follows that  
\begin{equation}
  \frac{\diff}{\diff t} \prodscal{w}{r}_{\Gamma(t)} 
 - \prodscal{w}{\partial_t r+\vel \cdot \grad{r}}_{\Gamma(t)} 
 +\prodscal{\diffSurf \gradS{w}}{\gradS{r}}_{\Gamma(t)}
 =0, \quad \forall r\in W.
 \label{weakformsurf}
\end{equation}
Following the derivation in~\cite{FrZa19} we also have that for $t\in I$
\begin{align} \label{eq:weakformNS}
\prodscal{\rho (\partial_t \vel+(\vel \cdot \grad){\vel})}{\bfv}_{\tdomain } 
+a(t, \vel,\bfv)-b(t,\bfv,p)+&b(t,\vel,q)
-\prodscal{\sigma(w) \meancurv+\gradS{\sigma(w)}}{\langle \bfv \rangle}_{\Gamma(t)} 
\nonumber \\
&=l(t,v,q), \quad \forall v\in \bfV, q\in Q.
\end{align}
Recalling the definition of the mean curvature vector $\meancurv$, i.e., 
\begin{equation}
\meancurv= \Delta_\Gamma \bfx_\Gamma, \quad  \bfx_\Gamma: \Gamma \ni \bfx \rightarrow \bfx \in {\R}^d,
\end{equation}
assuming $\Gamma(t)$ is sufficiently smooth, and using integration by parts we have 
\begin{equation}\label{eq:surftenforce}
\prodscal{\sigma(w) \meancurv  + \gradS{\sigma(w)}} {\bfv }_{\Gamma(t)}=\prodscal{\gradS \cdot (\sigma(w) \nabla_{\Gamma} \bfx_{\Gamma} )}{ \bfv }_{\Gamma(t)}=-\prodscal{\sigma(w)\nabla_{\Gamma} \bfx_{\Gamma}}{\gradS \bfv}_{\Gamma(t)},
\end{equation} 
see also Chapter 7.6.1 in \cite{GrRe11}. Using \eqref{eq:surftenforce} in \eqref{eq:weakformNS} we arrive at our weak formulation. 

\section{The mesh and finite element spaces}\label{sec:mesh_space}
\begin{figure}[ht]
\centering
\includegraphics[scale=0.35]{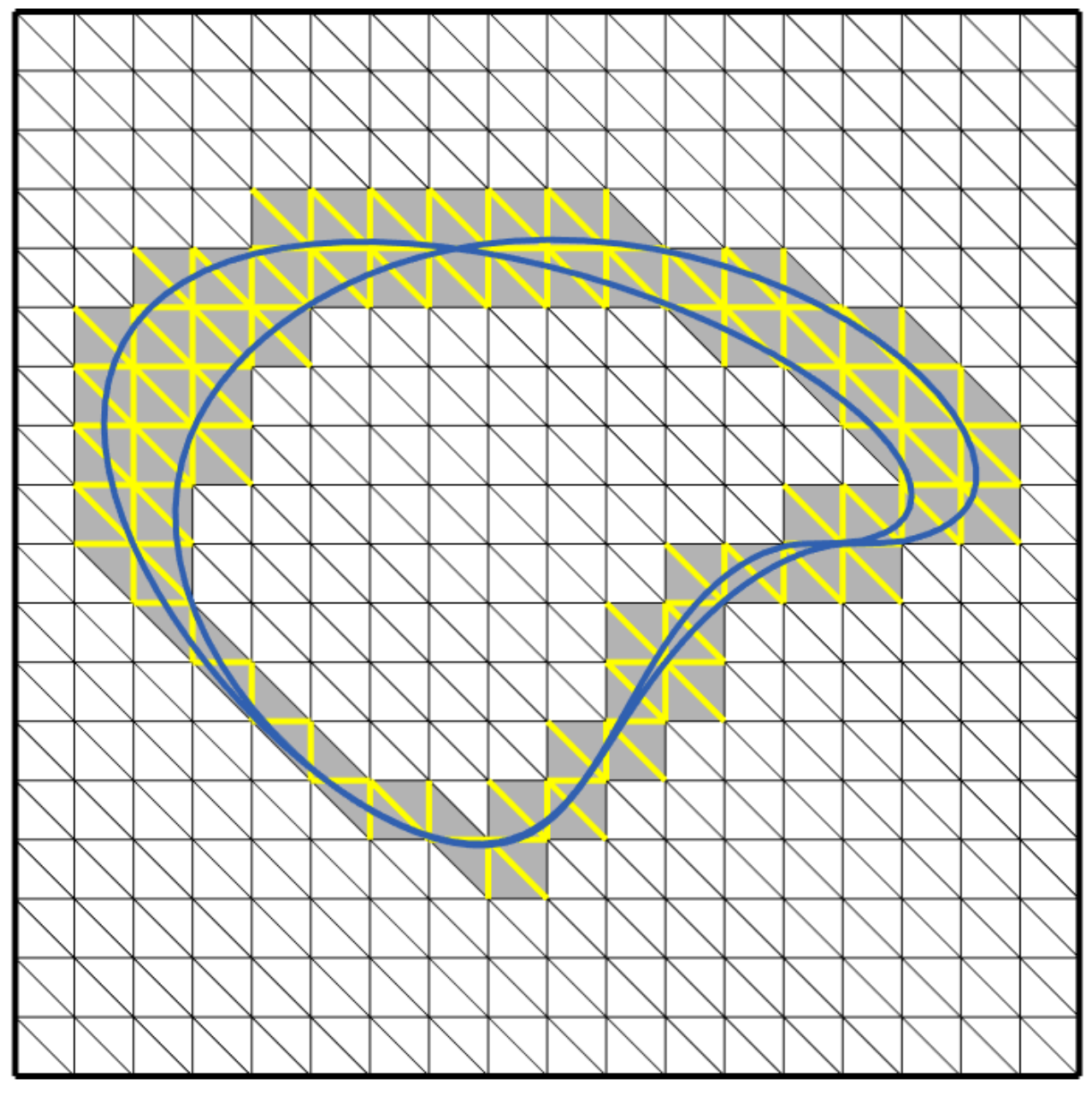}\hspace{0.7cm}
\includegraphics[scale=0.35]{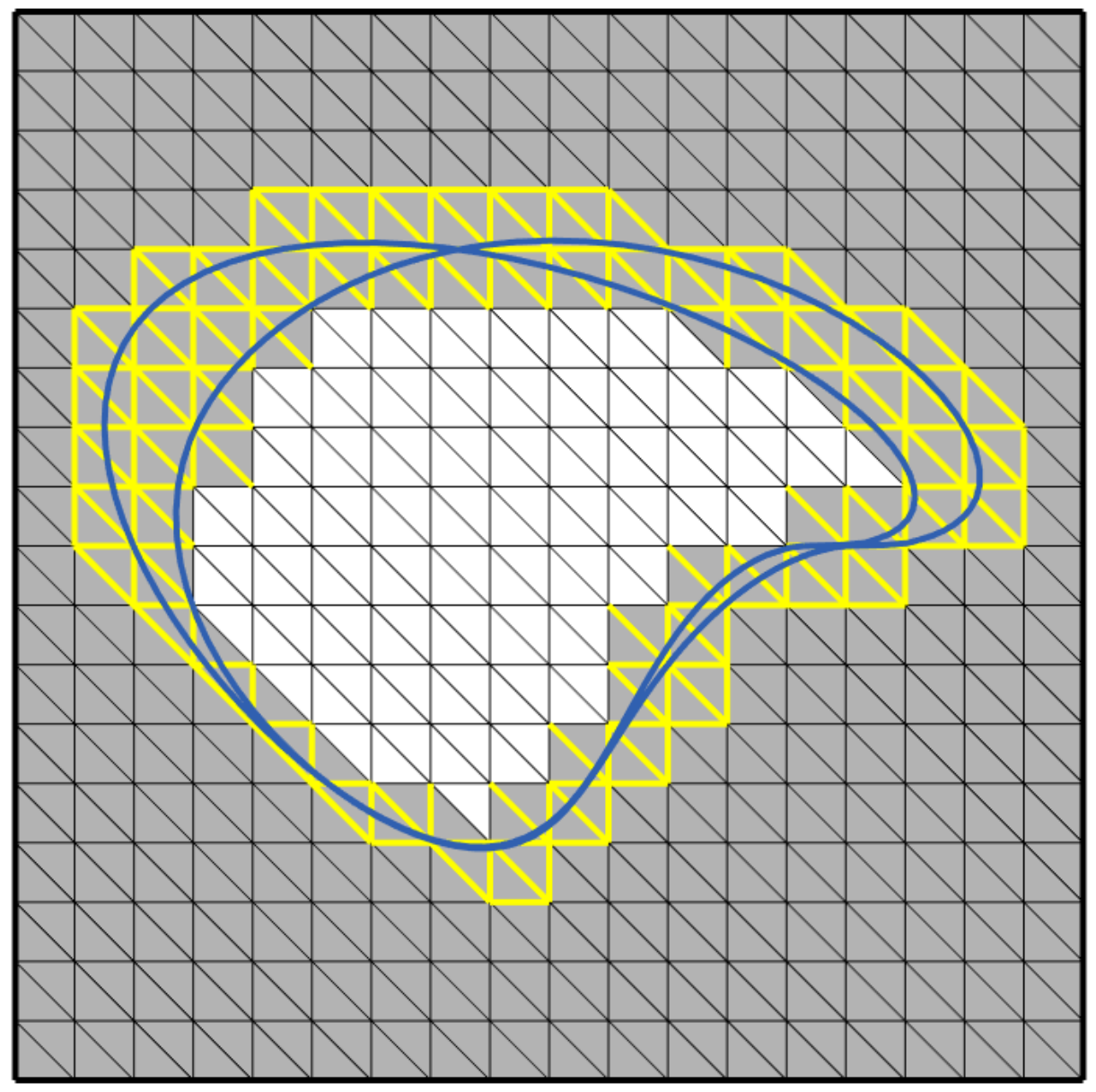}\hspace{0.7cm}
\includegraphics[scale=0.35]{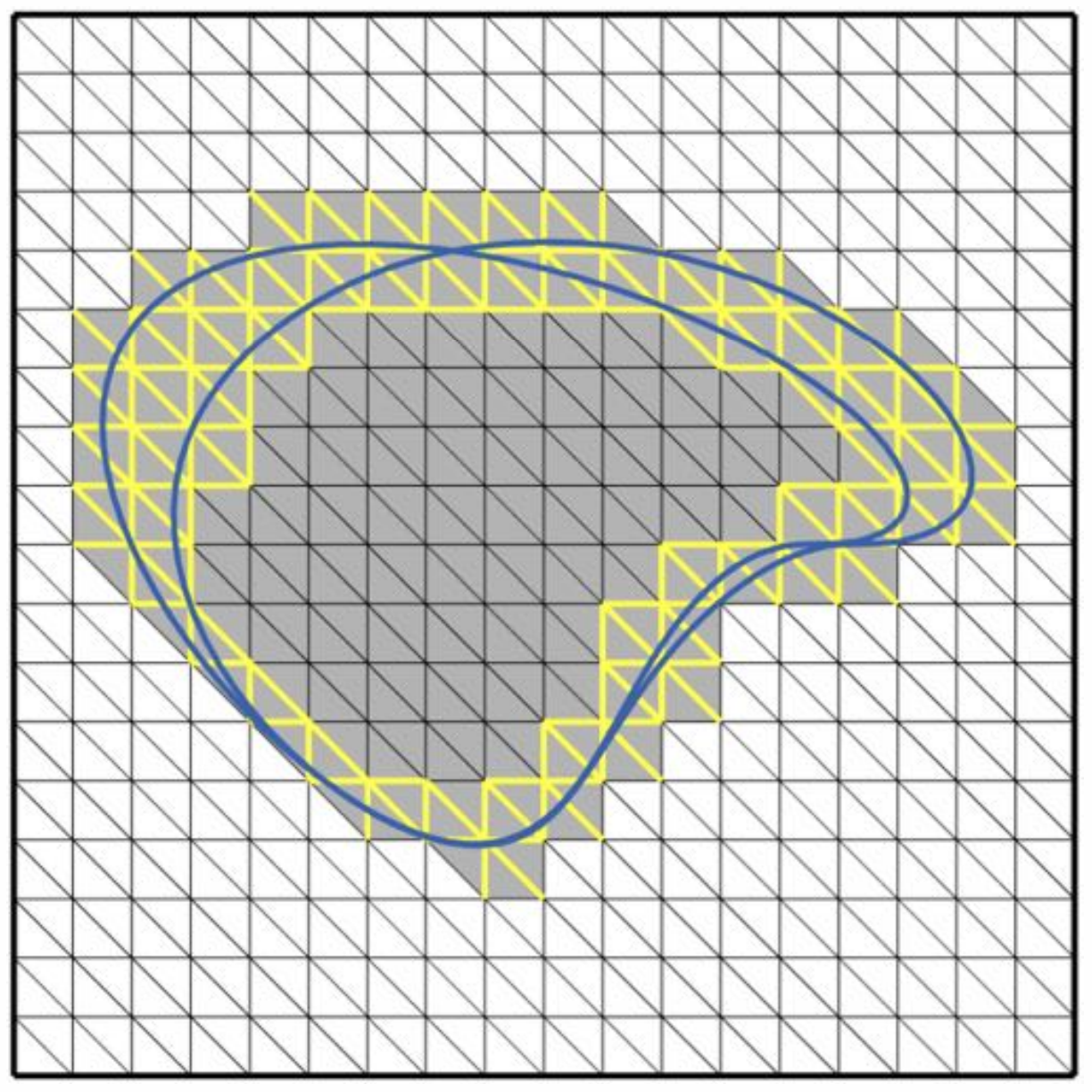}
\caption{Illustration of the sets introduced in Section~\ref{sec:mesh_space}. The two blue curves show the position of the interface $\Gamma(t)$ at the endpoints $t=t_{n}$ and $t=t_{n+1}$ of the time interval $I_n$.  The shaded domain shows $\mcN_{h,i}^n$ and edges in $\mcF_{h,i}^n$ are marked with yellow thick lines. Left: $i=0$. Middle: $i=1$. Right: $i=2$. 
\label{fig:illustmesh}}
\end{figure}

Here, we define the mesh and the finite element spaces we will need later when we define the finite element method. Let $\{\mesh \}$ be a regular family of simplicial meshes of $\Omega$ with $h\in(0,h_0]$. We denote by $h_K$ the diameter of element $K\in \mesh$ and we let $h$ be the piecewise constant function that is equal to $h_K$ on element $K$.  The mesh $\mesh$, which we will refer to as the background mesh, conforms to the fixed polygonal boundary $\dw$ but does not need to conform to the interface and hence at any time instance $t \in I$ the interface $\Gamma(t)$ may cut through the mesh arbitrarily, see Figure~\ref{fig:illustmesh}. However, the background meshes we use are not completely independent of the interface as we typically use graded meshes with more elements in interface regions. 
For each time instance $t\in I$ we denote by $\mcT_{h,i}(t)$ the collection of elements that exhibit a nonempty intersection with $\Gamma(t)$ for $i=0$ and with $\Omega_i(t)$, for $i=1,2$, i.e.,  
\begin{align}\label{eq:khi}
  \mcT_{h,0}(t) &=\enstq{K \in \mcT_{h} }{K \cap \Gamma(t) \neq \emptyset},
\quad
	\mcT_{h,i}(t) = \enstq{K \in \mcT_h}{ K \cap \Omega_i(t) \neq \emptyset}, \quad i=1,2, \ t \in I.
\end{align}

Let $0=t_0 < t_1 < \dots < t_N = T$ be a partition of $I$ into time intervals $I_n = (t_{n}, t_{n+1}]$ of length $\Delta t_n = t_{n+1} - t_{n}$,  $n = 0,1,\dots, N-1$. 
For each time interval $I_n$ we define active meshes $\mcT_{h,i}^n$, $i=0,1,2$ and corresponding domains $\mcN_{h,i}^n$. The active mesh $\mcT_{h,i}^n$ contains those elements in $\mesh$ that cover the domain $\mcN_{h,i}^n$ where 
\begin{align}
	\mcN_{h,i}^n &= \bigcup_{t \in I_n} \bigcup_{K \in \mcT_{h,i}(t)} K, \quad i=0,1,2. 
\end{align}
Denote by $\mathcal{F}_{h}^n$ the set of faces in the active mesh $\mcT_{h,0}^n$. Faces in $\mathcal{F}_{h}^n$ that are shared by two elements in the active mesh $\mcT_{h,i}^n$ constitute a set denoted by $\mathcal{F}_{h,i}^n$.  For an illustration, in two space dimensions, of the sets introduced in this section see Figure~\ref{fig:illustmesh}.

Let $P_{k}(D)$ be the space of polynomials in $D \subset \R^d$, $d=1,2,3$, of degree less or equal to $k\geq 0$. Denote by $\Bspace^m$ the space of continuous piecewise polynomials of degree less than or equal to $m\geq 1$ defined on the background mesh $\mesh$, i.e., 
\begin{equation}
  \Bspace^m=\enstq{v \in C^0(\bar{\Omega})}{v \restreinta{K} \in P_m(K), \ \forall K \in \mesh}.
\end{equation}
On the space-time slab $I_n \times \mcN_{h,0}^n$  we define the space
\begin{align}\label{eq:spaceVhi}
	W_{h,m}^{n,k} & = P_k(I_n)  \otimes \Bspace^m\restreinta{\mcN_{h,0}^n}, \quad k\geq 0, m\geq 1. 	
\end{align}
In this paper, we will propose a finite element method for the discretization of the surfactant transport equation using linear polynomials in both space and time, $m=k=1$. We use the notation $W_h^n=W_{h,1}^{n,1}$. Note that a function  $w_h \in W_{h}^{n}$ can be represented in the following form
\begin{align}\label{eq:repofwh}
	w_h(t, \bfx) =  \sum_{j=0}^1 w_{h,j} \left(\frac{t-t_{n}}{\Delta t_n} \right)^j, \quad t \in I_ n, \ \bfx \in \mcN_{h,0}^n
\end{align}
and for each j    
\begin{equation} \label{eq:sumk}
	w_{h,j} = \sum_{i} \textbf{$\xi$}_{i,j} \phi_i(\bfx)\restreinta{\mcN_{h,0}^n}, 
\end{equation}
where $\textbf{$\xi$}_{i,j} \in \R$ are coefficients and $\phi_i(\bfx)$ is the standard nodal basis function associated with mesh vertex $i$ in the finite element space $\Bspace^1$. 

\begin{figure}
\centering
	\includegraphics{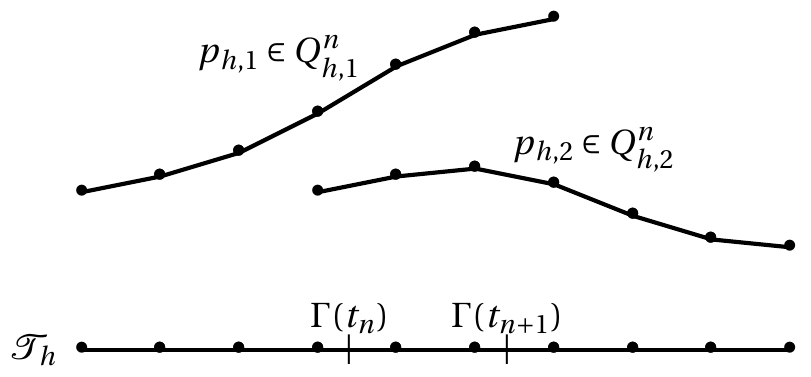}
	 \caption{Illustration of a function $p_h(t,\bfx)=(p_{h,1},p_{h,2}) \in Q^n_h$ at a time instance $t\in I_n$ in a one space dimensional model. The interface evolves from $\Gamma(t_{n})$ to $\Gamma(t_{n+1})$ during the time interval $I_n$. The domain to the left of the interface $\Gamma(t)$ is $\Omega_1(t)$ and the domain to the right is $\Omega_2(t)$. \label{fig:illustcutfemspace}} 
\end{figure}
For the discretization of the incompressible Navier--Stokes equations we start from the LBB stable Taylor-Hood element P2-P1 pair, ($\bfV_h, Q_h$), on the background mesh $\mesh$, i.e.,
\begin{align}
  \bfV_h &=  \left[\Bspace^2 \right]^d \\
  Q_h&=\Bspace^1 \cap Q=\enstq{q_h \in C^0(\bar{\Omega})}{q_h\restreinta{K} \in P_1(K), \ \forall K \in \mesh, \ \prodscal{\mu^{-1}q_h}{1}_{\Omega} = 0}.
  \end{align}
On each of the space-time slabs $I_n \times \mcN_{h,i}^n$, $i=1,2$, we define for $k\geq 0$ the spaces
\begin{align}
	\bfV_{h,i}^{n,k} & = P_k(I_n)  \otimes \bfV_{h}\restreinta{\mathcal{N}_{h,i}^n},  \quad i=1,2, 	\\
	Q_{h,i}^{n,k} & = P_k(I_n)  \otimes Q_{h}\restreinta{\mcN_{h,i}^n},  \quad i=1,2. \label{eq:spaceQhi}
\end{align}
Finally, let
\begin{align*}
 \bfV^{n,k}_{h} &=  
 \enstq{\bfv_h = (\bfv_{h,1},\bfv_{h,2})}{\bfv_{h,i} \in \bfV^{n,k}_{h,i}, \ i=1,2}, \\
  Q^{n,k}_h &=  
 \enstq{q_h = (q_{h,1},q_{h,2})}{q_{h,i} \in Q^{n,k}_{h,i}, \ i=1,2}.
\end{align*}
For $k=1$ we denote these spaces by $\bfV^n_h=\bfV^{n,1}_h$ and  $Q^n_h=Q^{n,1}_h$. Note that functions in $\bfV^n_h$ and $Q^{n}_h$ are double valued on elements in $\mcT_{h,0}^n$ since those elements exist in both active meshes $\mcT_{h,1}^n$ and $\mcT_{h,2}^n$. For an illustration, in one space dimension, of a function $p_h \in Q^n_h$ at a time instance $t\in I_n$ see Figure~\ref{fig:illustcutfemspace}.

Before we formulate the finite element method for the two-phase flow model we propose and study a space-time cut finite element method for the surface PDE modeling surfactant transport.

\section{The surfactant transport equation} \label{sec:surfactant}
In this section we assume that the velocity field $\vel$ is known and explicitly given and we consider the following convection-diffusion problem:
\begin{alignat}{2}
  \partial_t w + \vel \cdot \grad{w}  + (\divS{\vel}) w 
- \divS{( \diffSurf \gradS{w})} &  = f & & \quad \text{on }  \Gamma(t), \quad t\in (0,T],   \label{eq:convection_diffusion_surfactant}\\
w(0,\bfx) & = w^0 & & \quad \bfx \in \Gamma(0).\label{eq:initsurf}
\end{alignat}
Following the derivation in Section \ref{sec:weakform} we have for $w$ being a sufficiently smooth solution to the differential equation \eqref{eq:convection_diffusion_surfactant}-\eqref{eq:initsurf}, $r\in W$, and $t\in I$ that 
\begin{equation}
  \frac{\diff}{\diff t} \prodscal{w}{r}_{\Gamma(t)} 
 - \prodscal{w}{\partial_t r+\vel \cdot \grad{r}}_{\Gamma(t)} 
 +\prodscal{\diffSurf \gradS{w}}{\gradS{r}}_{\Gamma(t)}
 =\prodscal{f}{r}_{\Gamma(t)},
 \label{weak formulation 2}
\end{equation}
Assuming $w^-(t_n,\bfx)$ is given, integrating in the time interval $I_n$, and using that $w^+(t_n,\bfx)=w^-(t_n,\bfx)$, we have 
 \begin{align}\label{eq:weakformintime}
& \int_{I_n}  \prodscal{f}{r}_{\Gamma(t)} \diff t=
\int_{I_n} \left[
 \frac{\diff}{\diff t} \prodscal{w}{r}_{\Gamma(t)} 
 - \prodscal{w}{\partial_t r+\vel \cdot \grad{r}}_{\Gamma(t)} 
 +\prodscal{\diffSurf \gradS{w}}{\gradS{r}}_{\Gamma(t)}
  \right] \diff t
\nonumber \\
 &
  = 
  \prodscal{w}{r}_{\Gamma(t_{n+1})} 
 - \prodscal{w^-(t_n,\bfx)}{r}_{\Gamma(t_{n})} 
 +\int_{I_n} \left[ 
 -  \prodscal{w}{\partial_t r+\vel \cdot \grad{r}}_{\Gamma(t)} 
 +\prodscal{\diffSurf \gradS{w}}{\gradS{r}}_{\Gamma(t)}
\right] \diff t.
\end{align}
In particular, note that with $r = 1$ in \eqref{eq:weakformintime} we get 
  \begin{align}\label{eq:conservcont}
   \int_{\Gamma(t_{n+1})} w(t_{n+1}, \bfx) \diff s - \int_{\Gamma(t_{n})} w^-(t_n, \bfx) \diff s  - \int_{I_n}  \int_{\Gamma(t)} f \diff s \diff t=0.
\end{align}

\subsection{The Finite element method}\label{sec:surfactant_FEMt}
We are now ready to present a discretization of the convection-diffusion equation on the interface. Assume the active mesh $\mcT_{h,0}^n$ and the solution $w_h^-(t_{n},\bfx)$ from the previous space-time slab are known: Find $w_h \in W_{h}^n$  such that
\begin{align}\label{eq:weakformsurf}
	A(\vel, w_h,r_h)  + S(w_h,r_h)= \int_{I_n}  \prodscal{f}{r_h}_{\Gamma(t)} \diff t +\prodscal{w_h^-}{r_h}_{\Gamma(t_{n})}, \quad \forall r_h \in W_h^n,
\end{align}
where 
\begin{align}
A(\vel, w_h,r_h)& =  \prodscal{w_h}{r_h}_{\Gamma(t_{n+1})} 
 -\int_{I_n}  \prodscal{w_h}{\partial_t r_h+\vel \cdot \grad{r_h}}_{\Gamma(t)} 
 \diff t
+\int_{I_n} \prodscal{\diffSurf \gradS{w_h}}{\gradS{r_h}}_{\Gamma(t)}  \diff t
  \label{weakformulationFE2}
\end{align}
and
\begin{align}\label{eq:stabsurf}
	S(w_h,r_h) = \int_{I_n} s_w(t,w_h,r_h) \diff t. 
\end{align}
The last term $S$ is a stabilization term which is added to ensure that the method is stable and robust independent of how the interface is positioned relative the mesh. There are several choices for $s_w$, see e.g. \cite{BurHanLar15, LarZah19, OlsReuXu14-b, BurHanLarMasZah20} and references therein. In this work we use 
\begin{align}
	s_w(t,w_h, r_h) = 
		 \sum_{j=1}^m \left(c_{F,j}
		 h^{2j-1} 
		 \prodscal{\jump{ D^j_{\bfn_F} w_h}_F}{\jump{ D^j_{\bfn_F} r_h}_F}_{\mcF_{h,0}^n} 
		 +
		c_{\Gamma,j} h^{2j-1} 
		\prodscal{D^j_{\bfn_{{h}}} w_h}{ D^j_{\bfn_{{h}}} r_h}_{\Gamma_{h}(t)} \right), \label{eq:surfstab}
\end{align}
where $c_{F,j}$, and $c_{\Gamma,j}$ are positive constants. This stabilization was analyzed in connection with a CutFEM for surface PDEs on stationary surfaces in \cite{LarZah19}. For the space-time CutFEM it is important to define the set $\mathcal{F}_{h,0}^n$, on which ghost penalty stabilization is applied, appropriately (see Figure \ref{fig:illustmesh}). In particular note that in a space-time slab the set $\mathcal{F}_{h,0}^n$ does not depend on time. The proposed stabilization term yields an extension and therefore allows us to propose a cut finite element scheme based on directly approximating the space-time integrals in the weak formulation using quadrature rules, see the next section. The proposed stabilization can be used also when the polynomial degree $m>1$. For $m=1$ an alternative is to extend the stabilization proposed in \cite{BurHanLar15} to the space-time CutFEM, thus only having ghost penalty terms evaluated on the faces in $\mathcal{F}_{h,0}^n$, i.e.,   
\begin{align*}
s_w(t,w_h, r_h) = c_{F,j}   \prodscal{\jump{ D^j_{\bfn_F} w_h}_F}{\jump{ D^j_{\bfn_F} r_h}_F}_{\mcF_{h,0}}.
\end{align*}
If the problem is convection dominated other stabilization terms such as streamline diffusion stabilization can be added, see e.g. \cite{OlsReuXu14-b, BurHanLarMasZah20}.

\subsubsection{Space-time CutFEM with quadrature in time and implementation}\label{sec:implementationsurf}
We propose a weak form where integrals in time are replaced with a $n_q$-point quadrature rule with weights $\alpha_q^n$ and points $t_q^n$, $q=1, \cdots, n_q$. Thus, we propose the following weak formulation: Given $w_h^-(t_{n},\bfx)$, find $w_h\in W_h^n$ such that 
\begin{align}\label{eq:weakformsurfdis}
	A^n_h(\vel, w_h,r_h)  + \sum_{q=1}^{n_q} \alpha_q^n s_w(t_q^n, w_h,r_h)=\sum_{q=1}^{n_q} \alpha_q^n  \prodscal{f}{r_h}_{\Gamma(t_q^n)}  +\prodscal{w_h^-}{r_h}_{\Gamma(t_{n})}, \quad \forall r_h \in W_h^n,
\end{align}
where 
\begin{align}
A_{h}^n(\vel, w_h,r_h)& =  \prodscal{w_h}{r_h}_{\Gamma(t_{n+1})} 
 -\sum_{q=1}^{n_q} \alpha_q^n  \prodscal{w_h}{\partial_t r_h+\vel \cdot \grad{r_h}}_{\Gamma(t_q^n)} 
+\sum_{q=1}^{n_q} \alpha_q^n  \prodscal{\diffSurf \gradS{w_h}}{\gradS{r_h}}_{\Gamma(t_q^n)}.  
  \label{weakformulationFE2dis}
\end{align}
Recall the representation of functions in $W_h^n$, see~\eqref{eq:repofwh}, and note that for example for the second term in $A_h^n$ we have  
\begin{equation}
\sum_{q=1}^{n_q} \alpha_q^n (w_h, \partial_t r_h)_{\Gamma(t_q^n)} =
\sum_{q=1}^{n_q} \alpha_q^n\frac{1}{\Delta t_n}(w_{h,0},r_{h,1})_{\Gamma(t_q^n)} +
\sum_{q=1}^{n_q} \alpha_q^n \frac{t_q^n-t_{n}}{\Delta t_n^2}(w_{h,1},r_{h,1})_{\Gamma(t_q^n)}.
\end{equation}
Furthermore, since the set $\mathcal{F}_{h,0}^n$, see definition of $s_w$ in \eqref{eq:surfstab}, is on each space-time slab independent of time and the trial and the test functions, ($w_h$ and $r_h$) which are time dependent are polynomials, the quadrature rule can be chosen so that the integration of the face stabilization term is exact, i.e., 
\begin{align}
  \int_{I_n} \prodscal{\jump{ D^j_{\bfn_F} w_h(t,\bfx)}_F}{\jump{ D^j_{\bfn_F} r_h(t,\bfx)}_F}_{\mcF_{h,0}^n} dt  = \sum_{q=1}^{n_q} \alpha_q^n 
		 \prodscal{\jump{ D^j_{\bfn_F} w_h(t_q^n,\bfx)}_F}{\jump{ D^j_{\bfn_F} r_h(t_q^n,\bfx)}_F}_{\mcF_{h,0}^n}.
\end{align}
In the numerical examples we use Simpson's quadrature rule. Hence the quadrature points are $t_1^n=t_{n}$,  $t_3^n=t_{n+1}$, $t_2^n=\frac{t_{n}+t_{n+1}}{2}$, and the quadrature weights are $\alpha_1^n=\alpha_3^n=\frac{\Delta t_n}{6}$, $\alpha_2^n=\frac{4\Delta t_n}{6}$.

Note that taking $r_h=1$ in the proposed weak formulation \eqref{eq:weakformsurfdis} yields the following conservation property
\begin{align}\label{eq:conserverror}
 \prodscal{w_h}{r_h}_{\Gamma(t_{n+1})} - \prodscal{w_h^-}{r_h}_{\Gamma(t_{n})}-\sum_{q=1}^{n_q} \alpha_q^n \prodscal{f}{r_h}_{\Gamma(t_q^n)} =0.
 \end{align}

To implement the proposed method one needs to know the position of the interface at the quadrature points, $\Gamma(t_q^n)$, and also the set $\mathcal{F}_{h,0}^n$. Often we do not have the exact interface position. In this paper we use a piecewise planar approximation of the interface. Let $T^n=\{t_n, \{t_q^n\} , t_{n+1}\}$. For each $t\in T^n$, we either directly make an approximation of the level set function in the space $\Bspace^1$ if $\lsf$ is known explicitly at all time or we discretize the advection equation \eqref{eq:Level set} using the Crank-Nicolson scheme in time and continuous FEM in space with linear Lagrange elements and streamline diffusion stabilization, following~\cite{FrZa19}. By looking for sign change in each element $K\in \mesh$, i.e., checking the sign of this approximated level set function, $\lsf_h(t,\bfx)$, $\bfx \in \Omega$, we find those elements in the mesh that are cut by the interface and thus belong to the set $\mcT_{h,0}(t)$ for each $t\in T^n$(see \eqref{eq:khi}). We compute a piecewise planar approximation $\Gamma_h(t)$ at each time instance $t \in T^n$ by computing the intersection of the zero level set of $\lsf_h(t,\bfx)$ with each cut element $K\in \mcT_{h,0}(t)$.    

We also need the set of faces $\mathcal{F}_{h,0}^n$ where the stabilization is active. This set consist of all interior faces (faces with two neighbors) in the active mesh $\mcT_{h,0}^n$.  Hence we need to know if an element $K\in \mesh$ is in the active mesh $\mcT_{h,0}^n$ or not. We say that an element $K\in \mesh$ is in the active mesh if it is cut by the interface at some time instance $t \in T^n$ or if there are two time instances $t_k \in T^n$ and $t_l \in T^n$ where $K \cap \Omega_1(t_k) \neq \emptyset$ but $K \cap \Omega_1(t_l) = \emptyset$  (i.e. $K \cap \Omega_2(t_l) \neq \emptyset$ ).
Using our approximation of the signed distance function, $\lsf_h(t,\bfx)$,  we say that an  element $K\in \mesh$ belongs to $\mcT_{h,0}^n$ if it is in $\mcT_{h,0}(t)$ at some time instance $t \in T^n$ or the sign of $\lsf_h(t_k,K)$ is different than $\lsf_h(t_l,K)$ for two time instances $t_k, t_l \in T^n$. 

\begin{remark}
Note that if Reynolds' transport theorem \eqref{eq:Reynoldwr} holds, formulation \eqref{weak formulation 2} is equivalent to
  \begin{align}
    \prodscal{\partial_t  w+\vel \cdot \grad w+w \divS{\vel} }{r}_{\Gamma(t)} +\prodscal{\diffSurf \gradS{w}}{\gradS{r}}_{\Gamma(t)}=\prodscal{f}{r}_{\Gamma(t)},
	\label{weak formulation 1}
 \end{align}
and from \eqref{eq:weakformintime} we have
\begin{align}\label{eq:weakformintime2}
&\int_{I_n} \left[
 \frac{\diff}{\diff t} \prodscal{w}{r}_{\Gamma(t)} 
 - \prodscal{w}{\partial_t r+\vel \cdot \grad{r}}_{\Gamma(t)} 
 +\prodscal{\diffSurf \gradS{w}}{\gradS{r}}_{\Gamma(t)}
  \right] \diff t
\\  
 & 
  = \prodscal{\jump{w}_{t_n}}{r}_{\Gamma(t_n)} +\int_{I_n} \left[ \prodscal{\partial_t  w+\vel \cdot \grad w+w \divS{\vel}}{r}_{\Gamma(t)} +\prodscal{\diffSurf \gradS{w}}{\gradS{r}}_{\Gamma(t)}  \right] \diff t.
\end{align}
Hence, one may instead of the bilinear form $A_{h}^n$ in \eqref{weakformulationFE2dis} use
\begin{align}
  A_{h}^n& = \prodscal{w_h}{r_h}_{\Gamma(t_{n})} 
+ \sum_{q=1}^{n_q} \alpha_q^n \left[ \prodscal{\partial_t  w_h+\vel \cdot \grad w_h+w_h \divS{\vel}}{r_h}_{\Gamma(t_q^n)} +\prodscal{\diffSurf \gradS{w_h}}{\gradS{r_h}}_{\Gamma(t_q^n)}  \right].
  \label{weakformulationFE1dis}
\end{align}
This form was used in~\cite{HLZ16ST} together with a Lagrange multiplier enforcing the correct surfactant mass. We will see in the numerical examples that the two forms \eqref{weakformulationFE2dis} and \eqref{weakformulationFE1dis} may not give  equivalent discretizations since \eqref{eq:Reynoldwr} does not hold in our discrete setting. However, the proposed discretization which is based on \eqref{weakformulationFE2dis} results in discrete conservation of the surfactant mass without using a Lagrange multiplier, see equation \eqref{eq:conserverror}.
\end{remark}

\begin{remark}\label{remspacetimeFEM}
For the proposed method and its implementation the stabilization term \eqref{eq:stabsurf} is vital as it defines an extension of the approximate solution to the entire active space-time domain and hence allows to directly approximate the space-time integrals using quadrature rules. When the integrals in time are approximated by a quadrature rule, geometric computations such as the construction of the interface are done only at the quadrature points in time and space integrals are then computed at each quadrature point in time. We note that other space-time unfitted finite element methods have also been proposed see e.g.~\cite{CL15, OlRe14, OlReXu14a}. Those formulations are based on reformulating space-time integrals in the weak formulation to surface integrals over the space-time manifold in ${\R}^{d+1}$ (when $\Gamma$ is a surface in ${\R}^d$). In our formulation we never approximate surface integrals in ${\R}^{d+1}$. The implementation of the proposed space-time cut finite element method is straightforward starting from an implementation of CutFEM for stationary domains but relies on the stabilization terms and also on a method for representing and evolving the interface.
\end{remark}

\subsection{Numerical examples}
We consider one example in two space dimensions and two examples in three space dimensions and validate that the proposed unfitted finite element method is accurate and conserves the surfactant mass at all time instances $t_0, t_1, \ \cdots,  t_N$.  In the examples in this section, the velocity field $\vel(t,\bfx)$ is explicitly given, we always use a uniform time step size $\Delta t = T / N$, a uniform background mesh of the computational domain, and solve the resulting linear systems with a direct solver. Linear elements in both space and time are used. The proposed method uses the weak formulation \eqref {eq:weakformsurfdis} with $A_{h}^n(\vel, w_h,r_h)$ as in \eqref{weakformulationFE2dis} and is refered to as the conservative formulation. When $A_{h}^n(\vel, w_h,r_h)$ based on \eqref{weakformulationFE1dis} is used we refer to the method as the non-conservative method. The stabilization parameters $C_{F,1}$ and $C_{\Gamma,1}$ are in all examples equal to $10^{-2}$. 
We study the accuracy of the proposed cut finite element method by studying the following $L^2$-error at the final time $t_N=T$ 
\begin{align}
	\norm{e_w(t_N,\cdot)}_{L^2(\Gamma_h(t_N))} =\norm{w_{\text{exact}}(t_N,\cdot) - w_h(t_N,\cdot)}_{L^2(\Gamma_h(t_N))}.
	\label{eq:L2error}
\end{align}
We also report the following conservation error
\begin{equation}\label{eq:consverror}
		e_c(t_i) = \left|  \int_{\Gamma_h(t_i)} w_h(t_i, \bfx) \diff s - \int_{\Gamma_h(t_{0})} w_h(0, \bfx) \diff s  - \sum_{n=0}^{i-1} \left(\sum_{q=1}^{n_q}   \alpha_q^n  \int_{\Gamma_h(t_q^n)} f(t_q^n,\bfx) \diff s \right) =0 \right|, \quad i=0,1, \cdots,N,
\end{equation}
where we use the same quadrature rule as when we approximate the space-time integrals in the weak formulation, see Section \ref{sec:implementationsurf}.  Hence the quadrature points $t_q^n$ and weights $\alpha_q^n$ are from Simpson's quadrature rule. 

\subsubsection{Example 1 - 2D}
We consider the convection diffusion equation \eqref{eq:convection_diffusion_surfactant}-\eqref{eq:initsurf} with parameters and solution chosen as in \cite{DeElRa14}.  Initially the interface, $\Gamma(0)$, is a circle centered at origin with radius one. The interface is evolving with the velocity field :
\[
\beta (t,\bfx) = \frac{\pi}{4} \frac{\cos(2 \pi t)}{a(t)^2} (x,0), \quad  a(t) = \left( 1 + 0.25 \sin(2\pi t)  \right)^{\frac{1}{2}}.
\]
The corresponding level-set function is
\[
\phi(t,\bfx) = \frac{x^2}{a(t)^2} + y^2 -1.
\]
The diffusion coefficient $\diffSurf=1$, the right hand side $f$ and initial solution $w^0$ in equation \eqref{eq:convection_diffusion_surfactant}-\eqref{eq:initsurf} are chosen so that the exact solution is $w(t,\bfx) = x y \exp(-4t)$. 

\begin{figure}[ht]
\begin{center}
    \subfloat[$t=0$]
    	{
    	\scalebox{0.5}    	
	\centering	
	\includegraphics[scale=0.22]{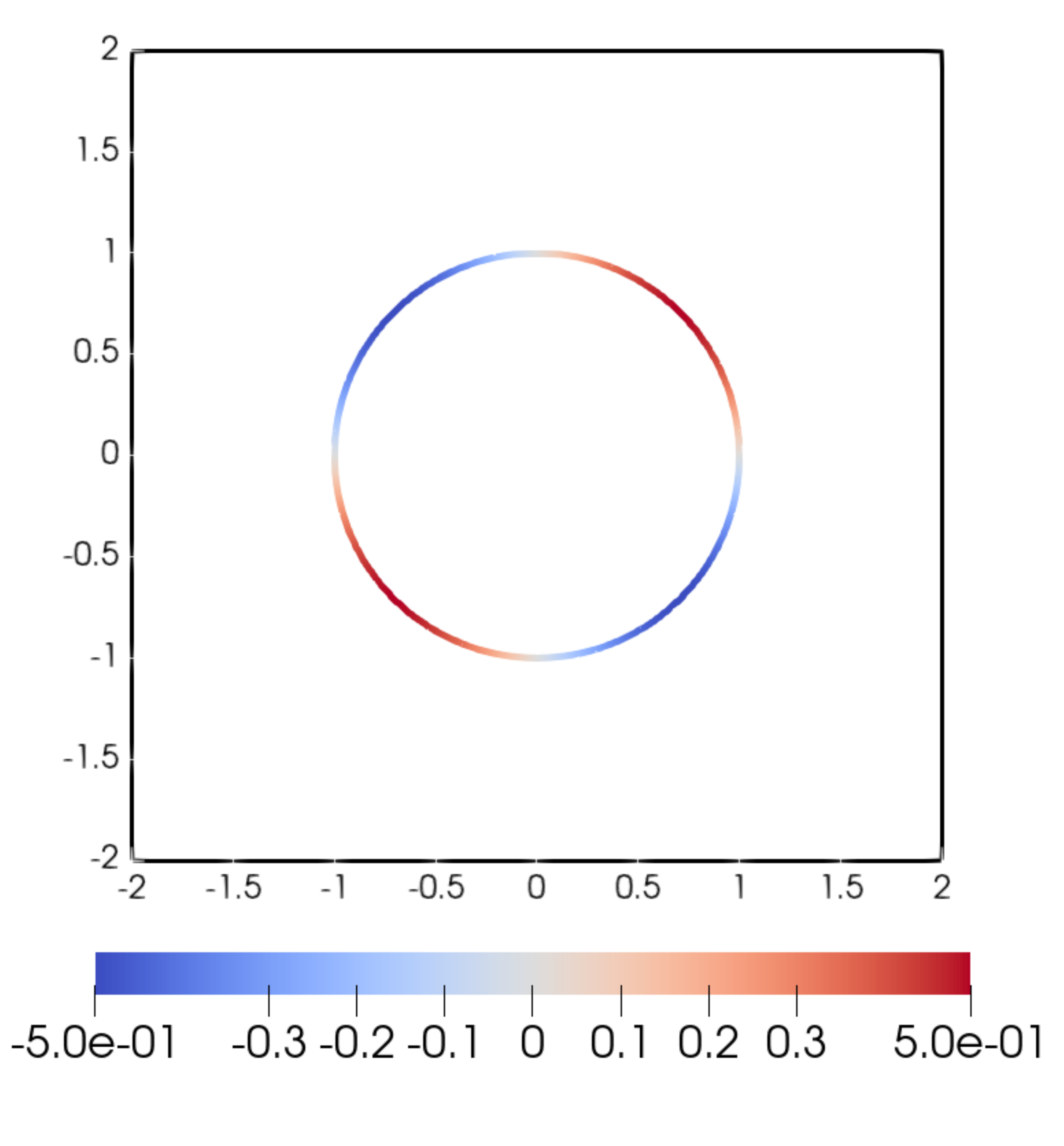}
        }  
          \subfloat[$t=0.25$]
    	{
    	\scalebox{0.5}    	
	\centering	
	\includegraphics[scale=0.22]{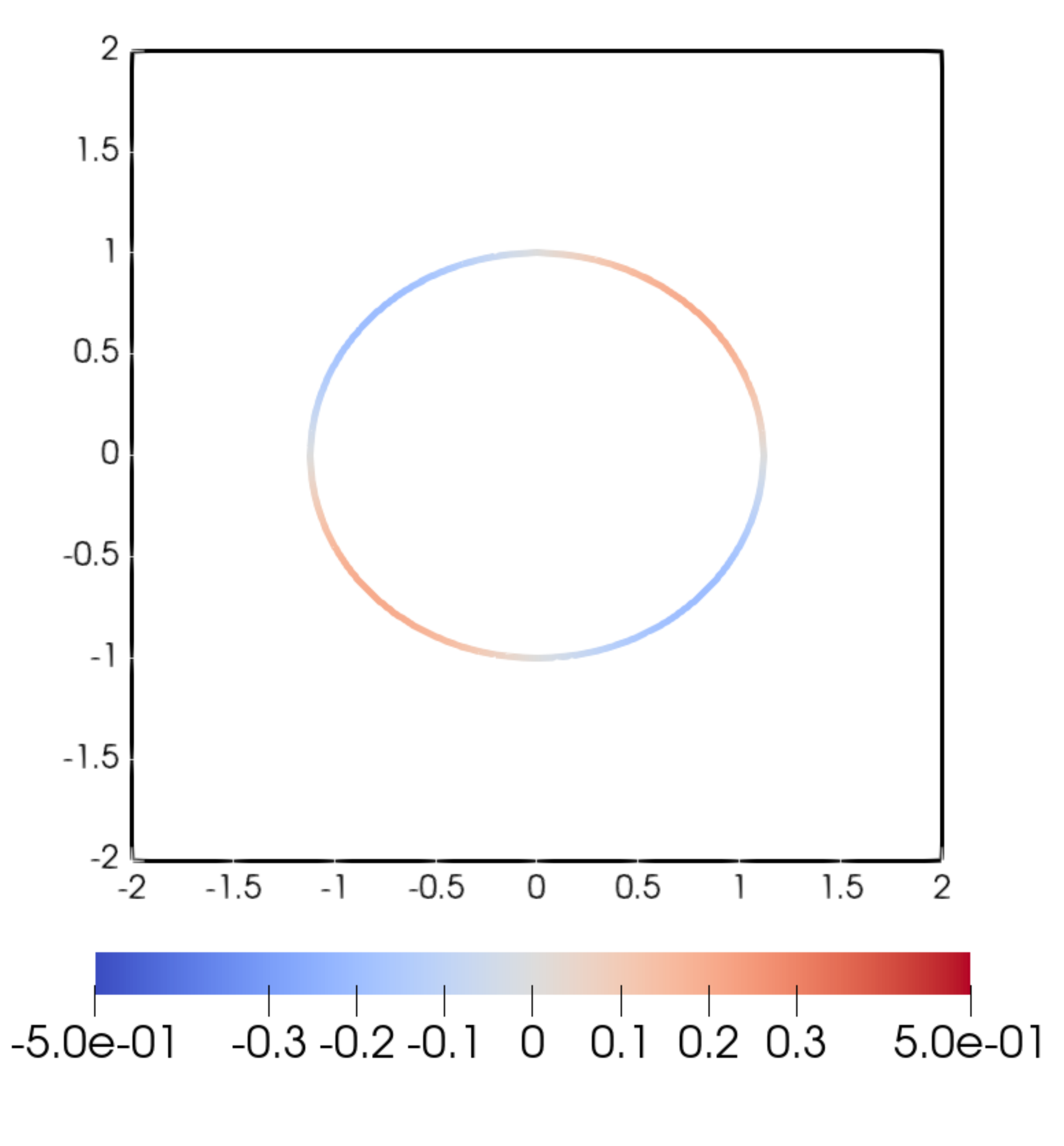}
        }           	
        \vspace{-0.4cm}
        \caption{Solution for Example 1 obtained using the proposed method with $h = 0.05$ and $\Delta t = h / 4$.}  
        \label{fig:solsurf_2D}
        \end{center}
\end{figure}

\begin{figure}[h!]
\centering	
\includegraphics[scale=0.9]{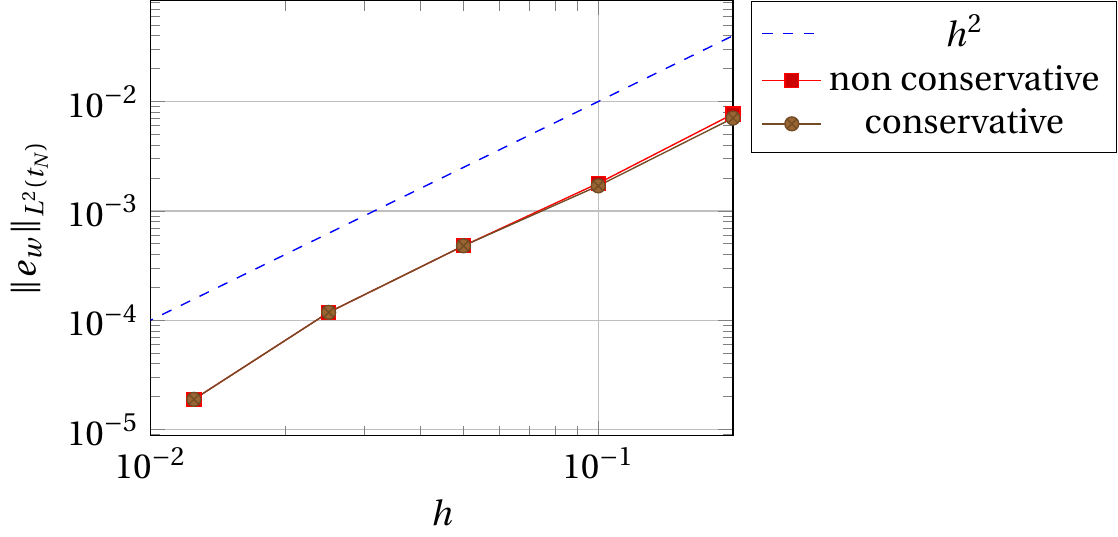}
\caption{Results for Example 1 using both the proposed scheme and the non-conservative formulation. The $L^2$-error~\eqref{eq:L2error} at $t_N=3$ versus mesh size $h$. The time step size used is $\Delta t = h / 4$. The dashed line indicates the expected second order convergence. }
\label{fig:convergence_bothFormulation_withF_2D}
\end{figure}

The computational domain is $\Omega = [-2, 2] \times [-2, 2]$ and we consider $t \in I = [0,3]$. We use the proposed finite element method and show the interface position and the concentration $w_h$ at two time instances in Figure \ref{fig:solsurf_2D}. We compute the $L^2$-error as defined in \eqref{eq:L2error} at the final time $t_N=3$ and show this error versus mesh size $h$ in Figure \ref{fig:convergence_bothFormulation_withF_2D}. Results for both the proposed scheme and the non-conservative formulation are shown. For both formulations, using $m=k=1$,  we see the expected second order convergence in the $L^2(\Gamma_h(t_N))$-norm, similar to the results reported for the unfitted finite element method in \cite{DeElRa14} (see Table 6.3 of \cite{DeElRa14}). 
\captionsetup[subfigure]{position=top}
\captionsetup[subfigure]{position=bottom}

\begin{figure}[h!]
\centering
     \subfloat[Conservative formulation]
    	{
    	\scalebox{0.5}    	
	\centering	
        \includegraphics[scale=0.8]{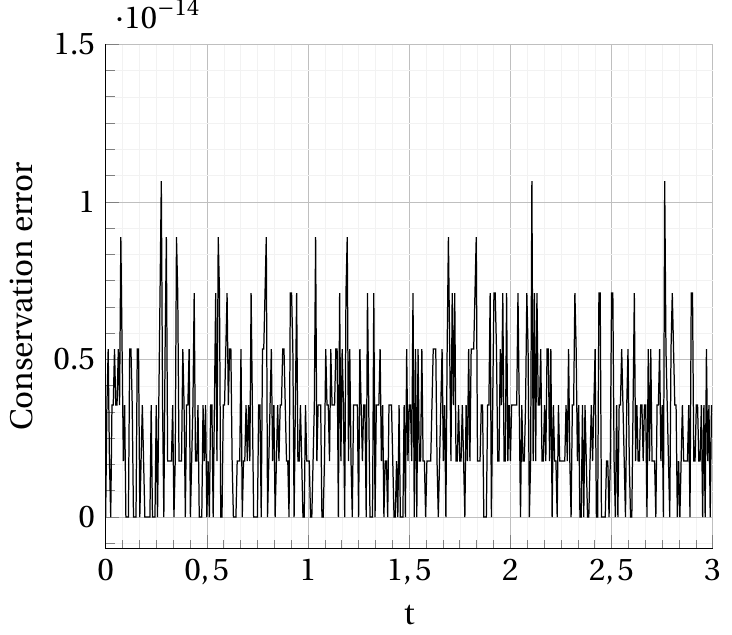}
                 \label{fig:surfactant_conservativeForm_conservation_withf_2D}

        } 
        \hspace{2cm}
        \subfloat[Non conservative formulation]
    	{
    	\scalebox{0.5}    	
	\centering	
	\includegraphics[scale=0.8]{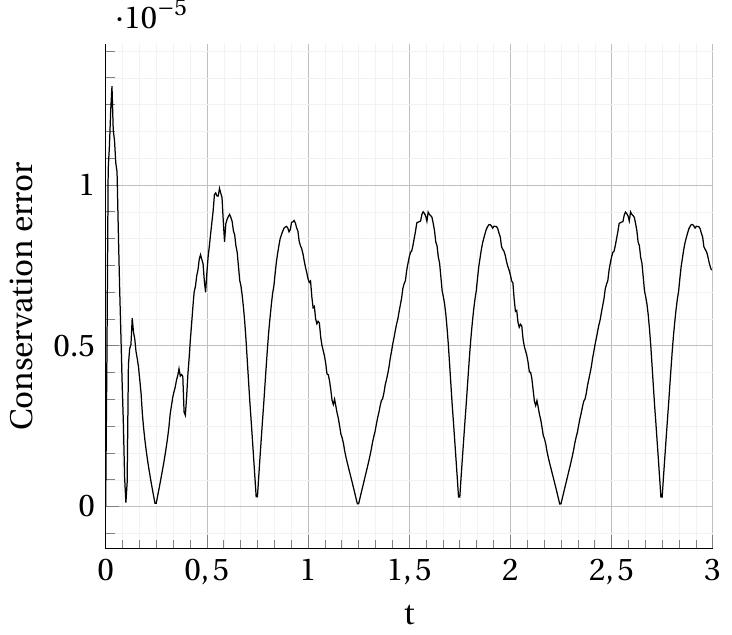}
	 \label{fig:surfactant_nonConservativeForm_conservation_withf_2D}
        }        	
        \caption{Results for Example 1 showing the conservation error \eqref{eq:consverror} using the proposed conservative formulation in (a) and the non-conservative formulation in (b). The mesh size is $h=0.05$ and the time step size is $\Delta t = h / 4$.  }  
\end{figure}

We also study the conservation of mass. 
The conservation error $e_c$, defined in \eqref{eq:consverror}, of the proposed scheme is shown in Figure \ref{fig:surfactant_conservativeForm_conservation_withf_2D} and results using the non conservative method are shown in Figure \ref{fig:surfactant_nonConservativeForm_conservation_withf_2D}. For the proposed formulation the magnitude of the conservation error is of the order of machine epsilon. In the non-conservative formulation the error $e_c$ depends on the drop deformation, the mesh size, the time step size, and the quadrature rule used. 

\subsubsection{Example 2 - accuracy test in 3D }
We follow Example 1 in \cite{OlReXu14a} and study the accuracy of the proposed space and time discretization. The computational domain is $\Omega = [-2, 2] \times [-2, 2]\times [-2, 2]$ and $I = [0,1]$. Initially the interface, $\Gamma(0)$, is a sphere centered at origin with radius equal to $1.5$ and the concentration on $\Gamma(0)$ is $w^0(\bfx)=1+xyz$. The velocity field is 
\[
\beta(t,\bfx) = \frac{3}{4}\exp\left(\frac{-t}{2}\right)\bfn.
\]
The corresponding level-set function is
\[
\phi(t,\bfx) = x^2 + y^2 +z^2 -1.5^2\exp(-t).
\]
The diffusion coefficient $\diffSurf=1$ and the right hand side of equation \eqref{eq:convection_diffusion_surfactant} is $f(t,\bfx)=(-1.5\exp(t)+\frac{16}{3}\exp(2t))xyz$ so the exact solution to equation \eqref{eq:convection_diffusion_surfactant}-\eqref{eq:initsurf} is $w(t,\bfx) =  (1+xyz)\exp(t)$. 

\begin{figure}[h!]
    \subfloat[Convergence in space for different fixed time steps $\Delta t$]
    	{
    	\scalebox{0.5}    	
	\centering	
        \includegraphics{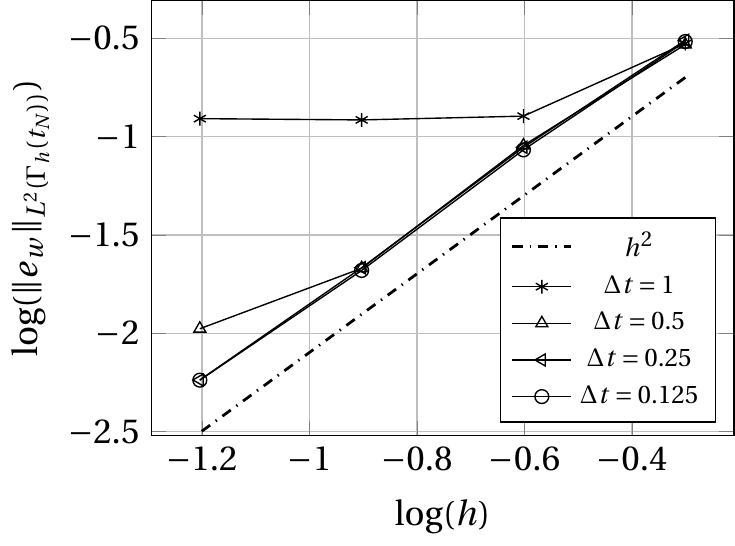}
        \label{fig:example3D_spaceConvergence}
        } 
\hfill
    \subfloat[Convergence in time for different fixed mesh sizes $h$]
    	{
	    	\scalebox{0.5}    
	\centering	
        \includegraphics{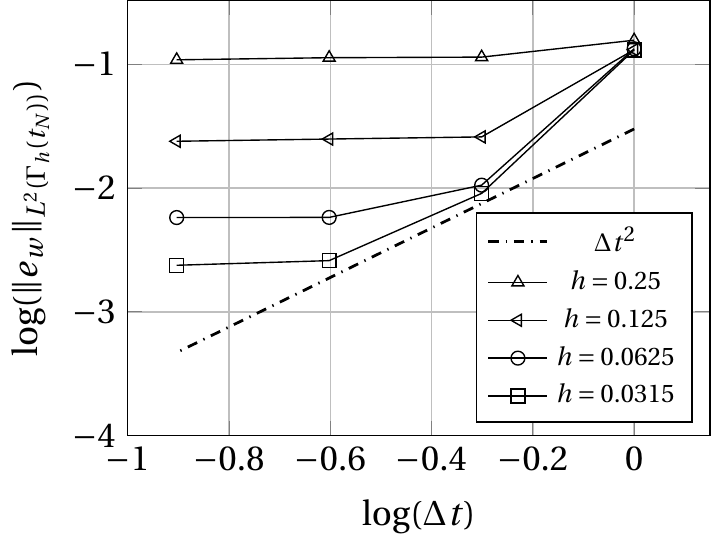}
\label{fig:example3D_timeConvergence}
}
\caption{The error in the $L^2(\Gamma_h(t_N))$-norm ($t_N=1$), see \eqref{eq:L2error}, for Example 2.}
\end{figure}

In Figure \ref{fig:example3D_spaceConvergence} we show the error in the $L^2(\Gamma_h(t_N))$-norm ($t_N=1$), see \eqref{eq:L2error}, versus mesh size $h$ for different fixed time step sizes $\Delta t$ while in Figure \ref{fig:example3D_timeConvergence} we show the error versus the time step size $\Delta t$ for different fixed mesh sizes $h$. We see that our scheme has a second order convergence rate. We also see that in this example the error is often dominated by the error in the space discretization. In our scheme in particular the geometric error,  the error coming from the approximation of the interface, is dominating in this example. The results in Figure \ref{fig:example3D_spaceConvergence}-\ref{fig:example3D_timeConvergence} can be compared with Figure 7.1 in \cite{OlReXu14a} where a different second order space-time unfitted finite element method is used, see Remark \ref{remspacetimeFEM}. However, note that for the results we report here we have approximated the level set function and the interface on the same background mesh as our spaces are defined on while the results reported in \cite{OlReXu14a} for $h$ and $\Delta t$ use an approximation of the interface that is on one regular refinement of the space-time mesh (hence with $h/2$ and $\Delta t/2$).  We observe that therefore when the error in the space discretization dominates we have slightly larger errors than reported in \cite{OlReXu14a} but when the error in the time discretization dominates (seen when larger time step sizes are used) the errors from the proposed scheme are smaller.  We study the conservation error in the next example. 

\subsubsection{Example 3 - conservation test in 3D}
\begin{figure}[ht]
    \subfloat
    	{
    	\scalebox{0.5}    	
	\centering	
	\includegraphics[scale=0.23]{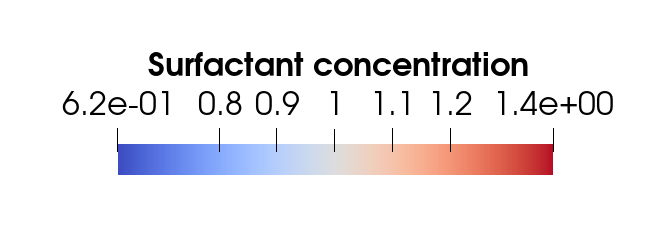}
        }  
          \subfloat
    	{
    	\scalebox{0.5}    	
	\centering	
	\includegraphics[scale=0.23]{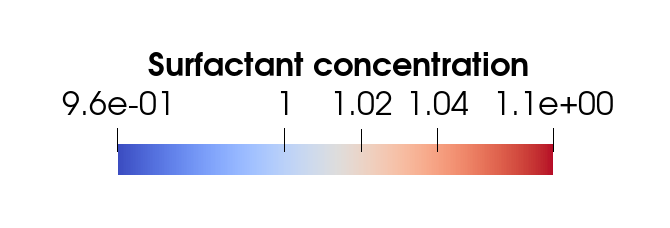}
        }          
        \subfloat
    	{
    	\scalebox{0.5}    	
	\centering	
	\includegraphics[scale=0.23]{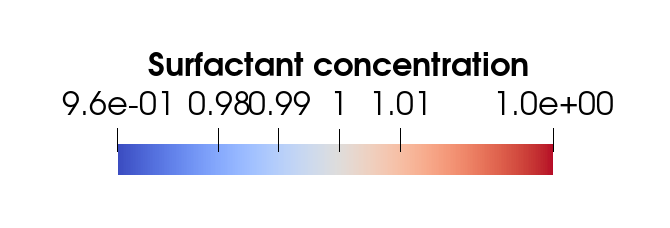}
        }     \\
        \setcounter{subfigure}{0}
    \subfloat[$t=0$]
    	{
    	\scalebox{0.5}    	
	\centering	
	\includegraphics[scale=0.25]{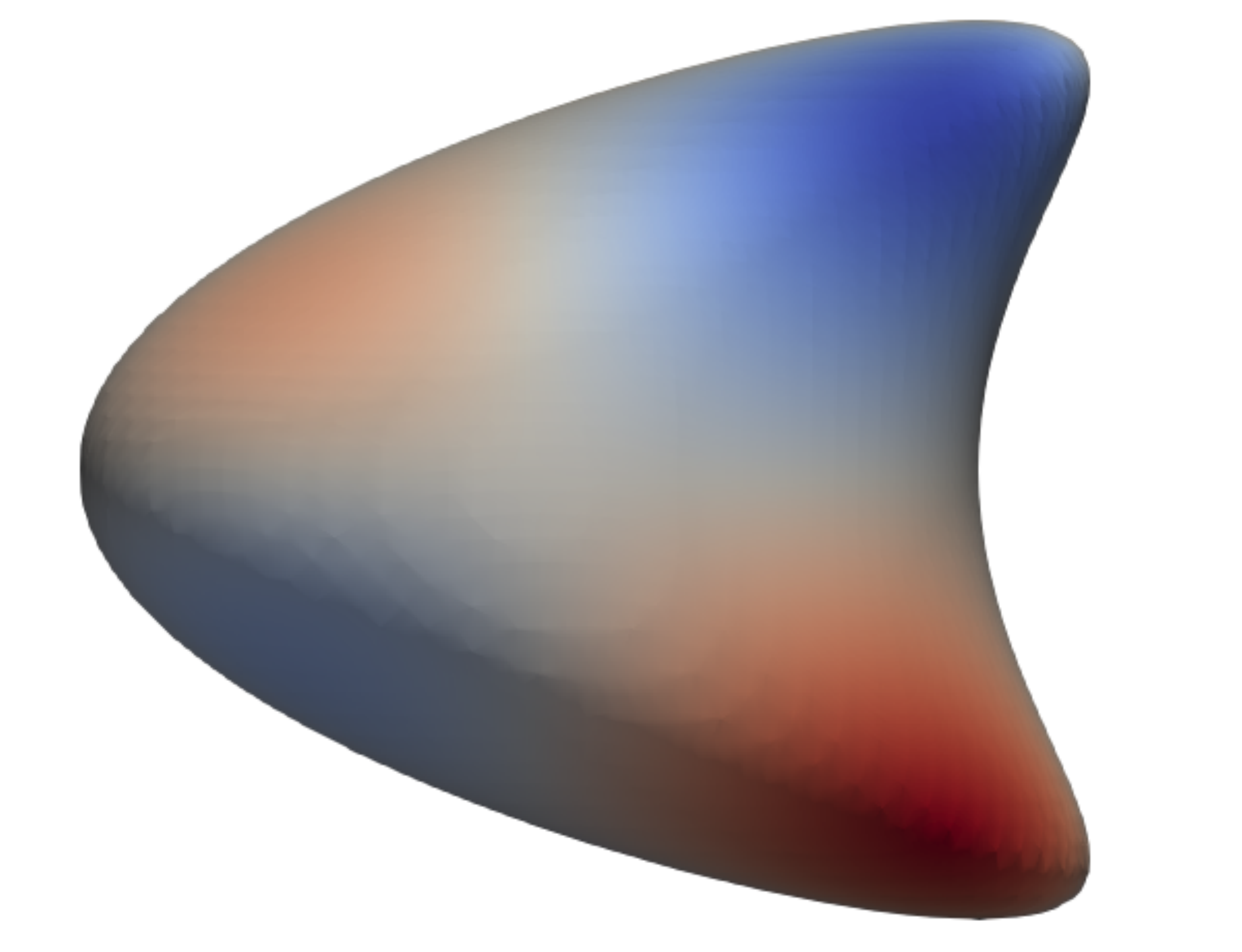}
        }  
          \subfloat[$t=2$]
    	{
    	\scalebox{0.5}    	
	\centering	
	\includegraphics[scale=0.25]{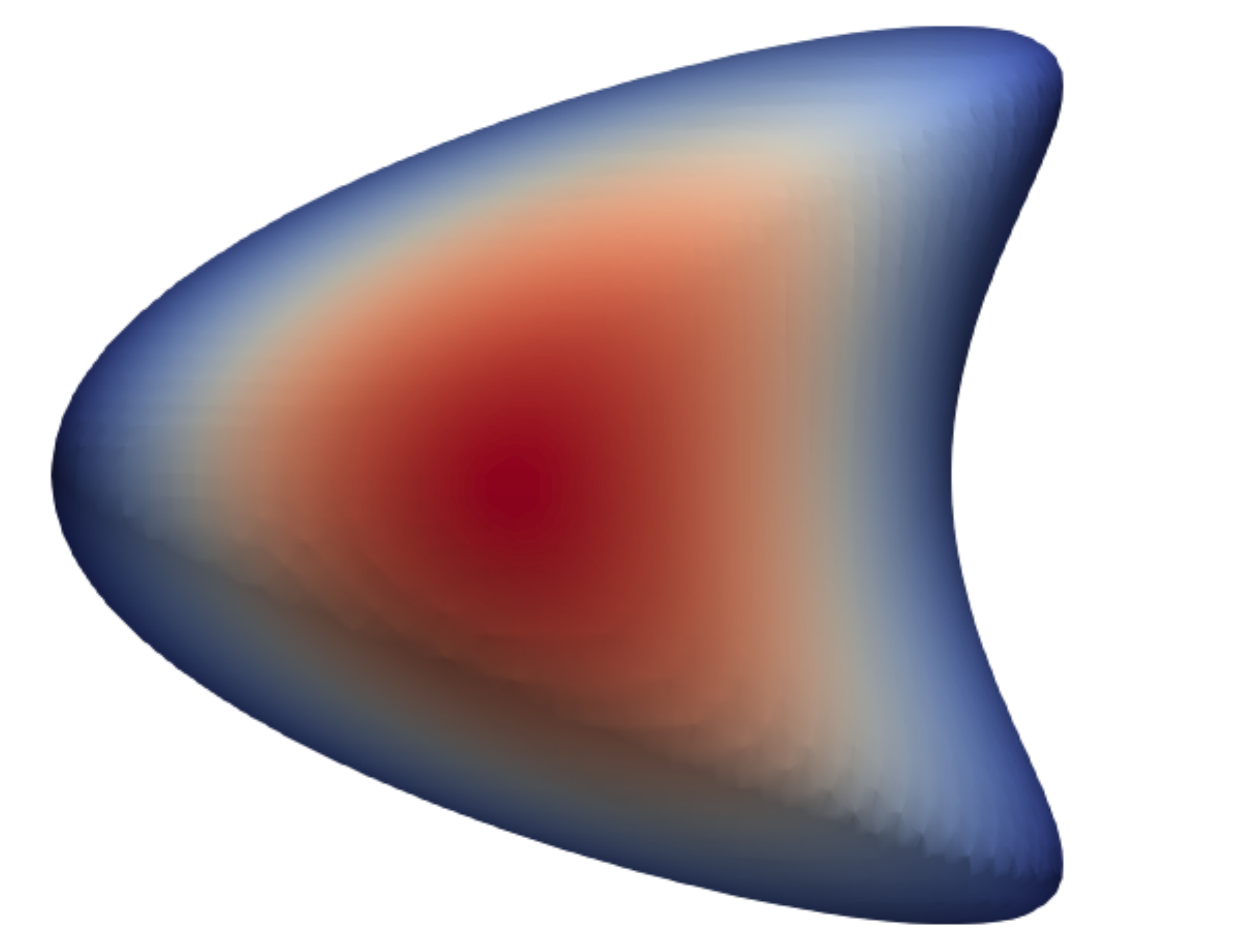}
        }          
        \subfloat[$t=4$]
    	{
    	\scalebox{0.5}    	
	\centering	
	\includegraphics[scale=0.25]{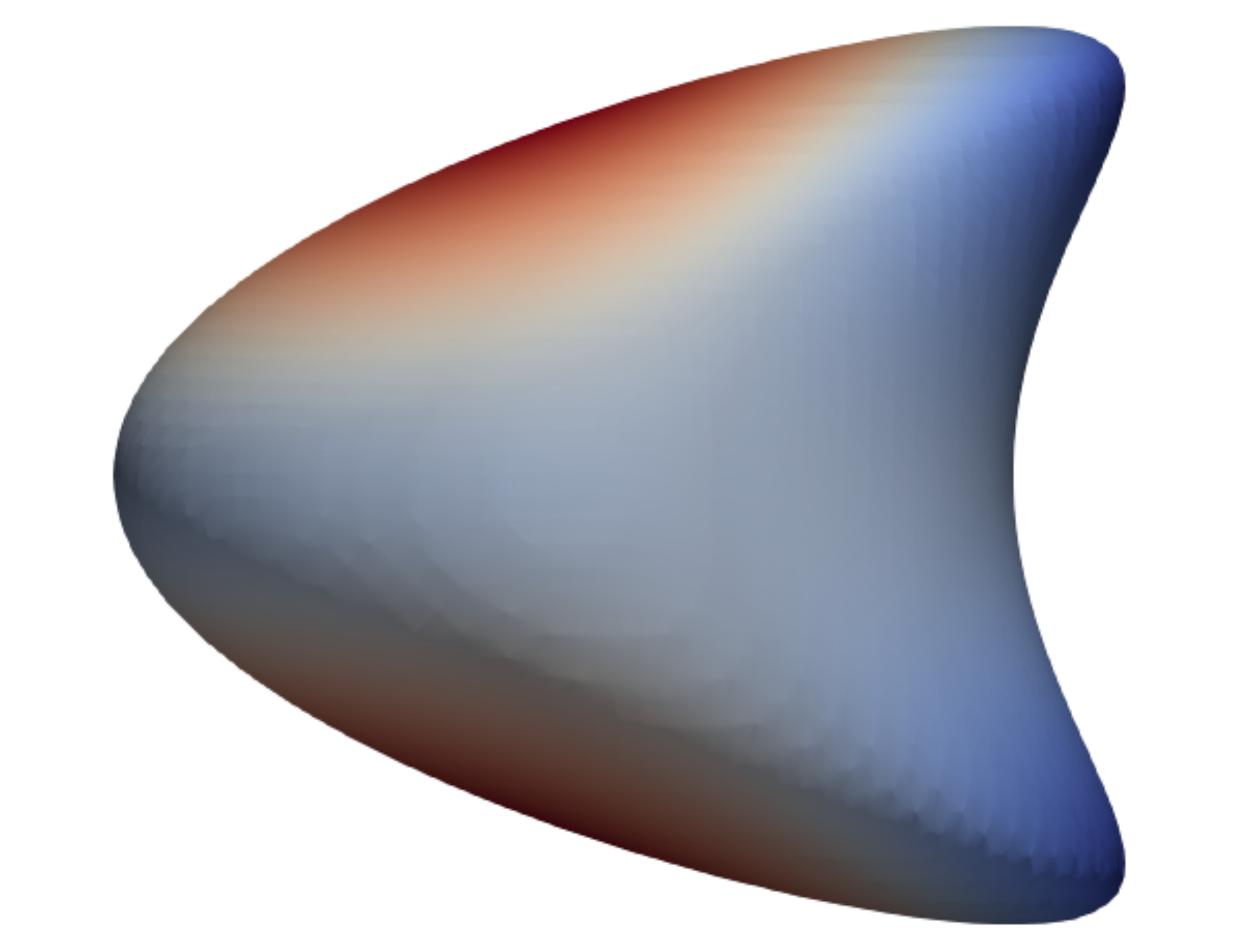}
        }      	
        \caption{Approximate solution of Example 3 using the proposed CutFEM with $h = 0.125$ and $\Delta t = 0.0625$.}  
        \label{fig:solution 3D example2 conservation}
\end{figure}

In this last example, also from \cite{OlReXu14a}, we study the conservation property of the proposed scheme in an example in three space dimensions.  The initial interface $\Gamma(0)$ is given by the zero level set of the following function
\[
\phi(0,\bfx) = (x - z^2)^2 + y^2 +z^2 - 1
\]
and the initial concentration on $\Gamma(0)$ is $w^0(\bfx)=1+xyz$.
The velocity field with which the interface is also transported is  
\[
\beta(t,\bfx)  = (0.1x\cos(t), 0.2y\sin(t),0.2z\cos(t))
\]
and the right hand side $f$ is zero. The computational domain is $\Omega = [-2, 2] \times [-2, 2]\times [-2, 2]$ and $I = [0,4]$. We choose $\Delta t = h /2$.
An approximate solution at different time instances using the proposed CutFEM is shown in Figure \ref{fig:solution 3D example2 conservation} for $h = 0.125$.   

\begin{figure}[ht]
    \subfloat
    	{
    	\scalebox{0.5}    	
	\centering	
	\includegraphics[scale=0.8]{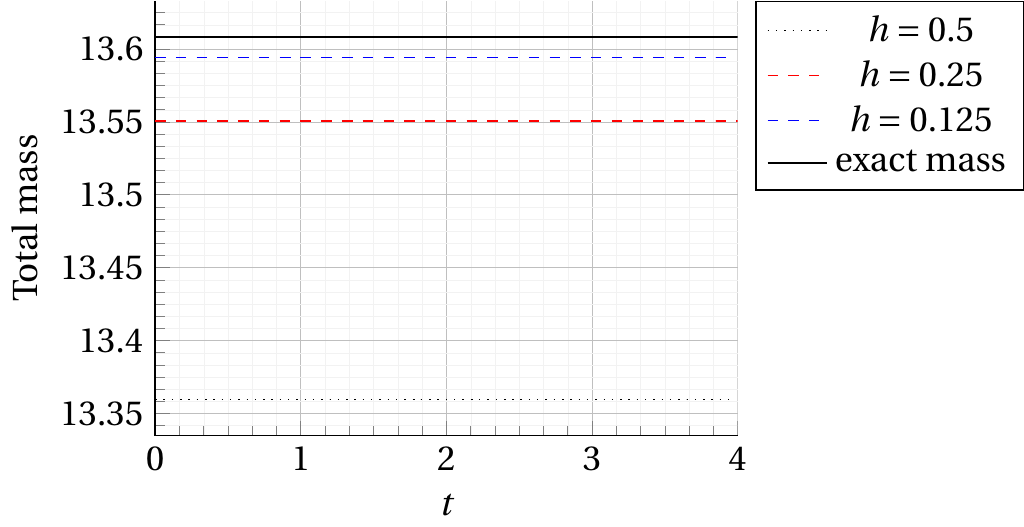}
        }  
          \subfloat
    	{
    	\scalebox{0.5}    	
	\centering	
	\includegraphics[scale=0.8]{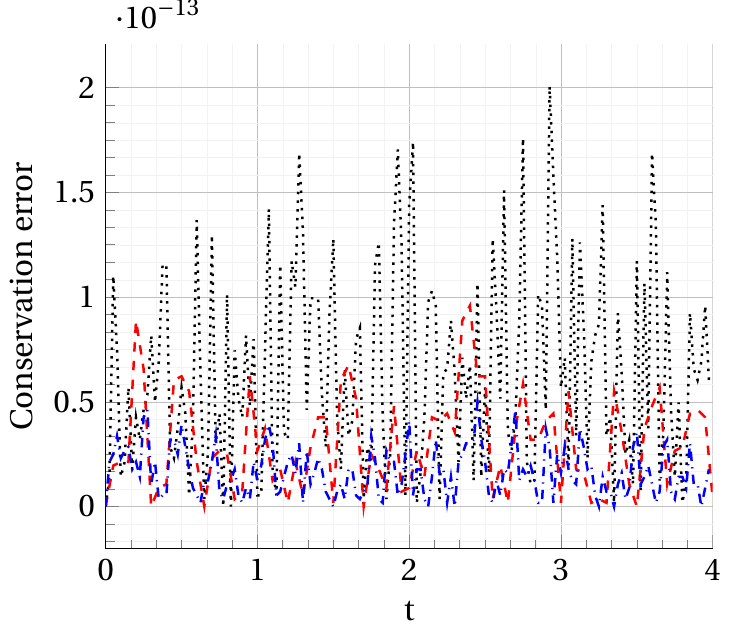}
        }     
        \caption{Results for Example 3 using the proposed space-time CutFEM. Left: The discrete total mass $M_h(t) = \int_{\Gamma_h(t)} w_h(t,\bfx) \diff s$ versus time. Right: The conservation error $e_c$ defined in \eqref{eq:consverror} (with $f=0$). }         \label{fig:3D_example3_concervation formulation 2}
 \end{figure}

We show the discrete total mass $M_h(t) = \int_{\Gamma_h(t)} w_h(t,\bfx) \diff s$ as function of time and the conservation error $e_c$ defined in \eqref{eq:consverror} (with $f=0$) in Figure \ref{fig:3D_example3_concervation formulation 2}.  Compared to the space-time unfitted finite element method in \cite{OlReXu14a} (see Figure 7.5. of  \cite{OlReXu14a}) the proposed method conserves the discrete mass.

\section{The Finite element method for the two-phase flow system}\label{sec: weakform_nummeth}
We are now ready to propose a discretization of the PDEs in Section \ref{sec:mathmod} describing a two-phase flow system with insoluble surfactants. 
Starting from the initial condition $\vel_h^-(t_{0},\bfx)=\vel^0(\bfx)$ in $\Omega_1(t_{0})\cup \Omega_2(t_{0})$ and $w_h^-(t_{0},\bfx)=w^0(\bfx)$ on $\Gamma(t_{0})$ we solve the following variational formulation, where space and time are treated similarly and interface and boundary conditions are imposed weakly, one space-time slab at a time: for $I_n$, $n=0,1, \cdots, N-1$  given $\vel_h^-(t_{n},\bfx)$ and $w_h^-(t_{n},\bfx)$ find $\vel_h \in \bfV^{n}_h$,  $p_h \in Q_h^{n}$, and $w_h \in W_h^{n}$ such that 
\begin{align}\label{eq:weakformcutfem}
F(\vel_h,p_h,w_h) &= B(\vel_h,p_h,w_h,\bfv_h,q_h,r_h) + \int_{I_n} c(t, \vel_h, \vel_h, \bfv_h)\diff t + \int_{I_n} \bff_\Gamma(t,\sigma(w_h),\bfv_h) \diff t  \nonumber \\ 
&+A(\vel_h, w_h,r_h)+ S(\vel_h,p_h,w_h,\bfv_h,q_h,r_h)- L(\bfv_h,q_h,r_h)=0, 
\end{align}
for all $\bfv_h \in \bfV^{n}_h$,  $q_h \in Q_h^{n}$, and $r_h \in W_h^{n}$. 
Here the form $A(\vel_h,w_h,r_h)$ is defined as in \eqref{weakformulationFE2}, 
\begin{align}\label{eq:fgamma}
\bff_\Gamma(t, \sigma(w_h),\bfv_h)&=\prodscal{\sigma(w_h)\nabla_{\Gamma} \bfx_{\Gamma}}{\gradS \langle \bfv_h \rangle}_{\Gamma(t)}, 
\end{align}  
\begin{align*}
c(t, \vel_h, \bfv_h,\bfw_h)&=\prodscal{\rho (\vel_h \cdot \grad {\bfv_h})}{\bfw_h}_{\tdomain},
\end{align*} 
and  
\begin{align}\label{eq:weakformAh}
B(\vel_h,p_h,w_h,\bfv_h,q_h,r_h)&=\int_{I_n} \prodscal{\rho \partial_t \vel_h}{\bfv_h}_{\tdomain} \diff t  +\prodscal{\rho \vel_h}{\bfv_h}_{\Omega_1(t_{n})\cup \Omega_2(t_{n})}
\nonumber \\
& 
+  \int_{I_n} \left[a(t, \vel_h,\bfv_h)-b(t,\bfv_h,p_h)+b(t,\vel_h,q_h)  \right] \diff t, 
\end{align}
where the forms $a(t,\bfu_h,\bfv_h)$ and $b(t,\bfv_h,p_h)$ are defined as in equation \eqref{eq:forma}-\eqref{eq:formb} with the penalty parameters $\lambda_\Gamma$, $\lambda_{\partial \Omega}$, and the weights $\w_1$, $\w_2$ in the averaging operators \eqref{eq:averagop} chosen as in \cite{FrZa19}. Further,   
\begin{align*}
	L(v_h,q_h,r_h) = \int_{I_n} l(t,\bfv_h,q_h,r_h) \diff t +\prodscal{\rho \vel_h^-}{\bfv_h}_{\Omega_1(t_{n})\cup \Omega_2(t_{n})} +\prodscal{w^-_h}{r_h}_{\Gamma(t_{n})},
\end{align*}
with $l(t,\bfv_h,q_h,r_h)$ defined as in \eqref{eq:forml} and we also have 
\begin{align*}
	 S(\vel_h,p_h,w_h,\bfv_h,q_h,r_h) = \int_{I_n} \left [ s_p(p_h,q_h) + s_\vel(\vel_h,\bfv_h) +s_w(t,w_h,r_h) \right ] \diff t,
\end{align*}
where the terms $s_p(p_h,q_h)$, $s_\vel(\vel_h,\bfv_h)$, and $s_w(t,w_h,r_h)$ are appropriate stabilization terms. We choose $s_w(t,w_h,r_h)$ as in \eqref{eq:surfstab} and 
\begin{align}\label{eq:sp}
s_p(p_h,q_h) = \sum_{i=1}^2 \sum_{F \in \mathcal{F}_{h,i}^n} 
C_p \mu_i^{-1} h^3 \prodscal{\jump{ \bfn_F \cdot \grad{p_{h,i}}}_F}{\jump{ \bfn_F \cdot \grad{q_{h,i}}}_F}_F,
\end{align}
\begin{align}\label{eq:su}
s_\vel(\vel_h,\bfv_h) = \sum_{l=1}^d \sum_{i=1}^2 \sum_{F \in \mathcal{F}_{h,i}^n}\sum_{m=1}^2
C_{\vel,m} (\mu_i,\rho_i) h^{2m-1} \prodscal{\jump{ D^m_{\bfn_F}{\vel_{h,i}^l}}_F}{\jump{D^m_{\bfn_F} {\bfv_{h,i}^l}}_F}_F.
\end{align}
Here $C_p$ and $C_{\vel,m}$ are positive constants and $\jump{v}_F$ denotes the jump of a function $v$ at the face F and is defined as $\jump{v}_F=v^+-v^-$, where $v^\pm=\lim_{t\rightarrow 0^+} v(t,\bfx \mp t \bfn_F)$, $\bfx\in F$, and $\bfn_F$ is a fixed unit normal to $F$. The stabilization terms $s_p$ and $s_\vel$ are similar as in \cite{HaLaZa14} but the sets $\mathcal{F}_{h,i}^n$ are different in the space-time method. Note that for each space-time slab the sets $\mathcal{F}_{h,i}^n$, $i=0,1,2$ do not depend on $t$, see Figure \ref{fig:illustmesh} but the trial and test functions do depend on time. These stabilization terms define an extension and therefore including them in the weak form allows us to formulate our space-time method based on quadrature in time and have a stable and robust discretization independent of how the interface cuts through the background mesh.

Note that from the weak formulation we find functions $p_h=(p_{h,1},p_{h,2}) \in Q_h^n$ and $\vel_h=(\vel_{h,1}, \vel_{h,2}) \in \bfV_h^n$ which are double valued in the interface region, see Section \ref{sec:mesh_space}. We use $(\vel_h, p_h)$ and define an approximate velocity field $\tilde{\vel}_h$ and pressure $\tilde{p}_h$ by: 
\begin{equation}\label{eq:approxvelp}
  \tilde{\vel}_h(t, \bfx)=\begin{cases} \vel_{h,1}  & \text{for $\bfx$ in } \Omega_1(t) \\
  \vel_{h,2}  & \text{for $\bfx$ in } \Omega_2(t)
  \end{cases}, \quad
   \tilde{p}_h(t, \bfx)=\begin{cases} p_{h,1}  & \text{for $\bfx$ in } \Omega_1(t) \\
  p_{h,2}  & \text{for $\bfx$ in } \Omega_2(t)
  \end{cases}, \quad  t\in I_n. 
\end{equation}

\subsection{Implementation} \label{sec: implem_coupledp}
The proposed discretization \eqref{eq:weakformcutfem} leads to a non linear system which we solve using Newton's method:  
\begin{enumerate}[-]
\item Choose initial starting guess $(\vel_{h,0}, p_{h,0}, w_{h,0}) \in \bfV^{n}_h \times Q_h^n \times  W_h^{n}$
\item While  $\norme{(\delta \vel, \delta p, \delta w)} > \varepsilon$
\begin{enumerate}[i.]
\item Solve: $DF(\vel_{h,0}, p_{h,0}, w_{h,0})(\delta \vel,\delta p, \delta w) = F(\vel_{h,0}, p_{h,0}, w_{h,0})$
\item Update ($\vel_{h,0}, p_{h,0}, w_{h,0}$) (new guess): $\vel_{h,0} = \vel_{h,0} - \delta \vel$, \ $p_{h,0} = p_{h,0} - \delta p$ and $w_{h,0} = w_{h,0} - \delta w$.
\end{enumerate}
\end{enumerate}
Here $DF(\vel_h,p_h,w_h)$ is the differential of $F$ at $(\vel_h,p_h,w_h)$,
\begin{align}\label{eq:DF}
DF(\vel_h,p_h,w_h)(\delta \vel,\delta p, \delta w)&=   
B(\delta \vel,\delta p, \delta w, \bfv_h,q_h,r_h) +\int_{I_n} \left [ c(t, \delta \vel, \vel_h, \bfv_h)+c(t, \vel_h, \delta \vel, \bfv_h) \right ] \diff t 
  \nonumber \\ 
& + \int_{I_n} \bff_\Gamma(t, \delta w \sigma'(w_h), \bfv_h) \diff t+A(\vel_h, \delta w, r_h)-\int_{I_n} \prodscal{w_h}{\delta \vel \cdot \grad{r_h}}_{\Gamma(t)} \diff t
  \nonumber \\ 
&
+S(\delta \vel,\delta p, \delta w, \bfv_h,q_h,r_h).
 \end{align}
Since we use the linear equation of state \eqref{eq:lineareqstate}, we have $\sigma'(w_h)= -\sigma_0\beta$ and 
\begin{align*}
\bff_\Gamma(t, \delta w \sigma'(w_h),\bfv_h)&=-\sigma_0 \beta \prodscal{ \delta w \nabla_{\Gamma} \bfx_{\Gamma}}{\gradS \langle \bfv_h \rangle}_{\Gamma(t)}.
\end{align*}  
If for example the Langmuir model is used which is a nonlinear model that describe the relation between the surface tension and the surfactant concentration, i.e.,  
\begin{equation}
\sigma(w_h)=\sigma_0+\beta \ln(w_\infty-w),
\end{equation}  
where $\sigma_0 \in \R_{>0}$, $\beta \in \R_{>0}$, and $w_\infty \in \R_{>0}$ are given parameters, we instead have $\sigma'(w_h)= \frac{-\beta}{w_\infty-w_h}$. 

In the first step of the above algorithm we choose the initial guess $(\vel_{h,0}(t), p_{h,0}(t), w_{h,0}(t))$ to be equal to the solution from the previous space-time slab at time $t_{n}$, i.e., $(\vel_h^-(t_n),p_h^-(t_n),w_h^-(t_n))$. In the assembly of the linear system in step i.  we approximate each space-time integral in $F$ and $DF$  by first using a  quadrature rule in time and then at each quadrature point in time the integral in space is computed. We approximate all the integrals in time using Simpson's quadrature rule. See also Section \ref{sec:implementationsurf}.  

Here, the exact level set function is not known explicitly except at time $t=0$. We find an approximation $\lsf_h(t,\bfx)$ at each time instance $t \in T^n=\{t_n, \{t_q^n\} , t_{n+1}\}$ by discretizing the level set equation \eqref{eq:Level set} using the Crank-Nicolson scheme in time and continuous FEM with quadratic Lagrange elements in space and streamline diffusion stabilization following~\cite{FrZa19}. As in Section~\ref{sec:implementationsurf} we find a piecewise planar $\Gamma_h(t)$ as the zero level set of $\lsf_h(t,\bfx) \in \Bspace^1$, for all $t\in T^n$. We then also have approximations $\Omega_{h, 1}$ and $\Omega_{h,2}$ of the two subdomains $\Omega_{1}$ and $\Omega_{2}$. The contribution from integration on $\Omega_{h,i}(t_q^n) \cap K$ is divided into contributions on one or several triangles in two dimensions and tetrahedra in three dimensions depending on how the interface cuts element $K$.

\section{Numerical examples}\label{sec:numexpcoup}
We consider several examples in two space dimensions including drop-drop interactions and one example in three space dimensions. In the last simulation, in three space dimensions, we use the proposed scheme but instead of linear element in time we use piecewise constant polynomials in time, i.e. $k=0$ (see Section \ref{sec:mesh_space}). This corresponds to using the implicit Euler method in time (see Remark 3.2 in \cite{FrZa19}). The stabilization parameters are all chosen to be $10^{-2}$. The resulting linear systems are always solved with a direct solver.
\subsection{Rising drop}\label{example rising drop 2D}
\begin{figure}[ht]
	\centering
	\includegraphics{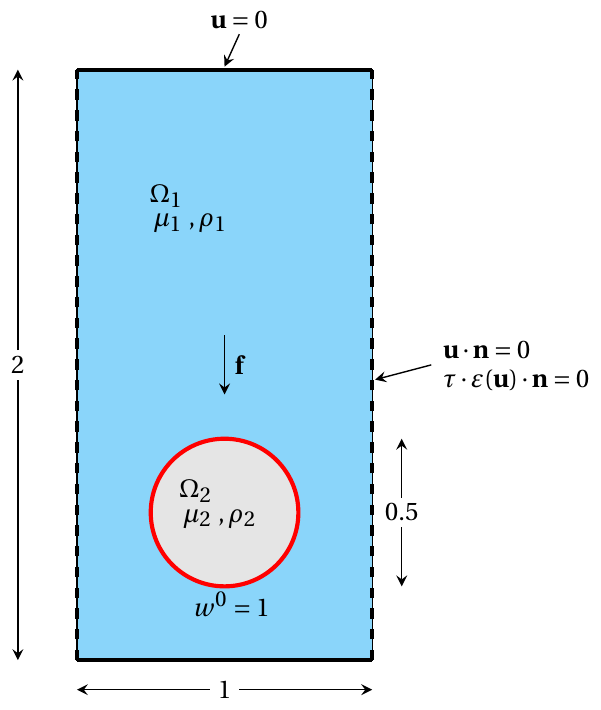}
	\caption{Initial configuration \label{fig:InitialConfiguration}}
\end{figure}
Consider a two-dimensional drop rising in a liquid column due to buoyancy with physical parameters as in Table \ref{table - set of parameters}.  The force of gravity is $\bff = \rho (0,-0.98)$. This is a benchmark test case from \cite{Benchmark} but we also add surfactant as in \cite{BaGaNu15}. The computational domain in space is $\Omega=[0,1] \times [0,2]$ and the drop is initially a circle centered at $(0.5, 0.5)$ with radius equal to $0.25$. The no-slip boundary condition, $\bfu = 0$, is imposed on the horizontal walls and the free slip condition, $\bfu \cdot \bfn = 0$ , $\tau \cdot 2\varepsilon(\bfu) \bfn = 0$, is imposed on the vertical walls (here $\bfn$ and $\tau$ are the unit normal- and tangent vectors to $\partial \Omega$). For an illustration of the initial setup see Figure \ref{fig:InitialConfiguration}.  
\begin{table}[ht]
\centering
\begin{tabular}{  p{1cm} p{1cm} p{1cm} p{1cm} p{1cm} p{1cm} p{1cm}p{1cm}}
\hline 
 $\rho_1$ & $\rho_2$ & $\mu_1$ & $\mu_2$ & $\sigma_0$ & $\beta$ & $w_0$ & $\diffSurf $\\ 
\hline
	 1000	&	100	&	10	&	1	&	24.5  & 0.5 &  1 & 0.1\\
\hline 
\end{tabular}
\captionof{table}{Parameters used in Example \ref{example rising drop 2D}. \label{table - set of parameters}}
\end{table}

\begin{figure}[h!]
\centering
\hspace*{-1cm}
\includegraphics[scale=0.35]{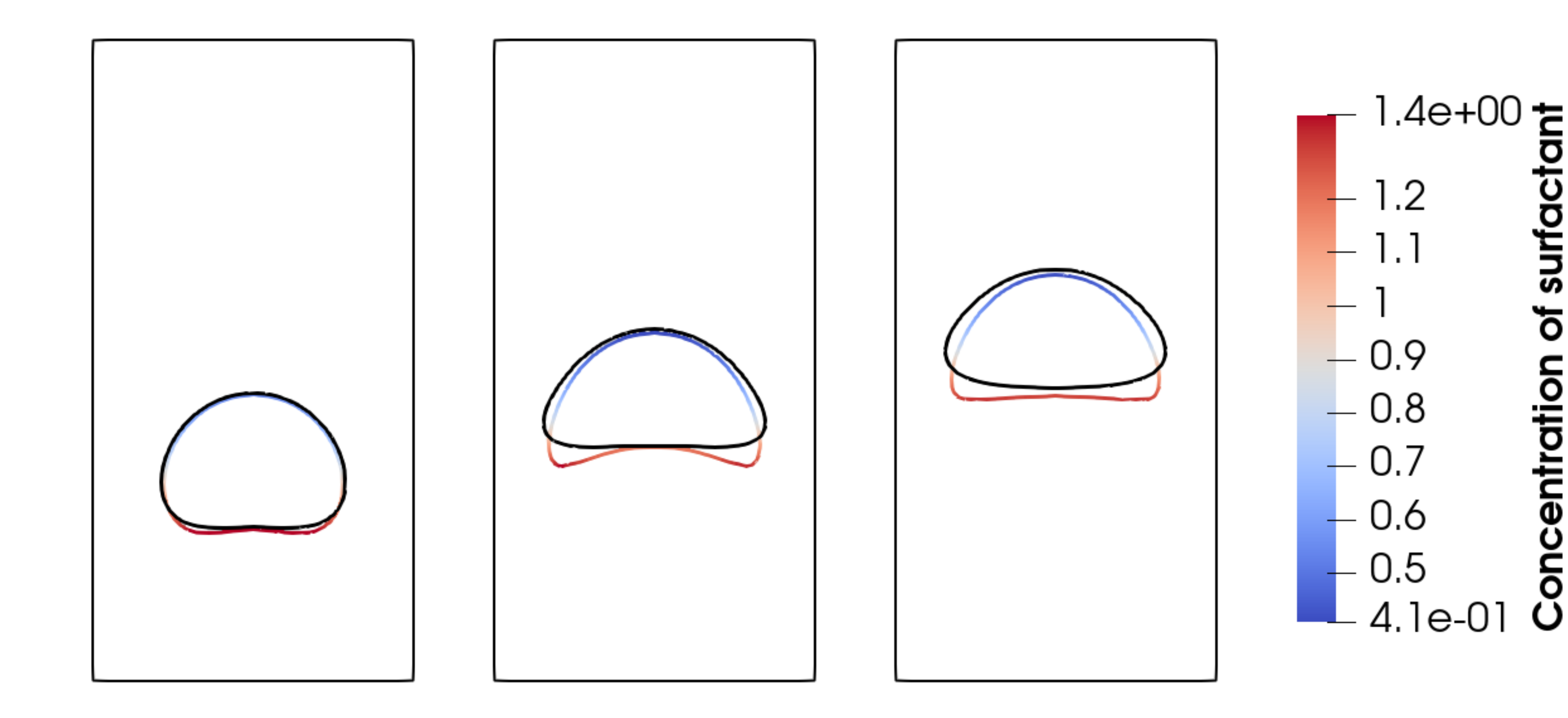}
        \caption{Drop shape and surfactant concentration in Example \ref{example rising drop 2D} at time instances $t=1,2,3$ for $\beta = 0.5$. The solid black line shows the drop shape in case of a clean interface, i.e., $\beta=0$. \label{fig:multiple shape}}   
\end{figure}

We use a background mesh with more elements in the region where the drop is evolving (the mesh can be seen in Figure \ref{jumpKink_Benchmark1}). The coursest mesh size is $h_\text{outer} = 1/40$ while the finest mesh size is $h_\text{inner} = 1/80$. We choose a time step size $\Delta t = h_\text{outer} / 4$. The evolution of the drop and the surfactant concentration is shown in Figure \ref{fig:multiple shape}. One can see that surfactant accumulate on the bottom and increase the deformation of the interface compared to the case of a clean interface, $\beta = 0$. 

\begin{figure}[h!]
    \subfloat[Rise velocity]
    	{
    	\scalebox{0.3}    	
	\centering	
	\includegraphics[scale=0.7]{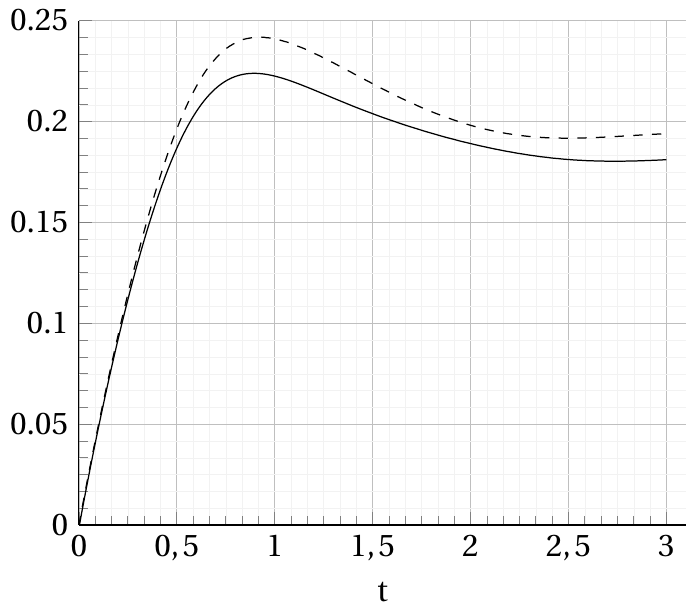}
        }
    \subfloat[Circularity]
    	{
    	\scalebox{0.3}    	
	\centering	
        \includegraphics[scale=0.7]{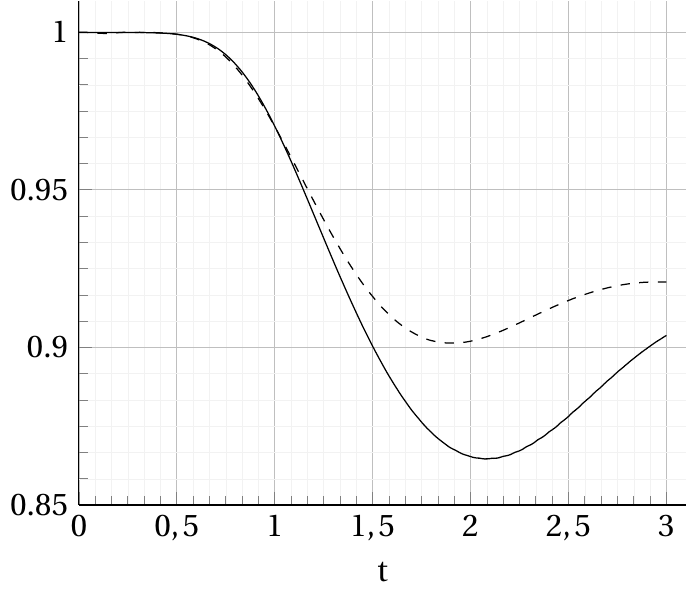}
        }
    \subfloat[Center of mass]
    	{
        \scalebox{0.3}    	
	\centering	
	\includegraphics[scale=0.7]{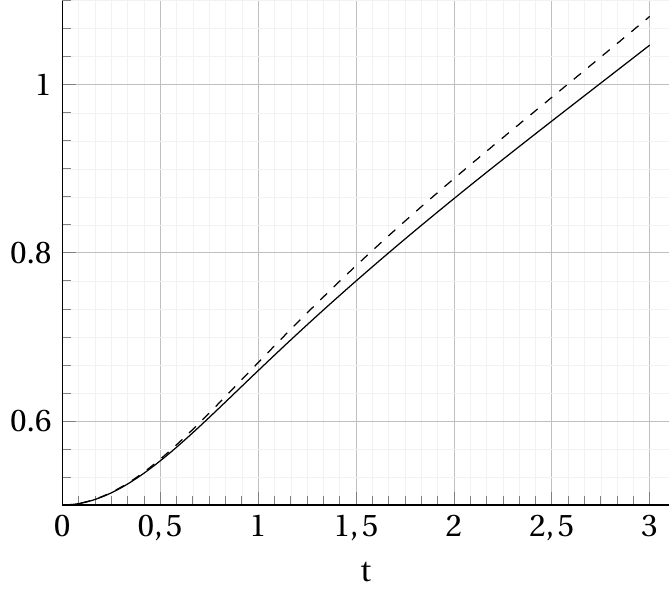}
        }                   	
        \caption{The rise velocity, circularity, and the center of mass for Example \ref{example rising drop 2D}. Comparison between surfactant being present $\beta =  0.5$ (plain line) and a clean interface, i.e., $\beta = 0$ (dashed line). }  
        \label{fig:bench1 quantities}
\end{figure}

We compute the following benchmark quantities:
\begin{itemize}
\item Center of mass 
\begin{align*}
\mathbf{X}_c = (x_c, y_c) = \frac{\int_{\Omega_2} \bfx \diff x}{\int_{\Omega_2} 1 \diff x}.
\end{align*}
Here, the second component $y_c$ is of interest.
\item Circularity 
\begin{align*}
c = \frac{P_a}{P_b},
\end{align*}
where $P_b$ is the perimeter of the drop and $P_a$ is the perimeter of the circle which has an area equal to that of the drop.
\item Rise velocity 
\begin{align*}
\mathbf{U}_c = (u^1_c, u^2_c) = \frac{\int_{\Omega_2} \bfu \diff x}{\int_{\Omega_2} 1 \diff x}.
\end{align*}
Here, we are interested in the velocity component $u^2_c$ which is in the direction opposite to the gravitational vector $\bff$.
\end{itemize}
In Figure \ref{fig:bench1 quantities} these quantities are shown versus time $t$ both when surfactant is present and for a clean interface which we studied carefully in~\cite{FrZa19}. When surfactant is present, the deformation increases and the rise velocity decreases. This impact also the center of mass. In Table \ref{table - result rising drop 2D} we report some characteristic values. Our results agree well with the results reported for the parametric finite element method in \cite{BaGaNu15} (See Table 2 and Figure 7 and 8 in \cite{BaGaNu15}).
\begin{table}[ht]
\centering
\begin{tabular}{ p{2cm} p{3cm} p{3cm}  }
\hline
 & without surfactant &  with surfactant\\
\hline
$c_\text{min}$	& 0.9015 & 0.8632	\\
$t_{c_\text{min}}$	& 1.9016 & 2.1125	\\
$u^2_{c, \text{max}}$ & 0.2417	& 0.2239	\\
$t_{u^2_{c, \text{max}}}$ & 0.9203	& 0.8969	\\
$y_{c}(t=3)$ & 1.0817& 1.0473 \\
\hline   
\end{tabular}
\captionof{table}{Some characteristic values for Example \ref{example rising drop 2D}. \label{table - result rising drop 2D}}
\end{table}

\begin{figure}[h!]
\centering
\includegraphics[scale=0.7]{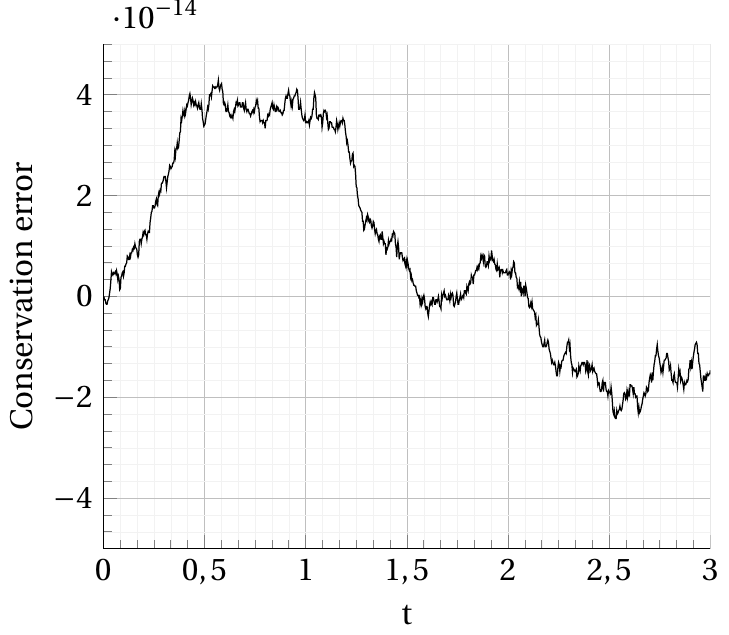}
      \caption{Conservation error for surfactant mass, Example \ref{example rising drop 2D}.} \label{fig:conserrorrising}   	
\end{figure}

\begin{figure}[h!]
  \begin{center}
\includegraphics[scale=0.25]{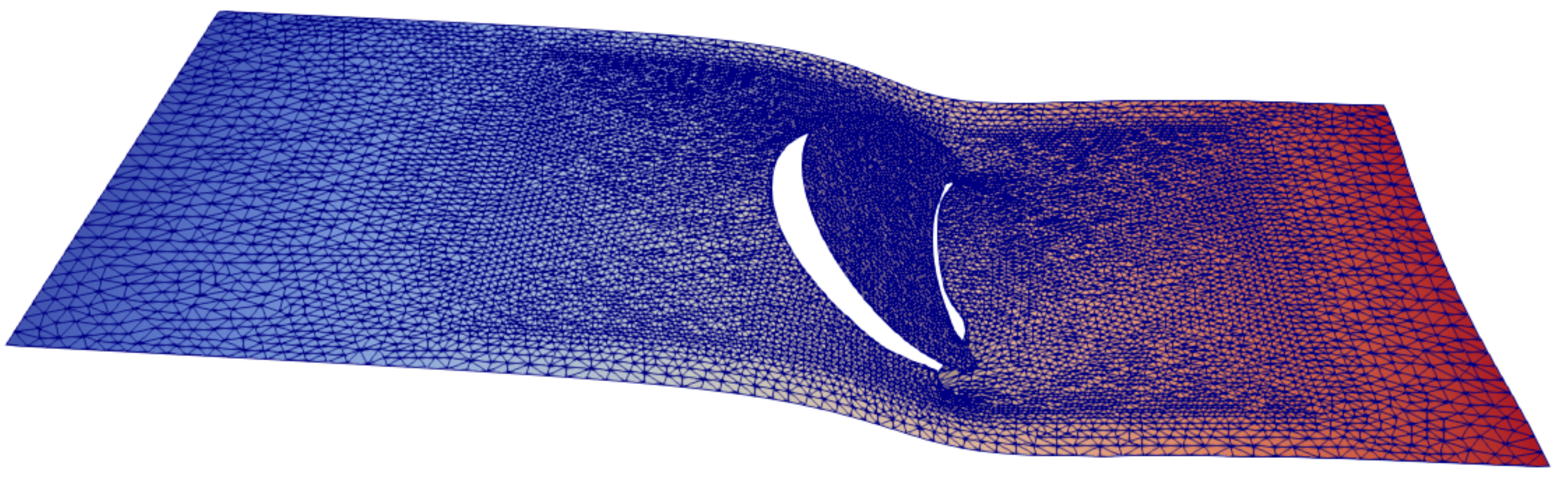}\\ 
\includegraphics[scale=0.15]{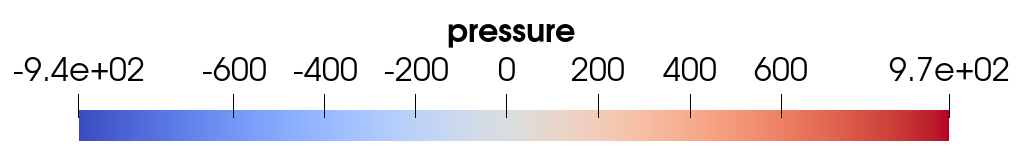} 
\caption{The approximated pressure at time $t=2$ in the presence of surfactant for Example \ref{example rising drop 2D}.\label{jumpKink_Benchmark1}}
\end{center}
\end{figure}

In Figure \ref{fig:conserrorrising} we show the conservation error which is of the same order as machine epsilon. Thus, the proposed scheme conserves the total surfactant mass. Finally, in Figure \ref{jumpKink_Benchmark1} we show the approximated discontinuous pressure at $t=2$. We see that the discontinuity across the evolving interface is accurately approximated. 

\subsection{Drop in shear flow}\label{sec:Dropinshear}
We now consider an initially circular drop, centered at the origin with radius one, in a shear flow experiment from \cite{LaiTseHua08}. We study the influence of surfactant on the deformation of the drop. The computational domain is $\Omega=[-5,5] \times [-2, 2]$ and the Dirichlet boundary condition $\bfu(t,\bfx) = (0.5 y, 0)$ is imposed on $\partial \Omega$. The physical parameters are chosen as in Table \ref{table - set of parameters_shearFlow} and $\bff=0$.
\begin{table}[ht]
\centering
\begin{tabular}{ p{3cm} p{1cm} p{1cm} p{1cm} p{1cm} p{1cm} p{1cm} p{1cm}p{1cm}}
\hline 
Test case  & $\rho_1$ & $\rho_2$ & $\mu_1$ & $\mu_2$ & $\sigma_0$ & $\beta$ & $\diffSurf$ & $w_0$\\ 
\hline
1	& 1	&	1	&	0.1	&	0.1	&	0.2  & 0.0 & 0.1&  1\\
2	& 1	&	1	&	0.1	&	0.1	&	0.2  & 0.25 & 0.1&  1\\
3	& 1	&	1	&	0.1	&	0.1	&	0.2  & 0.5 & 0.1&  1\\
\hline 
\end{tabular}
\captionof{table}{Parameters used in Example~\ref{sec:Dropinshear}. \label{table - set of parameters_shearFlow}}
\end{table}

\begin{figure}[h!]
	\centering
 	\includegraphics[width=0.75\linewidth]{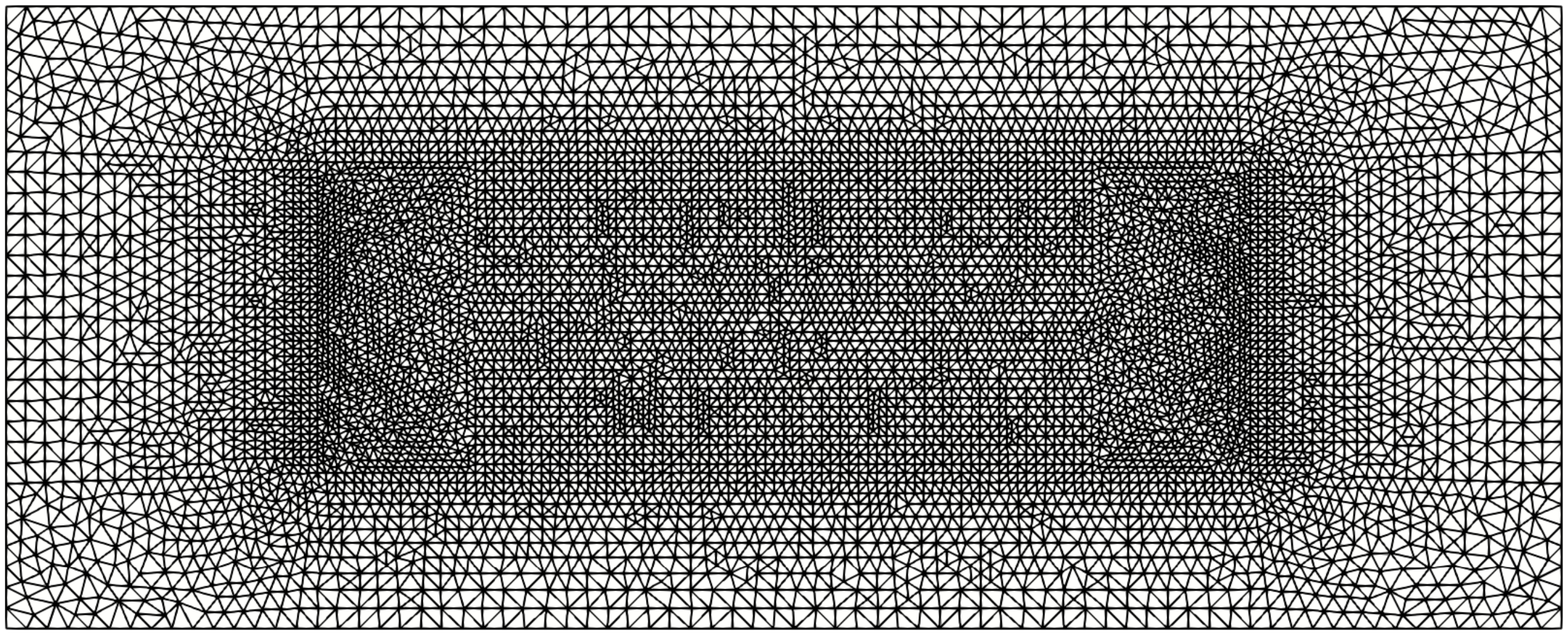}
	 \caption{Mesh used for the simulations of Example~\ref{sec:Dropinshear}. It contains 7669 nodes and 15126 elements.}  
	 \label{fig:shearMesh}
\end{figure}

\begin{figure}[h!]
    \subfloat
    	{
    	\scalebox{0.25}    	
	\centering	
	\includegraphics[width=0.25\linewidth]{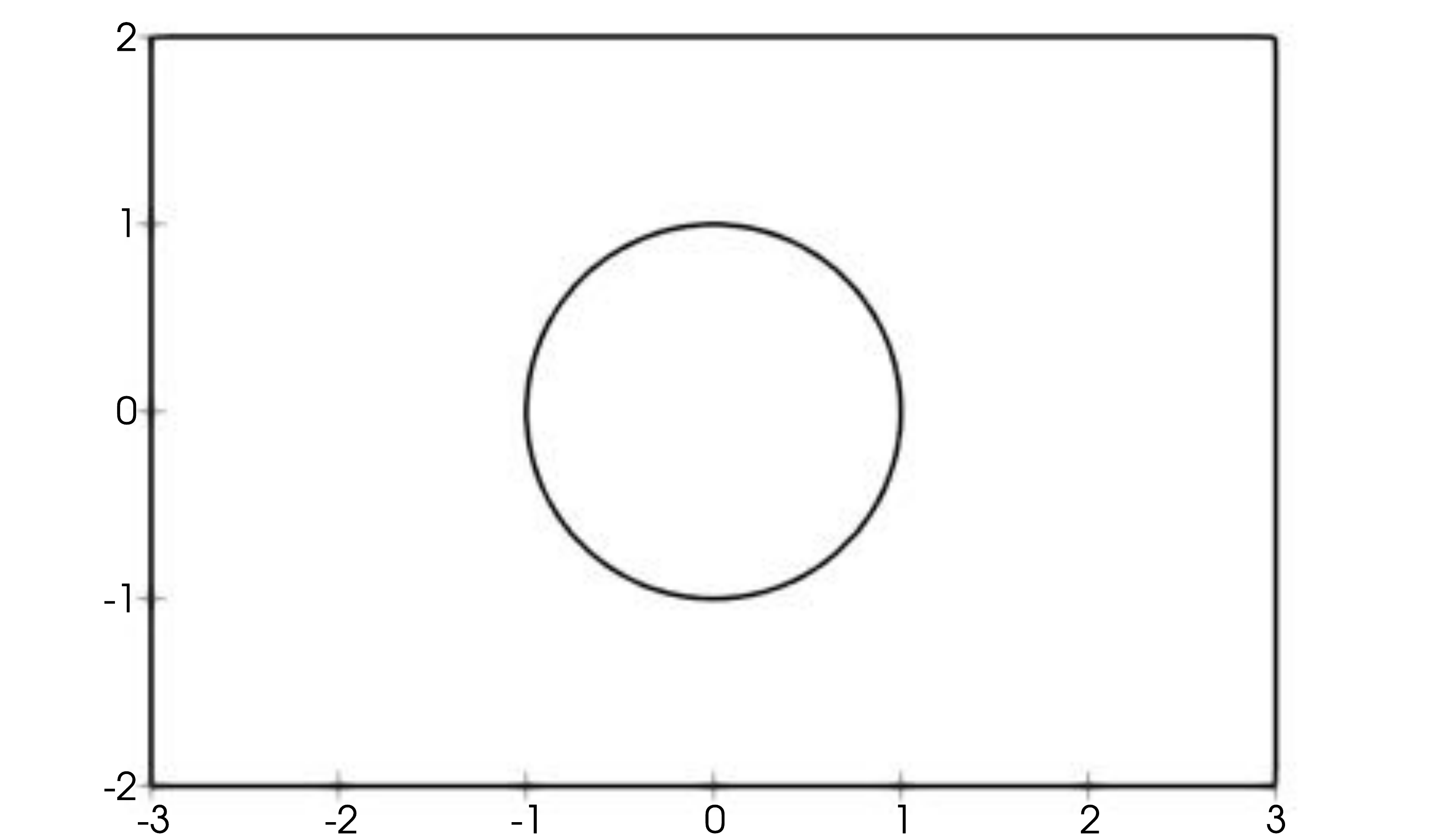}
        }
    \subfloat
    	{
    	\scalebox{0.25}    	
	\centering	
        \includegraphics[width=0.25\linewidth]{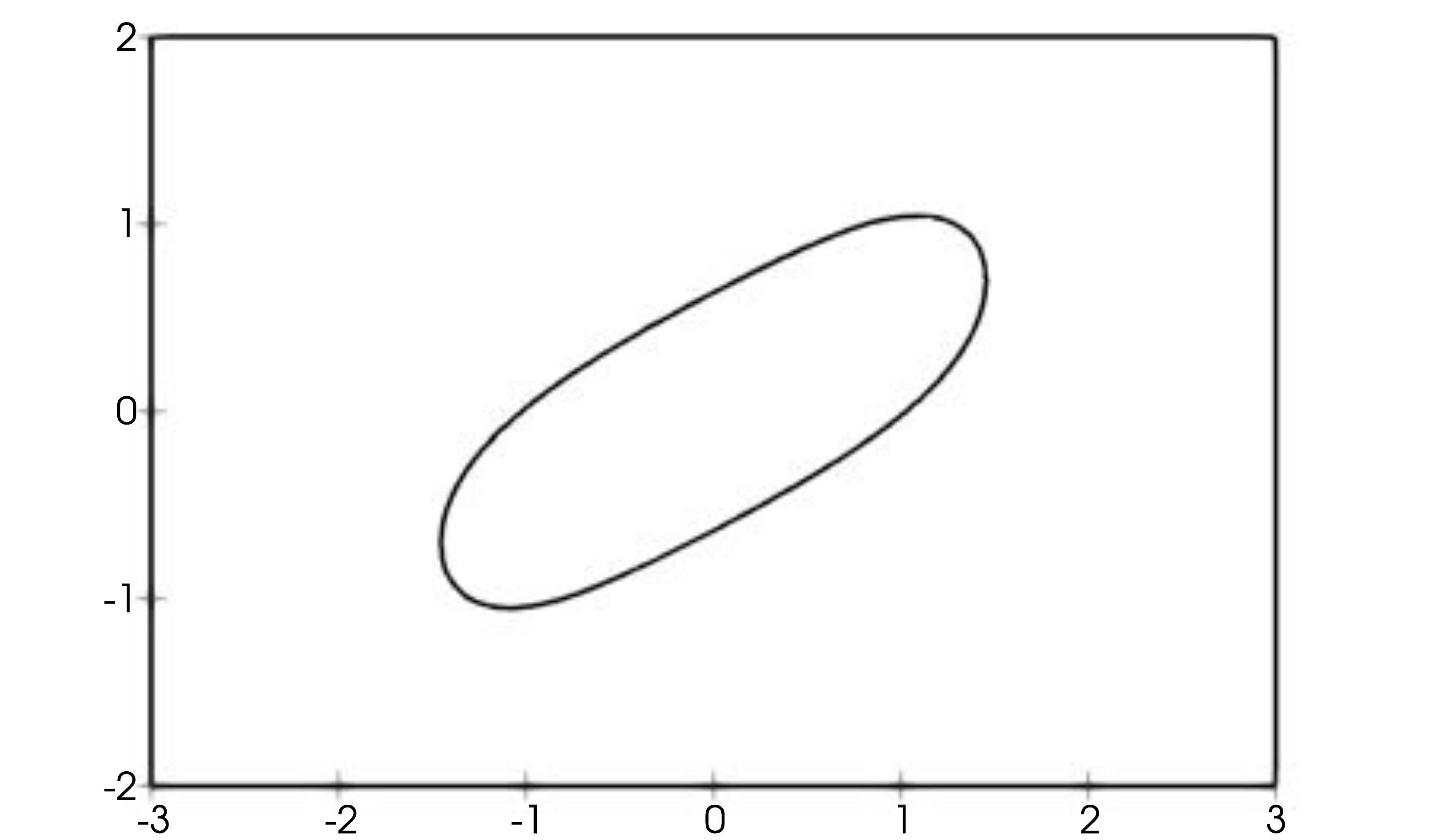}
        }
    \subfloat
    	{
        \scalebox{0.25}    	
	\centering	
	\includegraphics[width=0.25\linewidth]{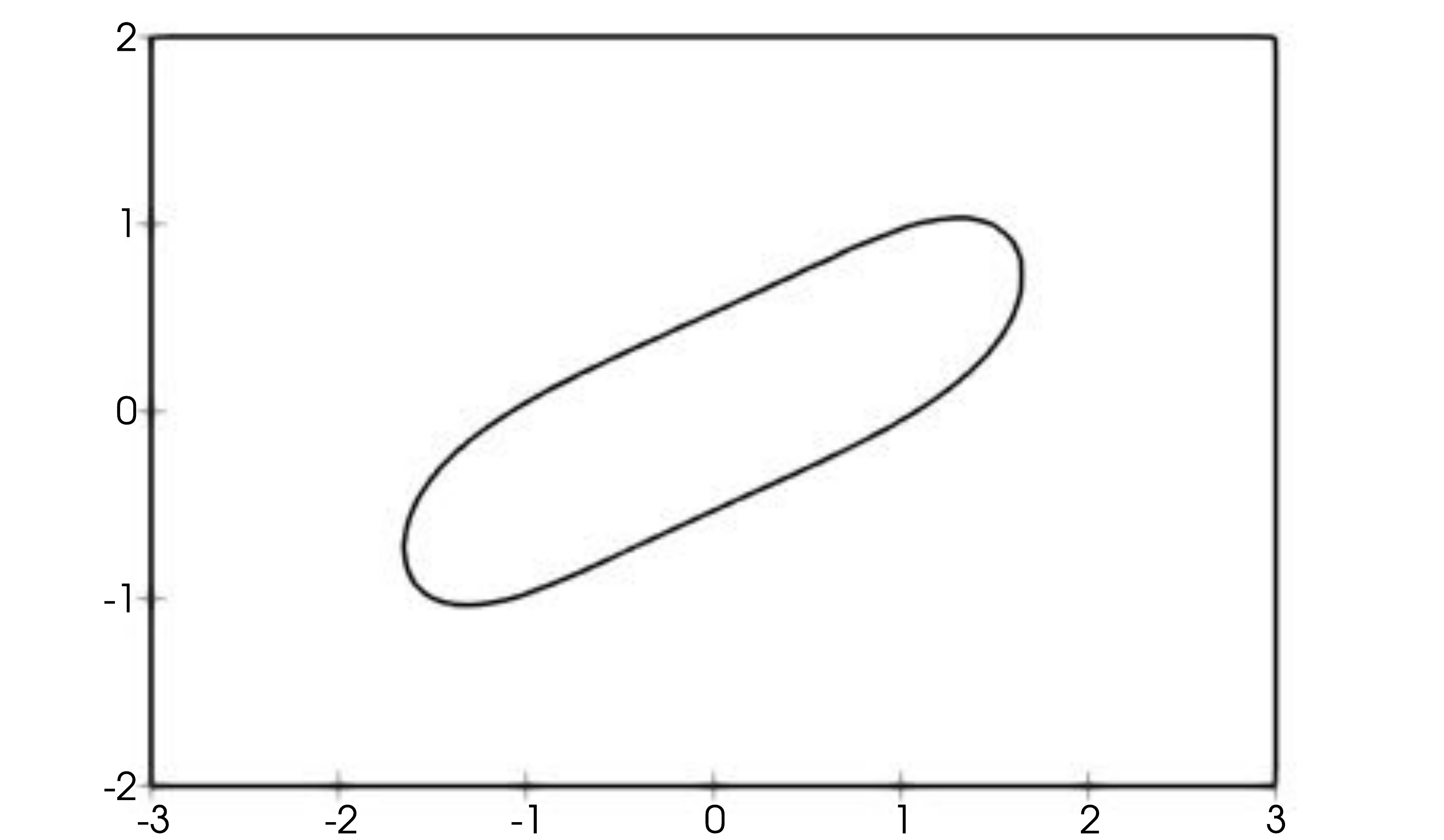}
        }            
     \subfloat
    	{
        \scalebox{0.25}    	
	\centering	
	\includegraphics[width=0.25\linewidth]{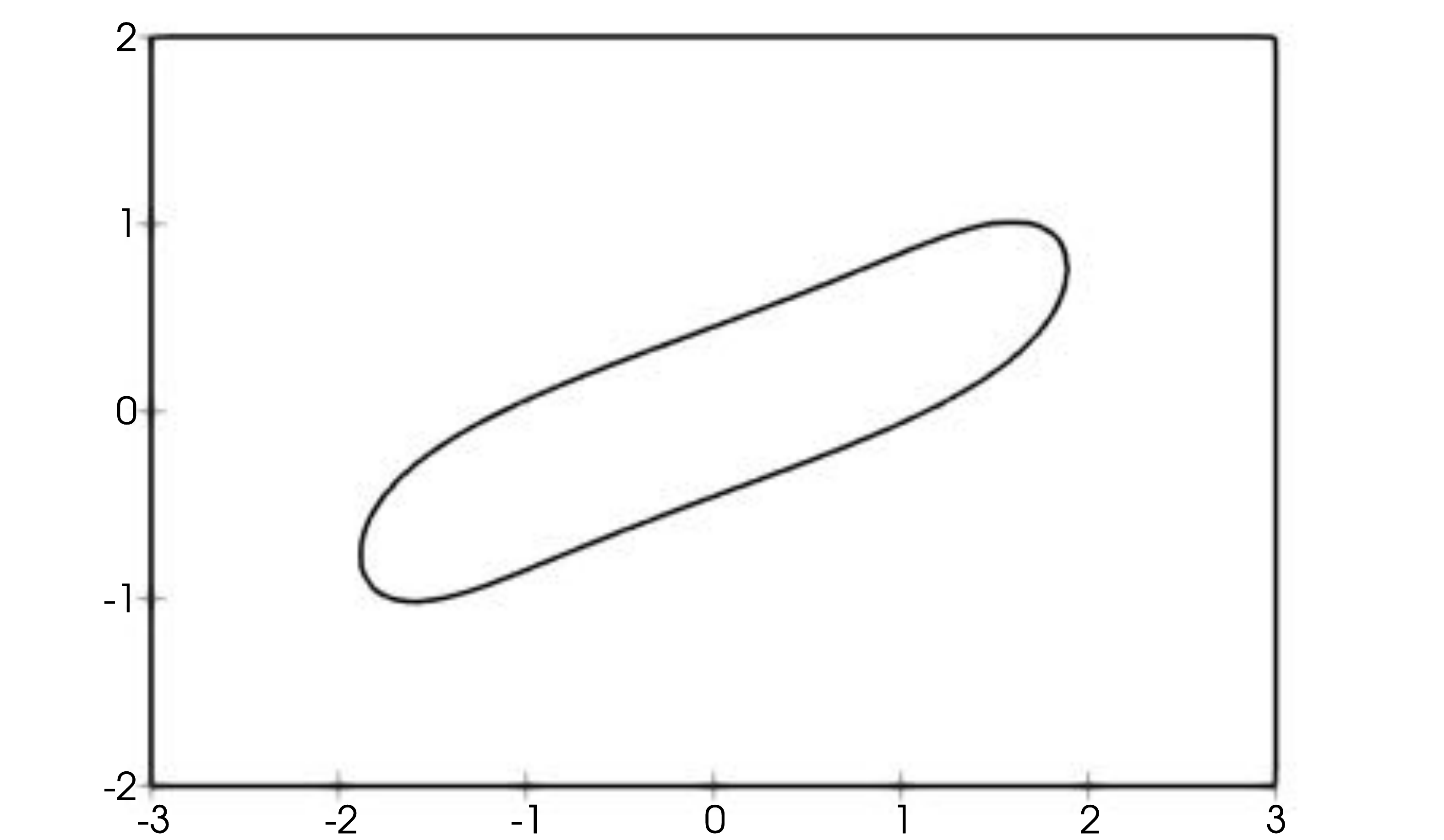}
        }    \\               	
    \subfloat
    	{
    	\scalebox{0.25}    	
	\centering	
	\includegraphics[width=0.25\linewidth]{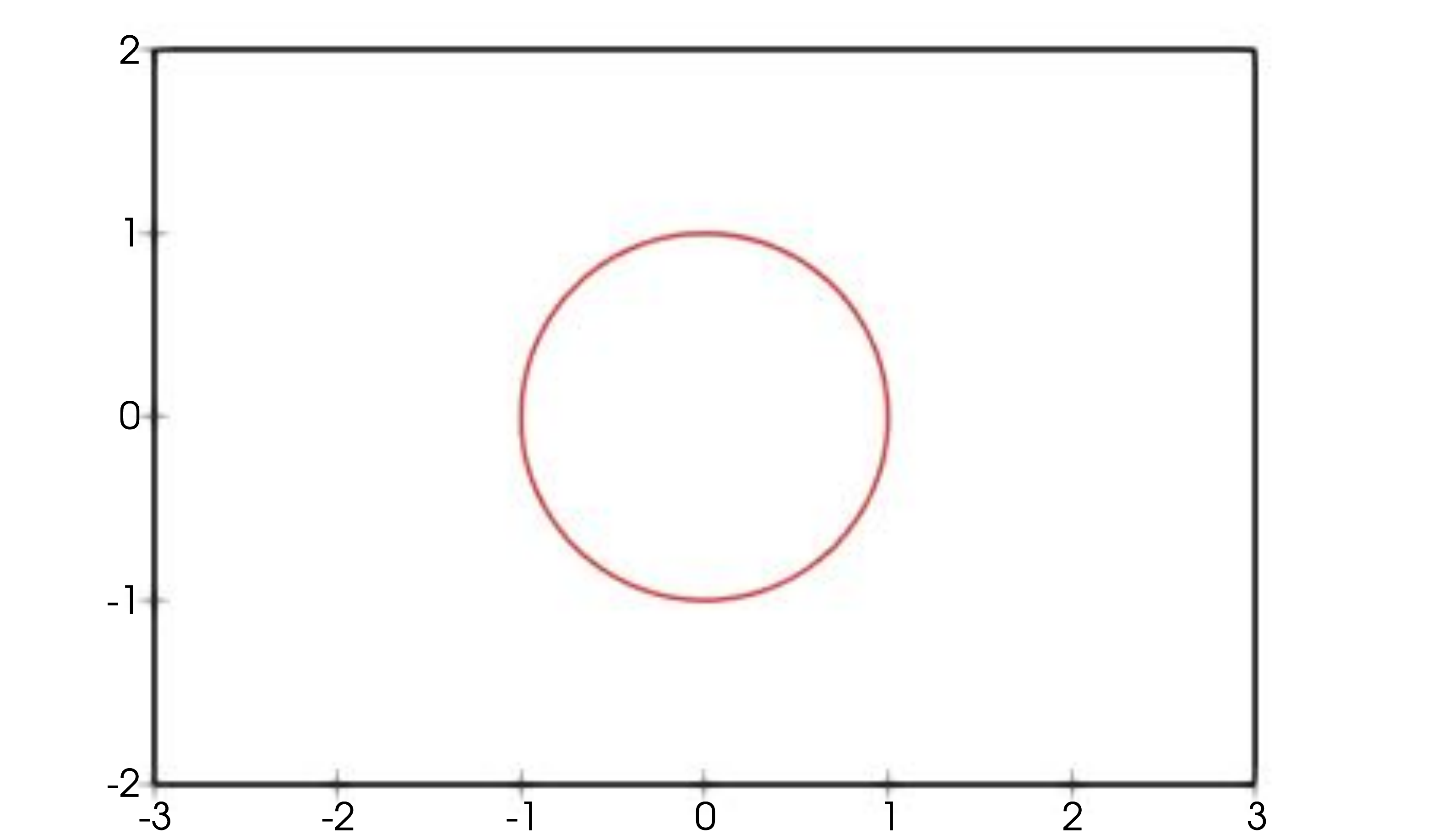}
        }
    \subfloat
    	{
    	\scalebox{0.25}    	
	\centering	
        \includegraphics[width=0.25\linewidth]{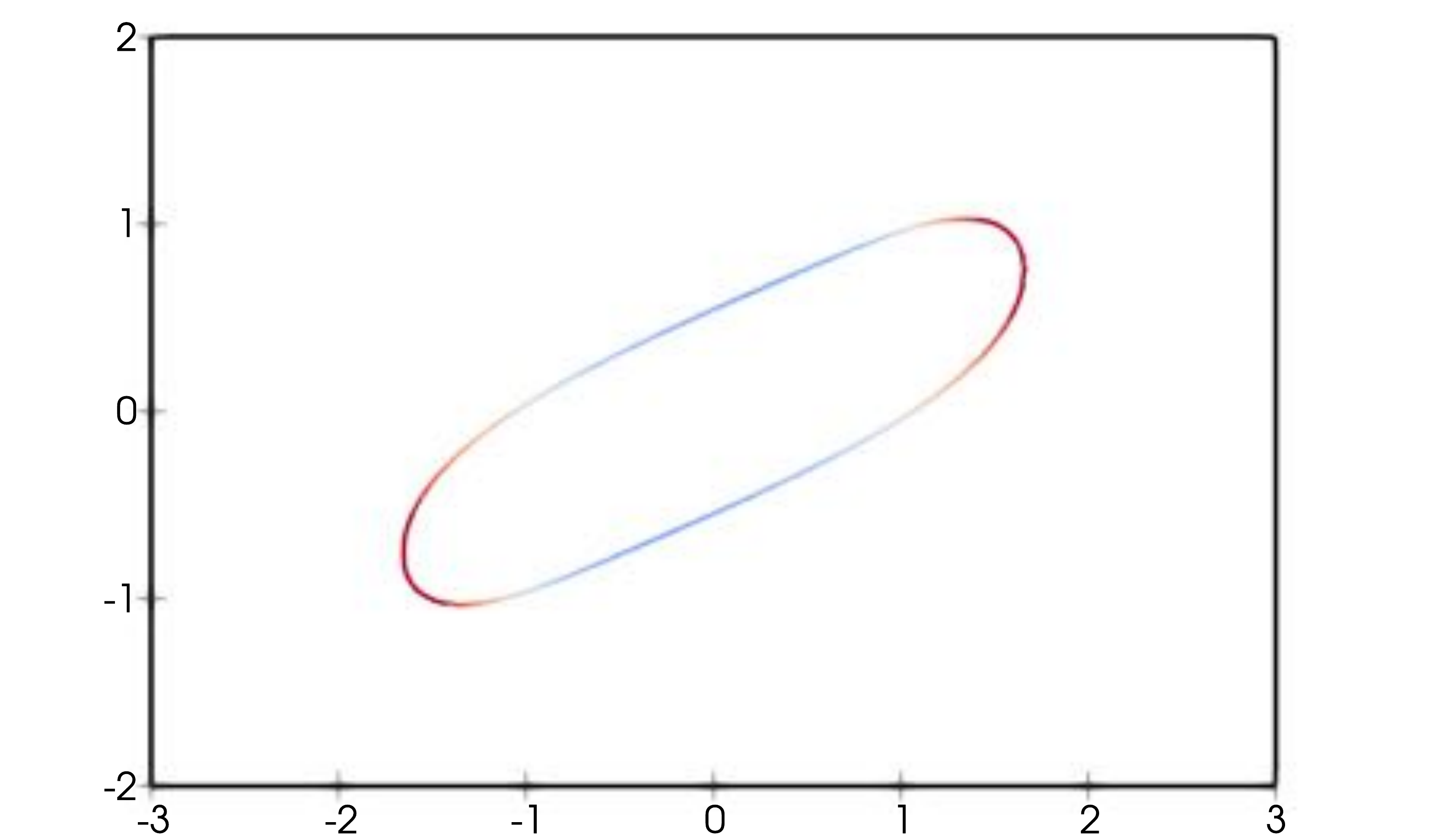}
        }
    \subfloat
    	{
        \scalebox{0.25}    	
	\centering	
	\includegraphics[width=0.25\linewidth]{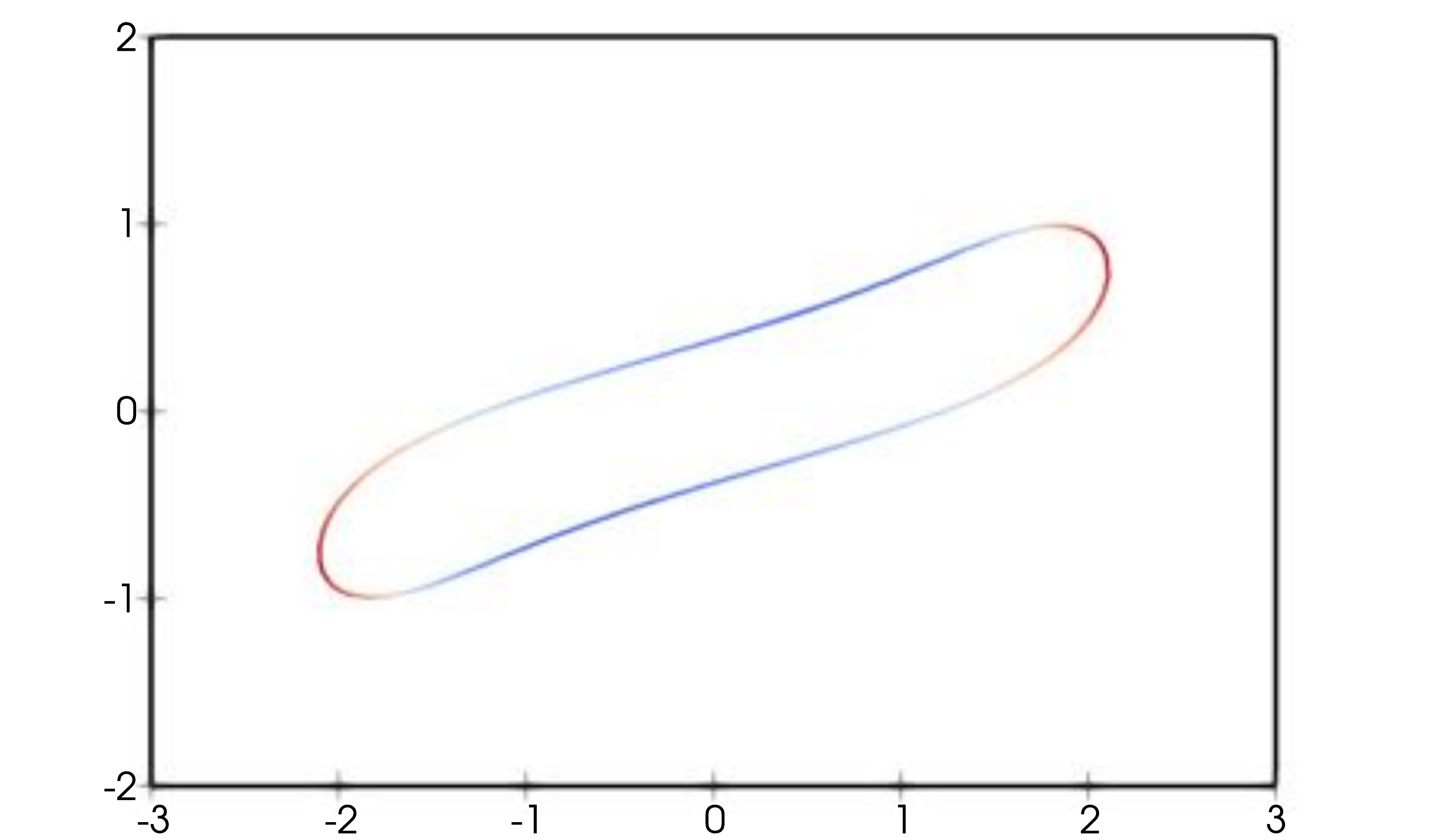}
        }            
     \subfloat
    	{
        \scalebox{0.25}    	
	\centering	
	\includegraphics[width=0.25\linewidth]{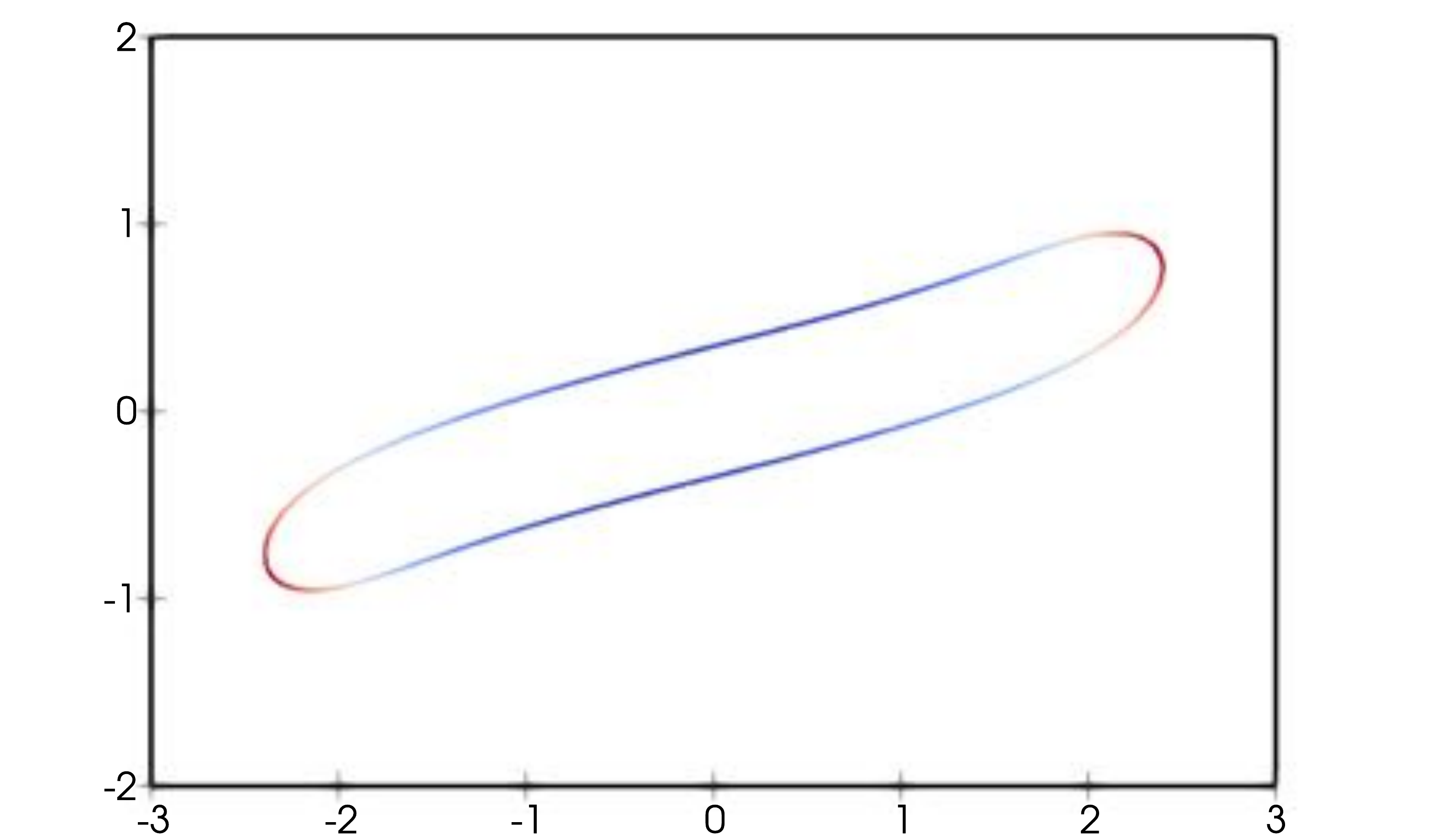}
        }        \\           	
        \setcounter{subfigure}{0}
    \subfloat[$T=0$]
    	{
    	\scalebox{0.25}    	
	\centering	
	\includegraphics[width=0.25\linewidth]{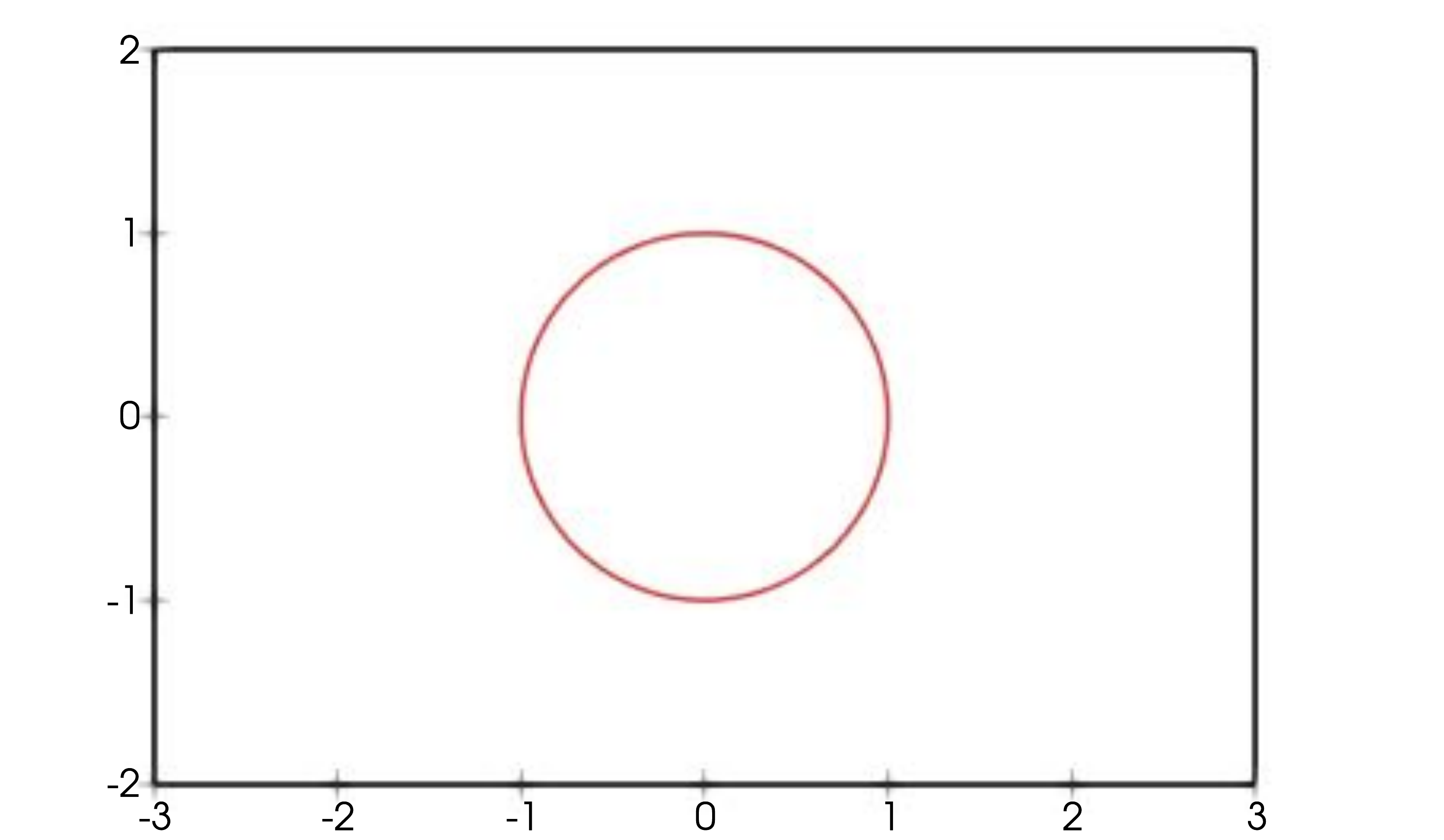}
        }
    \subfloat[$T=4$]
    	{
    	\scalebox{0.25}    	
	\centering	
        \includegraphics[width=0.25\linewidth]{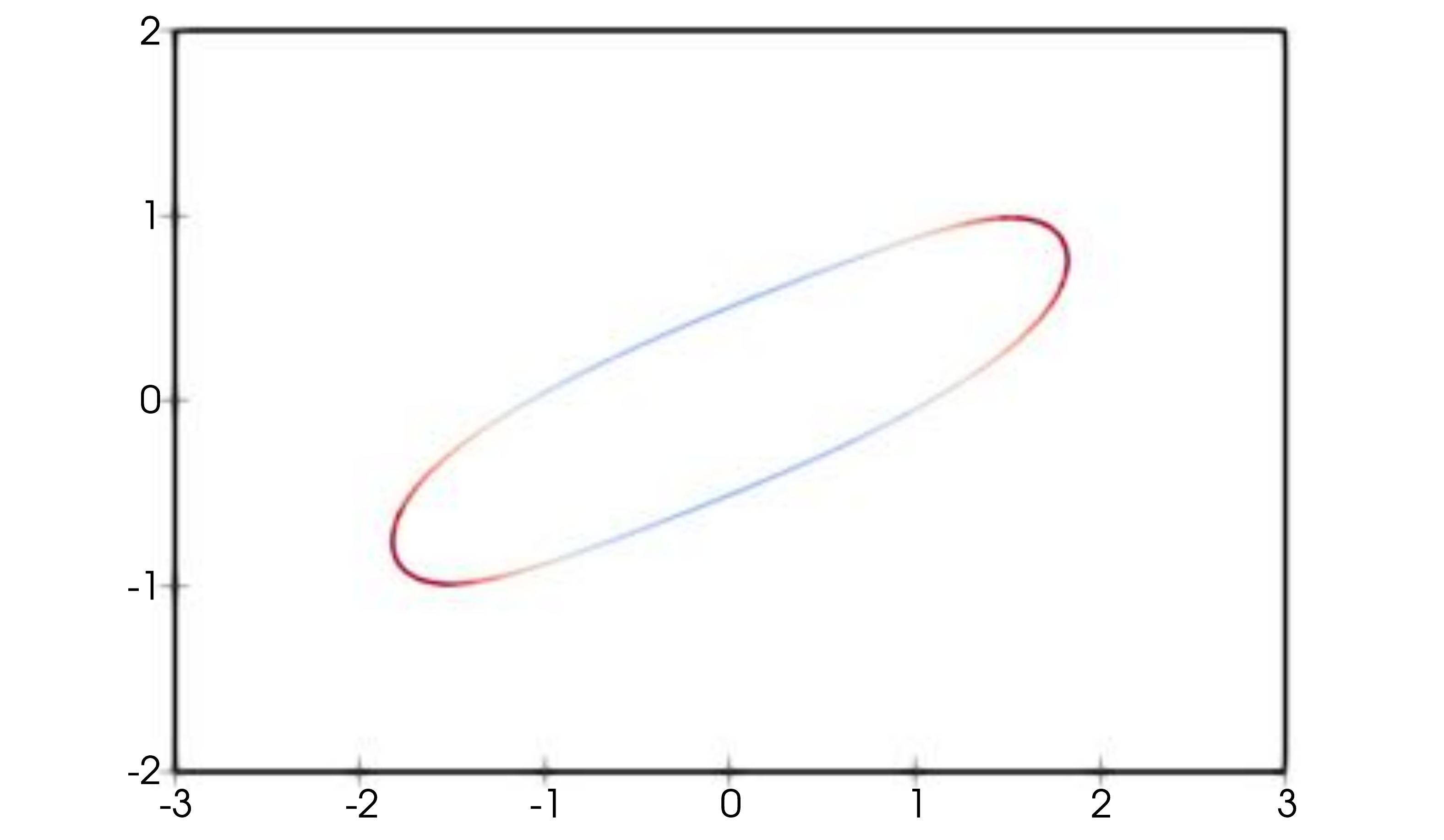}
        }
    \subfloat[$T=8$]
    	{
        \scalebox{0.25}    	
	\centering	
	\includegraphics[width=0.25\linewidth]{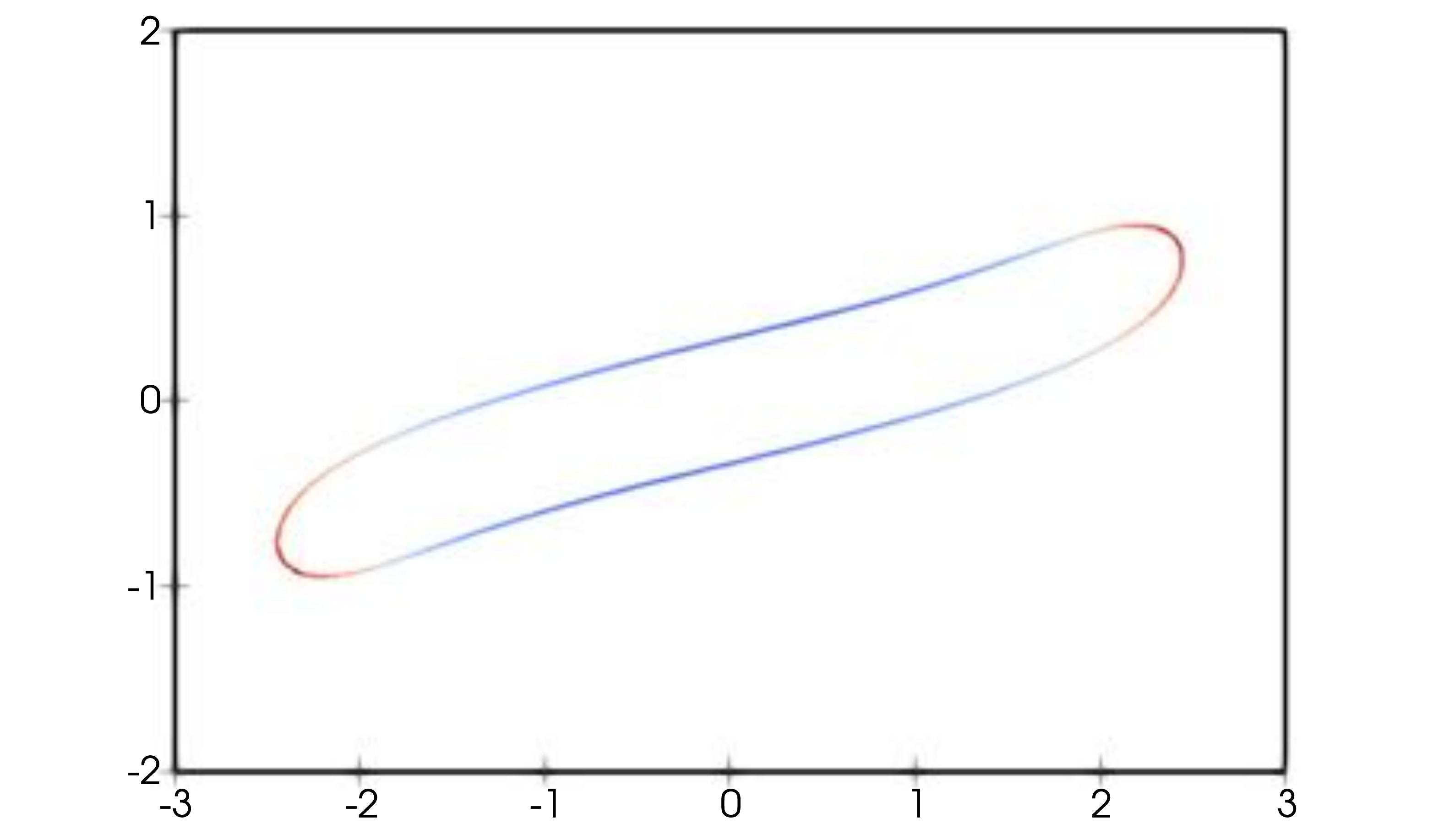}
        }            
     \subfloat[$T=12$]
    	{
        \scalebox{0.25}    	
	\centering	
	\includegraphics[width=0.25\linewidth]{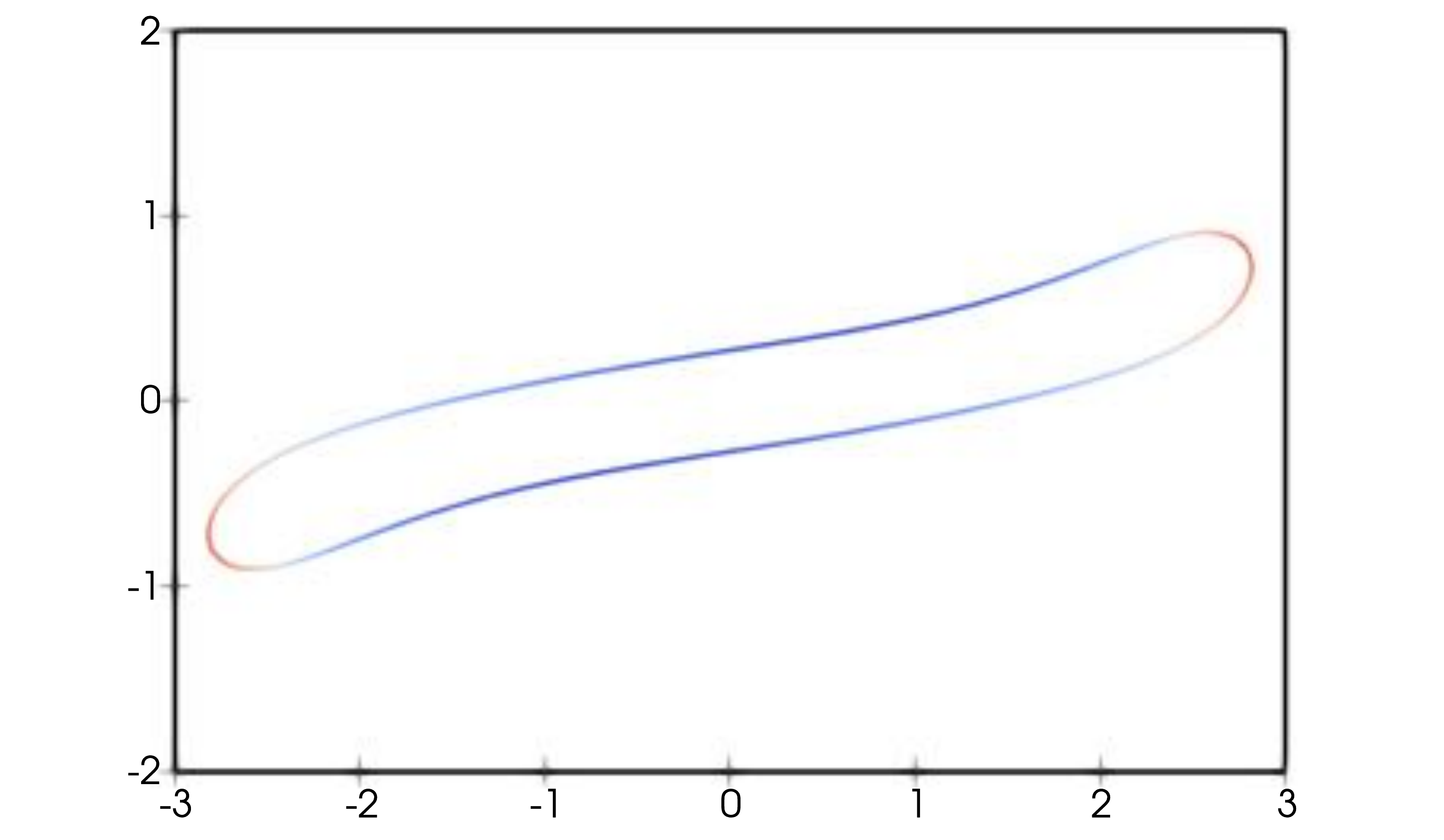}
        }                   \\	
        \centering	
        \vspace{-0.3cm}
     \subfloat
    	{
        \scalebox{0.4}    	
	\centering	
	\includegraphics[width=0.5\linewidth]{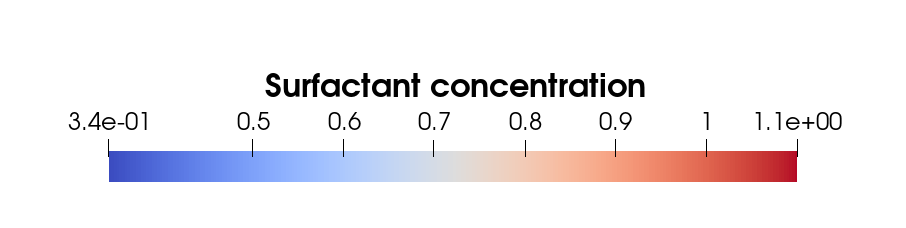}
        }       
        \vspace{-0.3cm}
        \caption{Interface position and surfactant concentration on a drop in a steady shear flow. From top to bottom: results for test case 1-3, ($\beta = 0, 0.25, 0.5$), of Example~\ref{sec:Dropinshear}.}
        \label{fig:shearFlow_T12}
\end{figure}

The simulations in this section have all been done on the mesh shown in Figure \ref{fig:shearMesh}. The coursest elements have size $h_{\text{outer}} = 0.13$ and the smallest element size is $h_{\text{inner}} = 0.05$. The time step has been chosen equal to $0.015$. 
In Figure \ref{fig:shearFlow_T12} the shape of the drop and the surfactant concentration at different time instances, t=0, 4, 8, 12, is shown. One can see that when $\beta$ increases, the drop becomes more elongated and narrower. Our results with the proposed CutFEM agree well with Figure 1  of \cite{LaiTseHua08} where the simulations are made with an immersed boundary method with $h=0.02$ and $\Delta t = h/8=0.0025$ and with Figure 11 of \cite{BaGaNu15}.

\begin{figure}[h!]
     \subfloat[Conservation of surfactant mass]
    	{
        \scalebox{0.25}    	
	\centering	
	\includegraphics[width=0.3\linewidth]{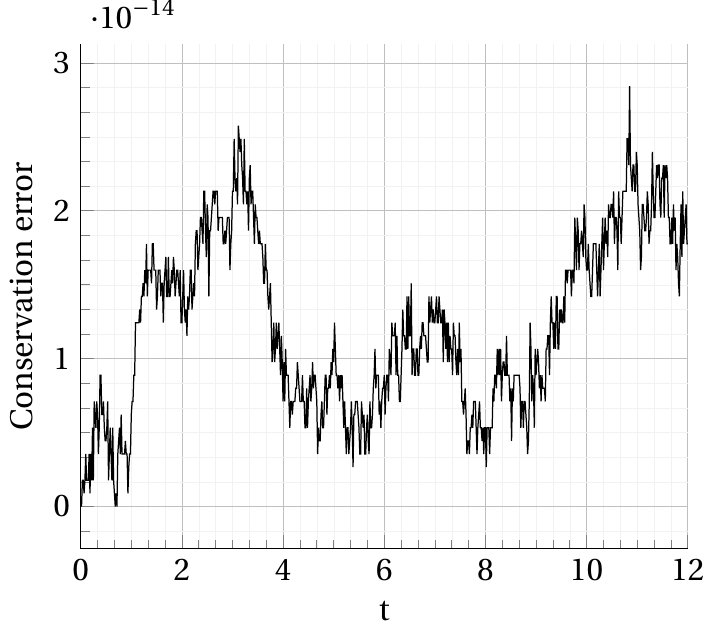}
	\label{fig:shear Evolution quantity - conservation}
        }          \hfill         
     \subfloat[Error in the area of the drop]
    	{
        \scalebox{0.25}    	
	\centering	
	\includegraphics[width=0.3\linewidth]{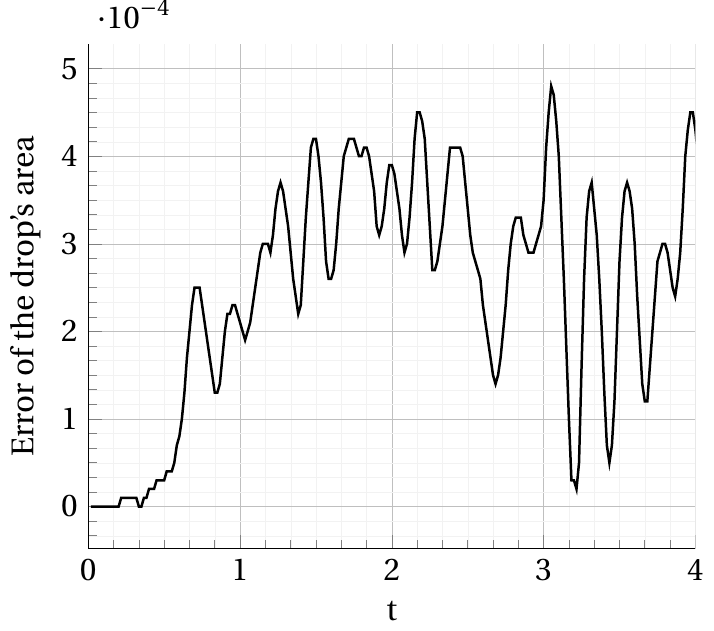}
		\label{fig:shear Evolution quantity - area}
        }       \hfill  
         \subfloat[Total length of the interface]
    	{
        \scalebox{0.25}    	
	\centering	
	\includegraphics[width=0.3\linewidth]{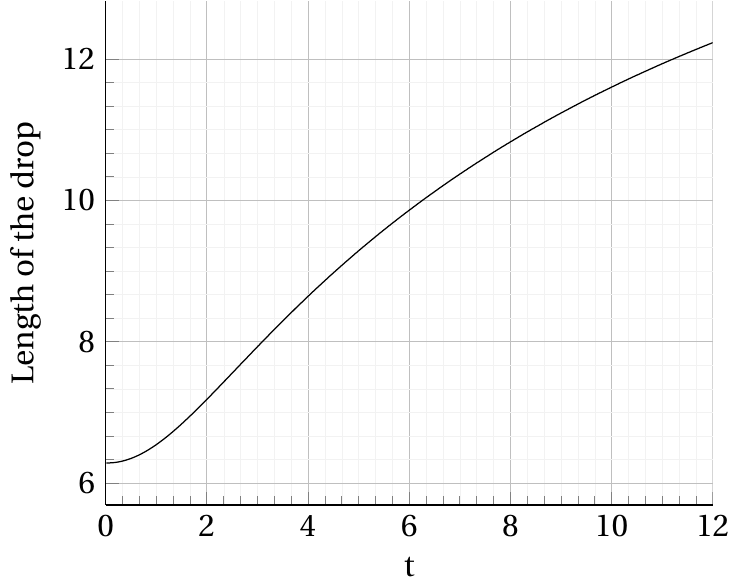}
	\label{fig:shear Evolution quantity - length}
        }    
        \caption{Results for test case 3, i.e., $\beta = 0.5$, in Example~\ref{sec:Dropinshear}}  
        \label{fig:shear Evolution quantity}
\end{figure}

In Figure \ref{fig:shear Evolution quantity} we show different quantities of interest for test case 3 when $\beta=0.5$. First we see in Figure \ref{fig:shear Evolution quantity - conservation} that the total surfactant mass is conserved during the simulation. We show the error in the area of the drop in Figure \ref{fig:shear Evolution quantity - area}. The level set method and the discretization we use does not conserve the area of the drop and for the mesh and time step size that has been used in the simulation the error in the conservation of the area is of order $10^{-4}$. Finally, Figure \ref{fig:shear Evolution quantity - length} shows the perimeter of the drop. These results can be compared with Figure 5 in \cite{LaiTseHua08}.

\newpage

\subsection{Pair of drops in shear flow}\label{sec:dropdrop}
In this example we study drop-drop interaction. Two circular drops with radius $0.5$ are intitially centered at $(-1,0.25)$ and $(1,-0.25)$. The computational domain is $\Omega=[-4,4]\times[-1,1]$ and Dirichlet boundary condition $\bfu(t,\bfx) = (y, 0)$ is imposed. The initial velocity is $\bfu(0,\bfx) = (y,0)$, $\bff=0$, and we consider two test cases with physical parameters chosen as in Table \ref{tab:set parameter double drops}.
\begin{table}[ht]
\centering
\begin{tabular}{ p{2cm}  p{1cm} p{1cm} p{1cm} p{1cm} p{1cm} p{1cm} p{1cm} p{1cm}p{2cm}}
\hline 
Test case   & $\rho_1$ & $\rho_2$ & $\mu_1$ & $\mu_2$ & $\sigma_0$ & $\beta$ & $w_0$ & $D_{\Gamma}$\\ 
\hline
1		& 1	&	2	&	0.05	&	0.1	&	0.25    & 0.0       &  1   &  2.5	\\
2		& 1	&	2	&	0.05	&	0.1	&	0.25 	   & 0.6      &  1   &  2.5 	\\
\hline 
\end{tabular}
\captionof{table}{Parameters used in Example~\ref{sec:dropdrop}. }
\label{tab:set parameter double drops}
\end{table}

\begin{figure}[h!]
	\hspace{-8cm}
 	\includegraphics[width=2\linewidth]{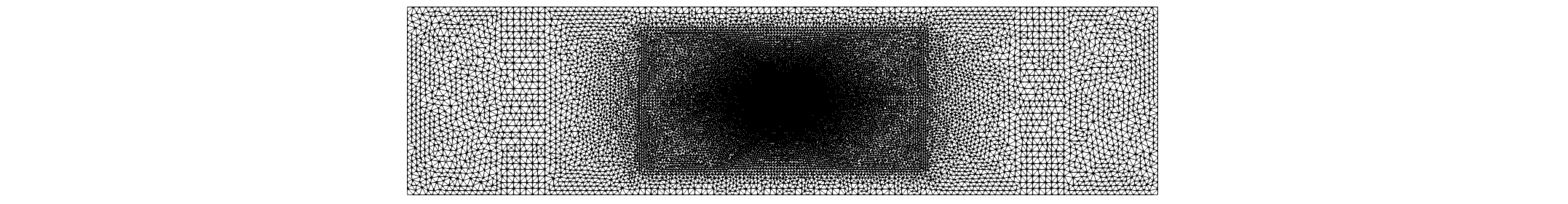}
	 \caption{Mesh used for the simulations of Example~\ref{sec:dropdrop}. The mesh contains 24328 nodes and 48354 elements.}  
	 \label{fig:meshDoubleDrop}
\end{figure}

\begin{figure}[h!]
\centering
    \subfloat[$T=0$]
    	{
    	\scalebox{0.3}    	
	\centering	
	\includegraphics[width=0.25\linewidth]{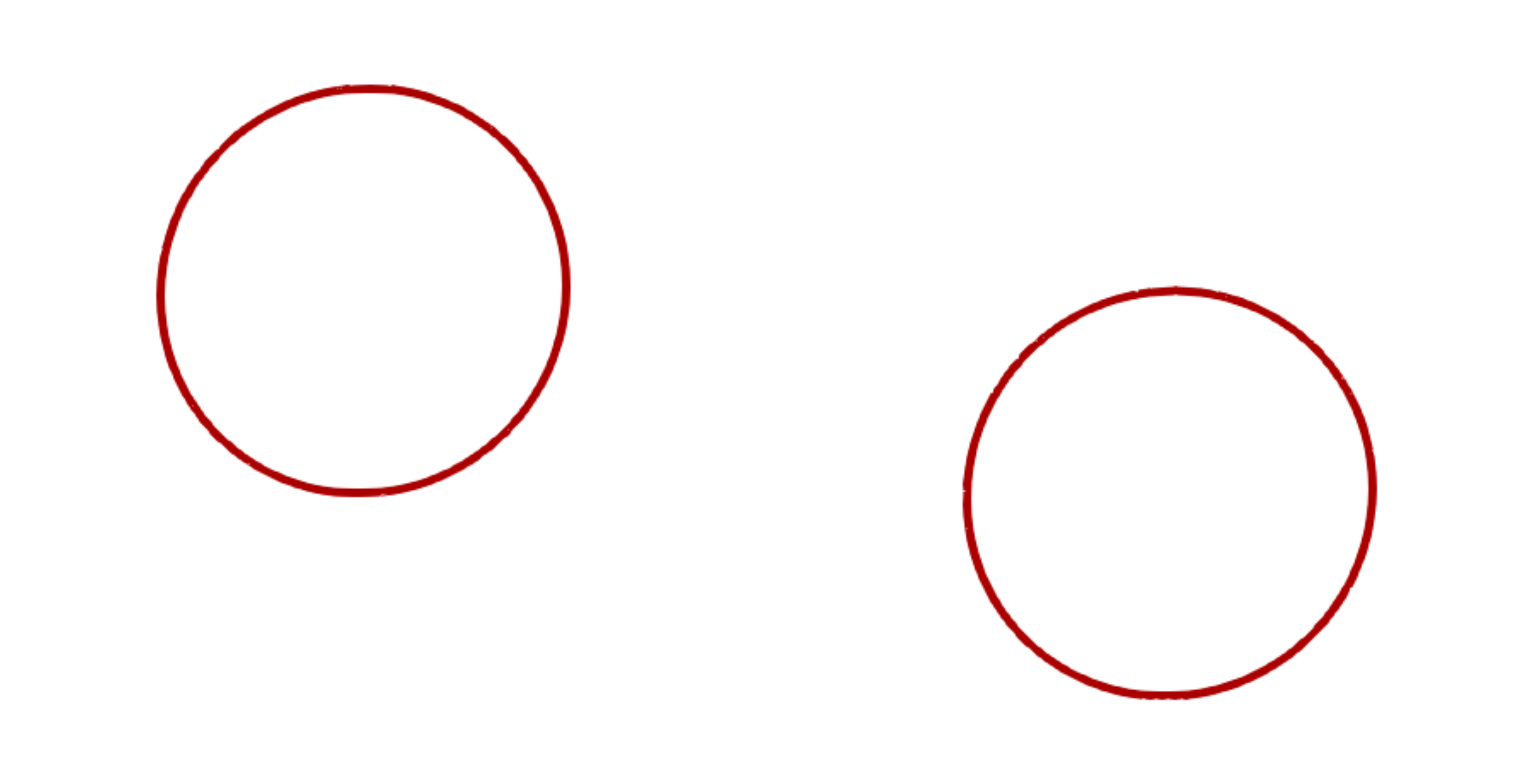}
        }
        \hfill
    \subfloat[$T=2$]
    	{
    	\scalebox{0.3}    	
	\centering	
	\includegraphics[width=0.25\linewidth]{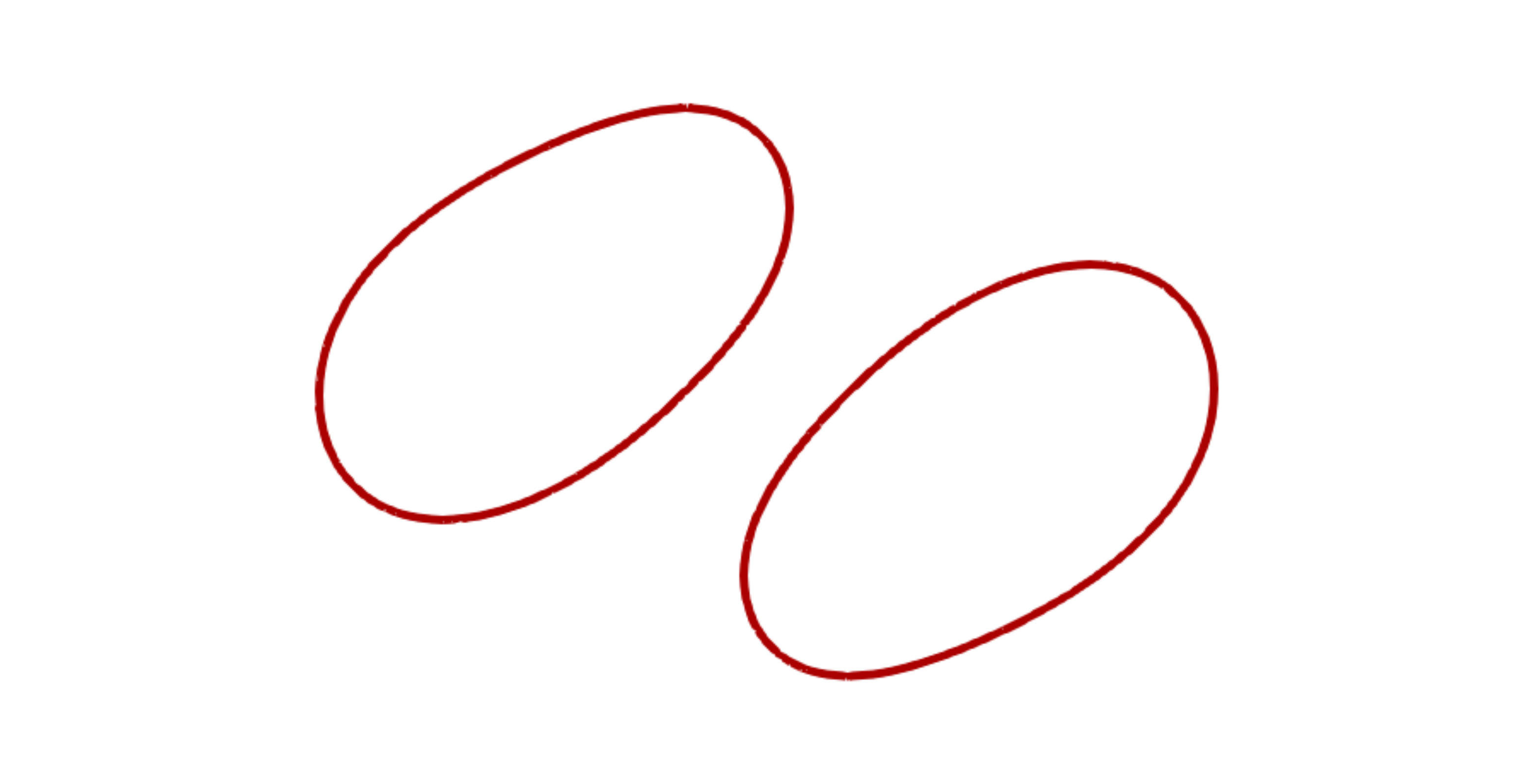}
        }
         \hfill
    \subfloat[$T=4$]
    	{
        \scalebox{0.3}    	
	\centering	
	\includegraphics[width=0.25\linewidth]{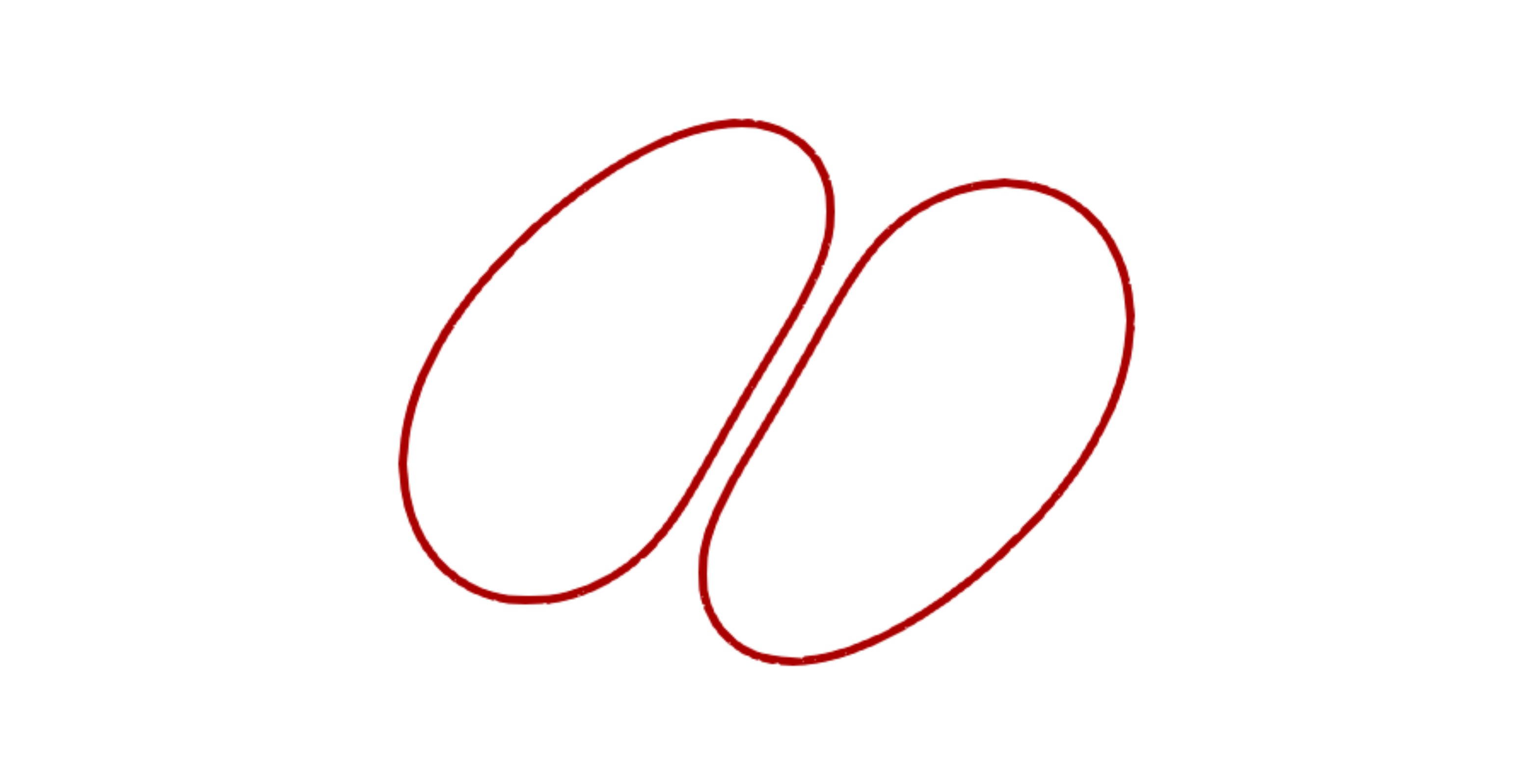}
        }       \\ 
     \subfloat[$T=6$]
    	{
        \scalebox{0.3}    	
	\centering	
	\includegraphics[width=0.25\linewidth]{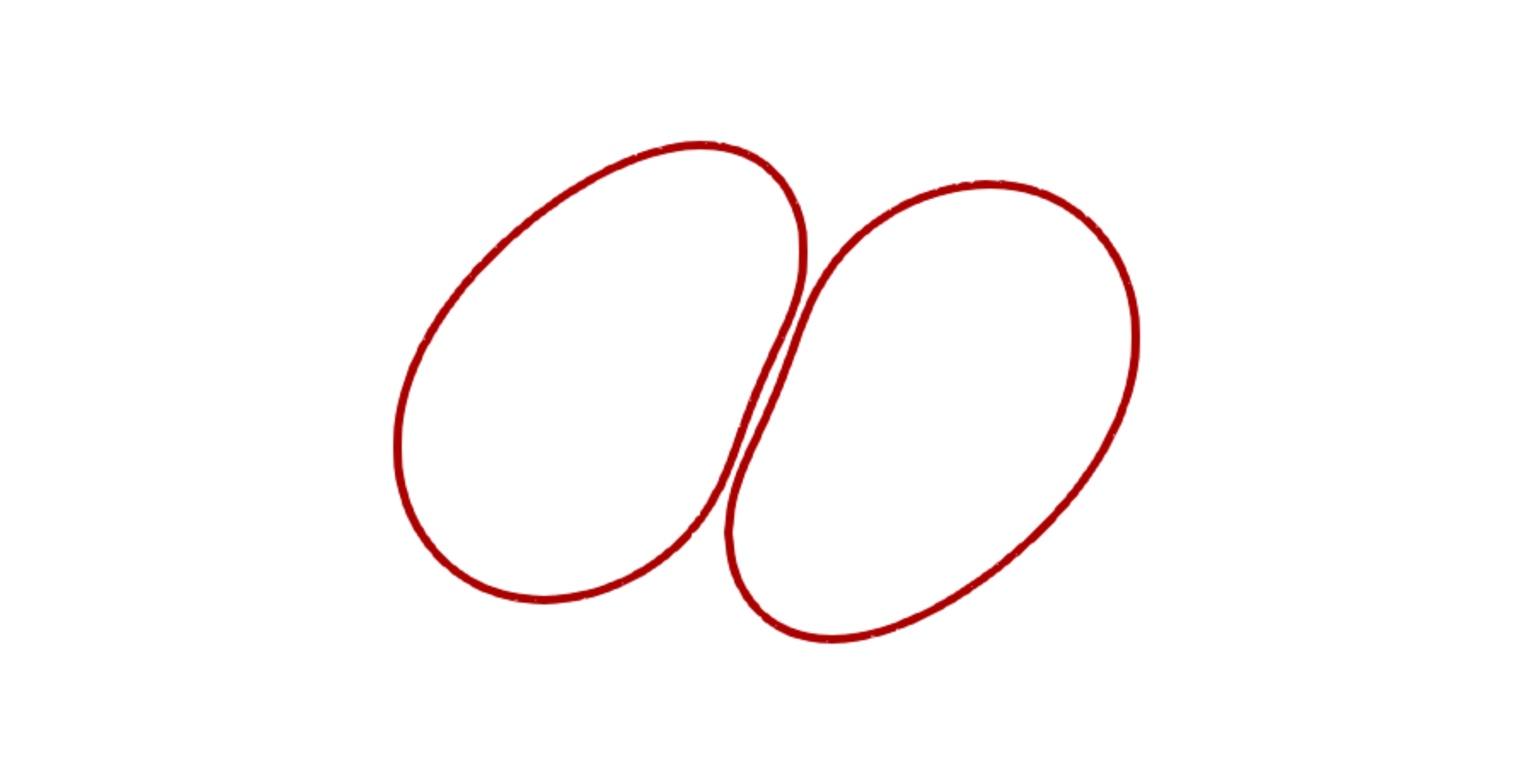}
        }      
         \hfill
         \subfloat[$T=8$]
    	{
        \scalebox{0.3}    	
	\centering	
	\includegraphics[width=0.25\linewidth]{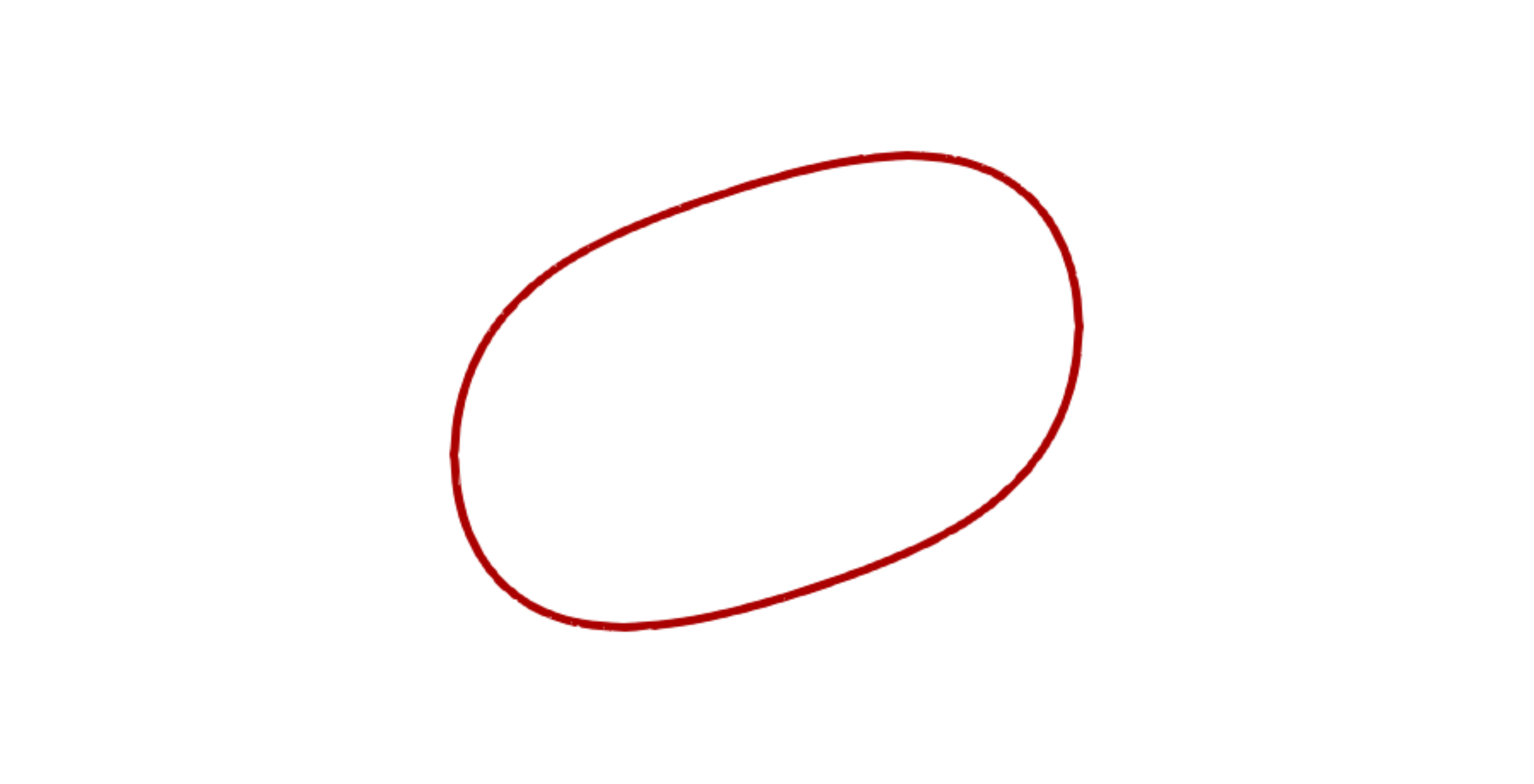}
        }      
         \hfill
        \subfloat[$T=10$]
    	{
        \scalebox{0.3}    	
	\centering	
	\includegraphics[width=0.25\linewidth]{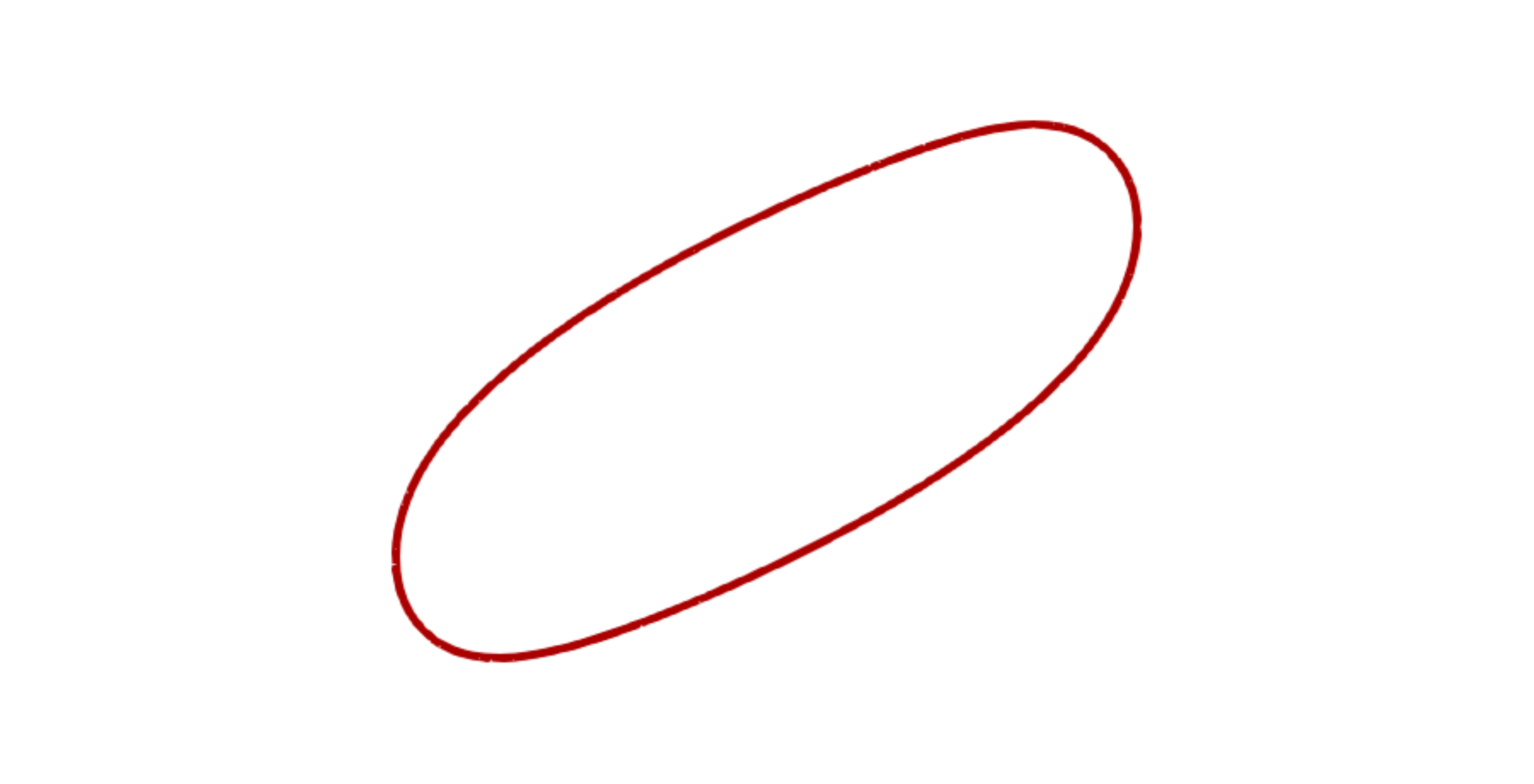}
        }           \\
        \vspace{-0.3cm}
        \caption{Test case 1 in Example~\ref{sec:dropdrop}: clean interfaces, the two drops merge and form a single drop. }
        \label{fig:2drops test case 1}
\end{figure}

\begin{figure}[h!]
\centering
    \subfloat[$T=0$]
    	{
    	\scalebox{0.3}    	
	\centering	
	\includegraphics[width=0.25\linewidth]{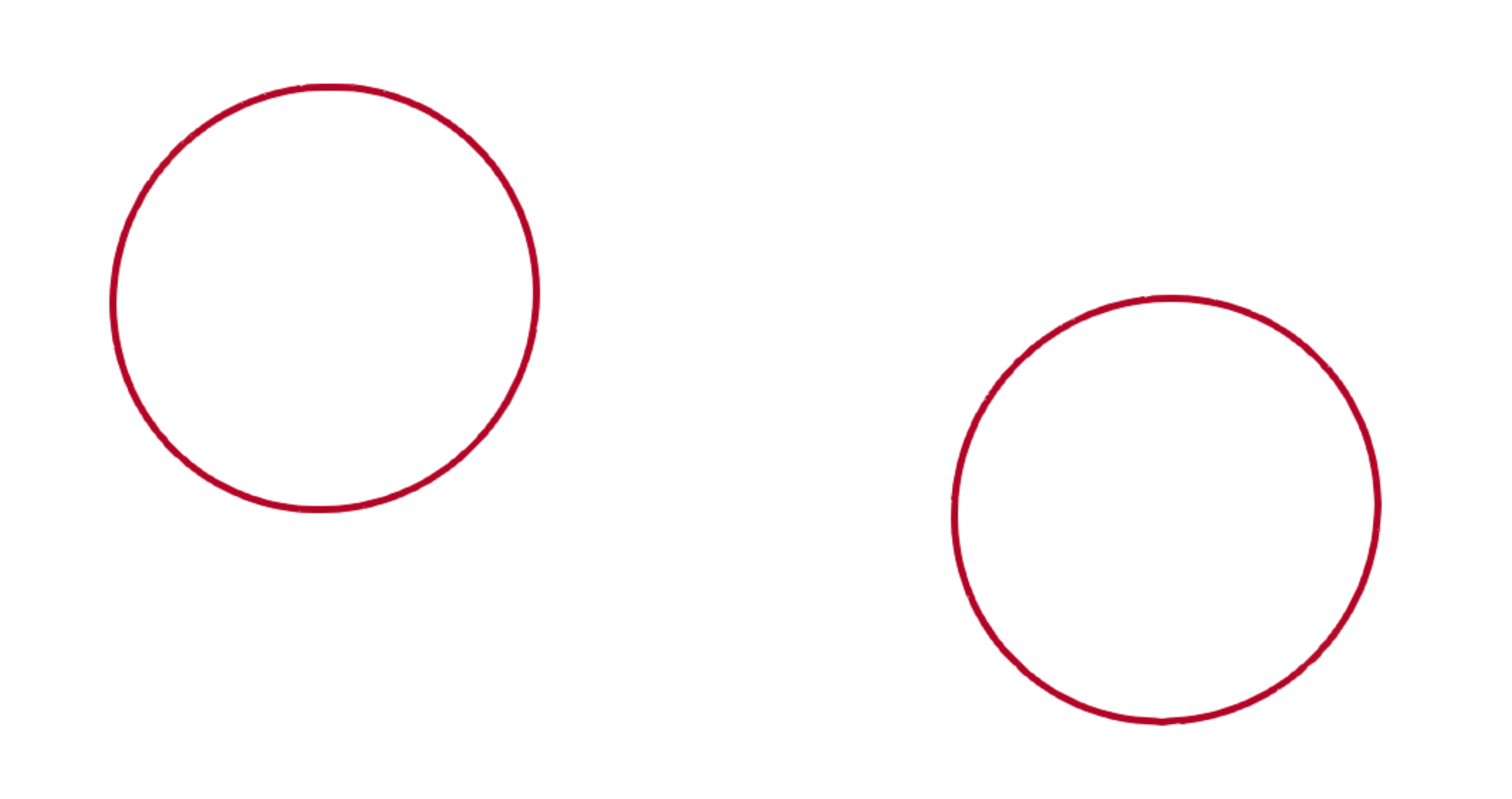}
        }
        \hfill
    \subfloat[$T=2$]
    	{
    	\scalebox{0.3}    	
	\centering	
	\includegraphics[width=0.25\linewidth]{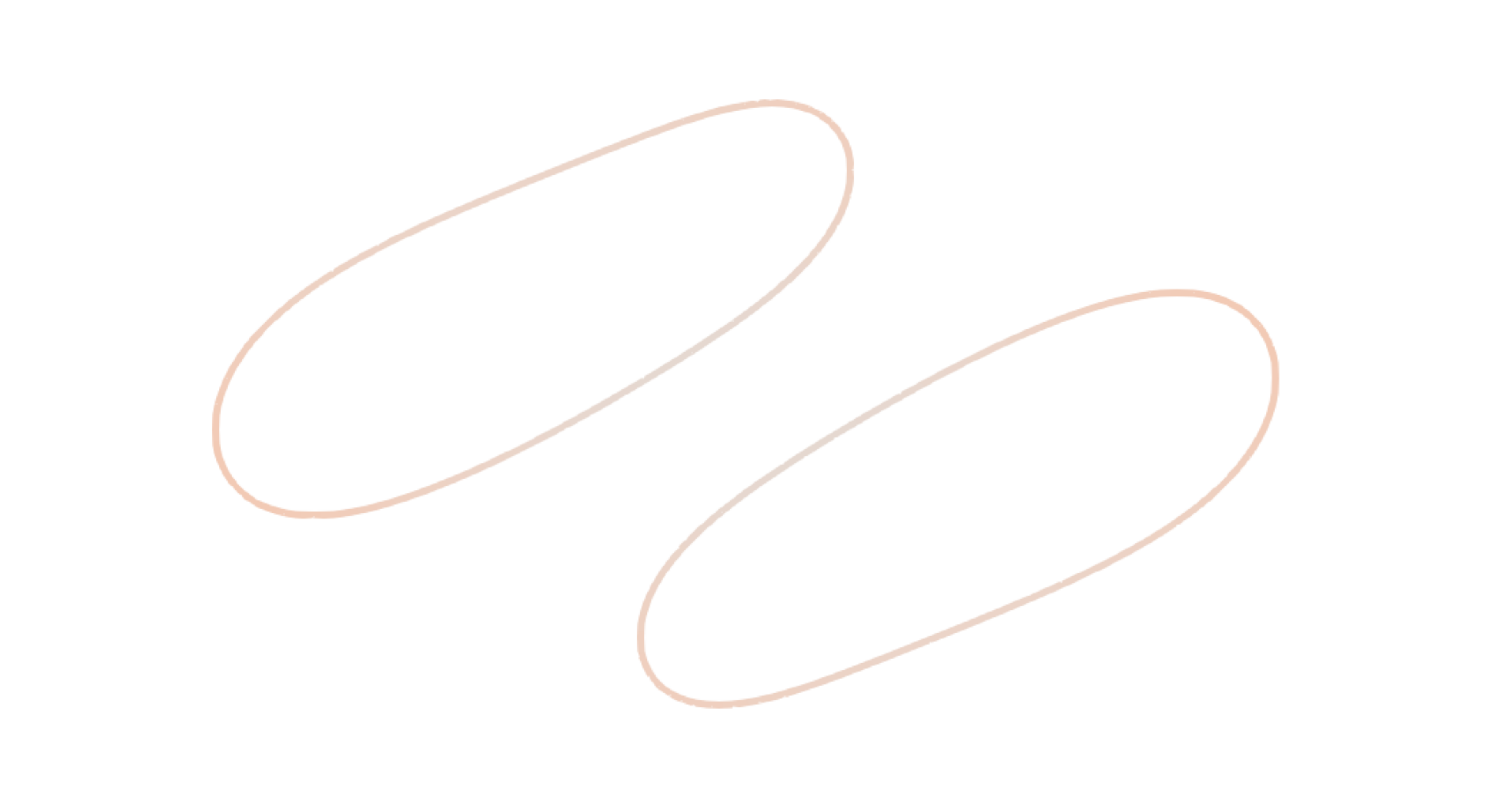}
        }
         \hfill
    \subfloat[$T=4$]
    	{
        \scalebox{0.3}    	
	\centering	
	\includegraphics[width=0.25\linewidth]{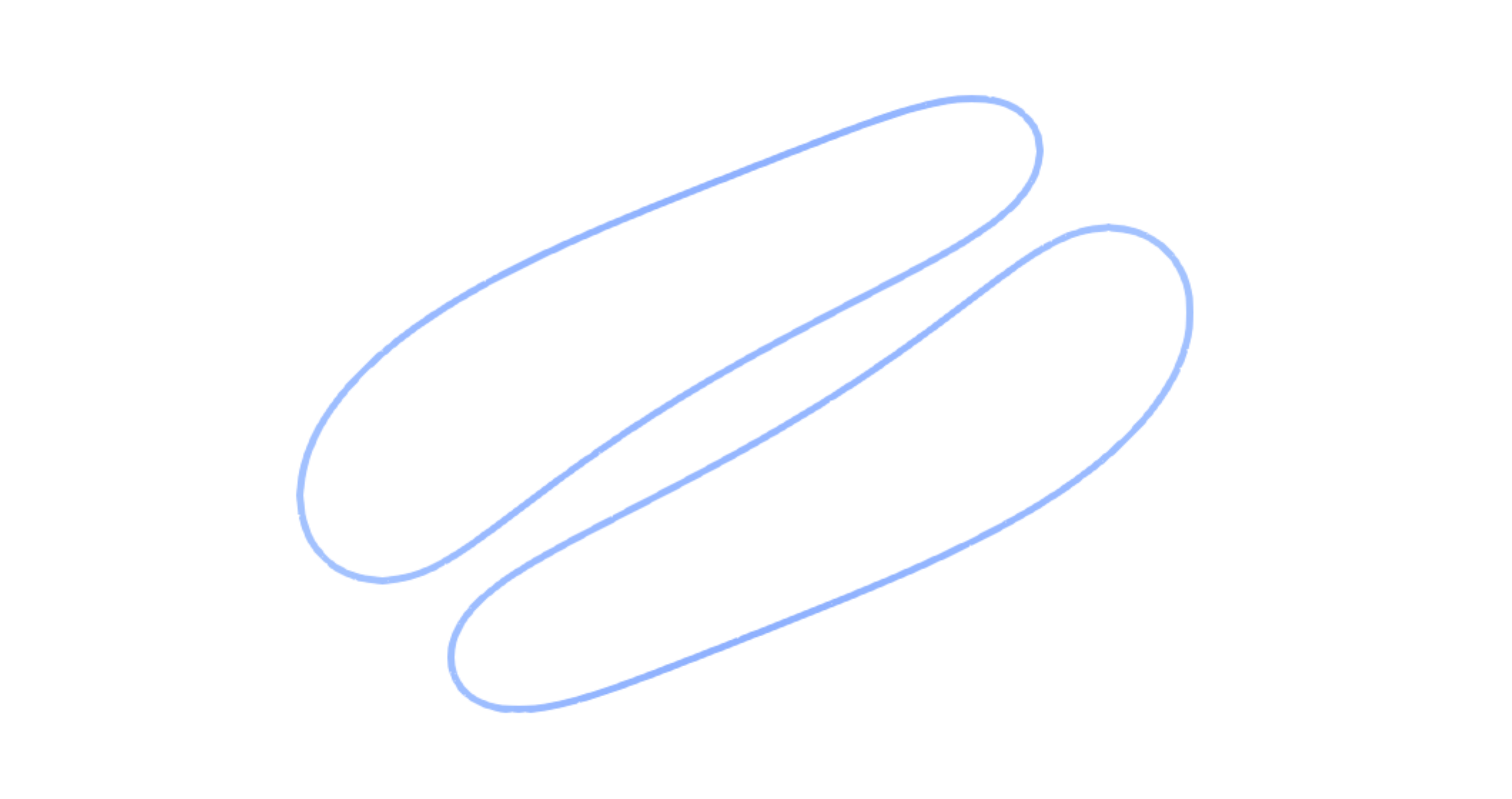}
        }       \\ 
     \subfloat[$T=6$]
    	{
        \scalebox{0.3}    	
	\centering	
	\includegraphics[width=0.25\linewidth]{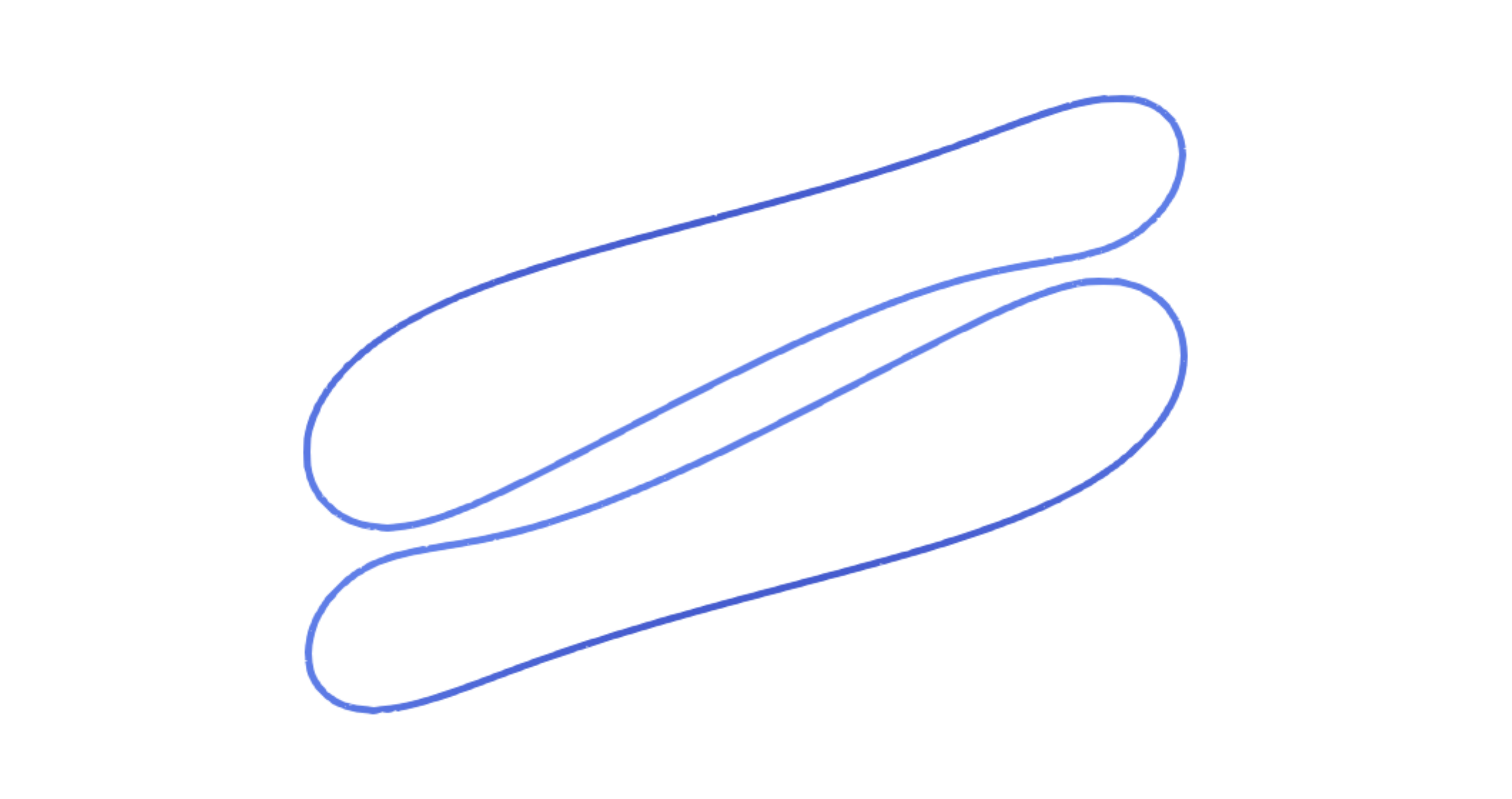}
        }      
         \hfill
         \subfloat[$T=8$]
    	{
        \scalebox{0.3}    	
	\centering	
	\includegraphics[width=0.25\linewidth]{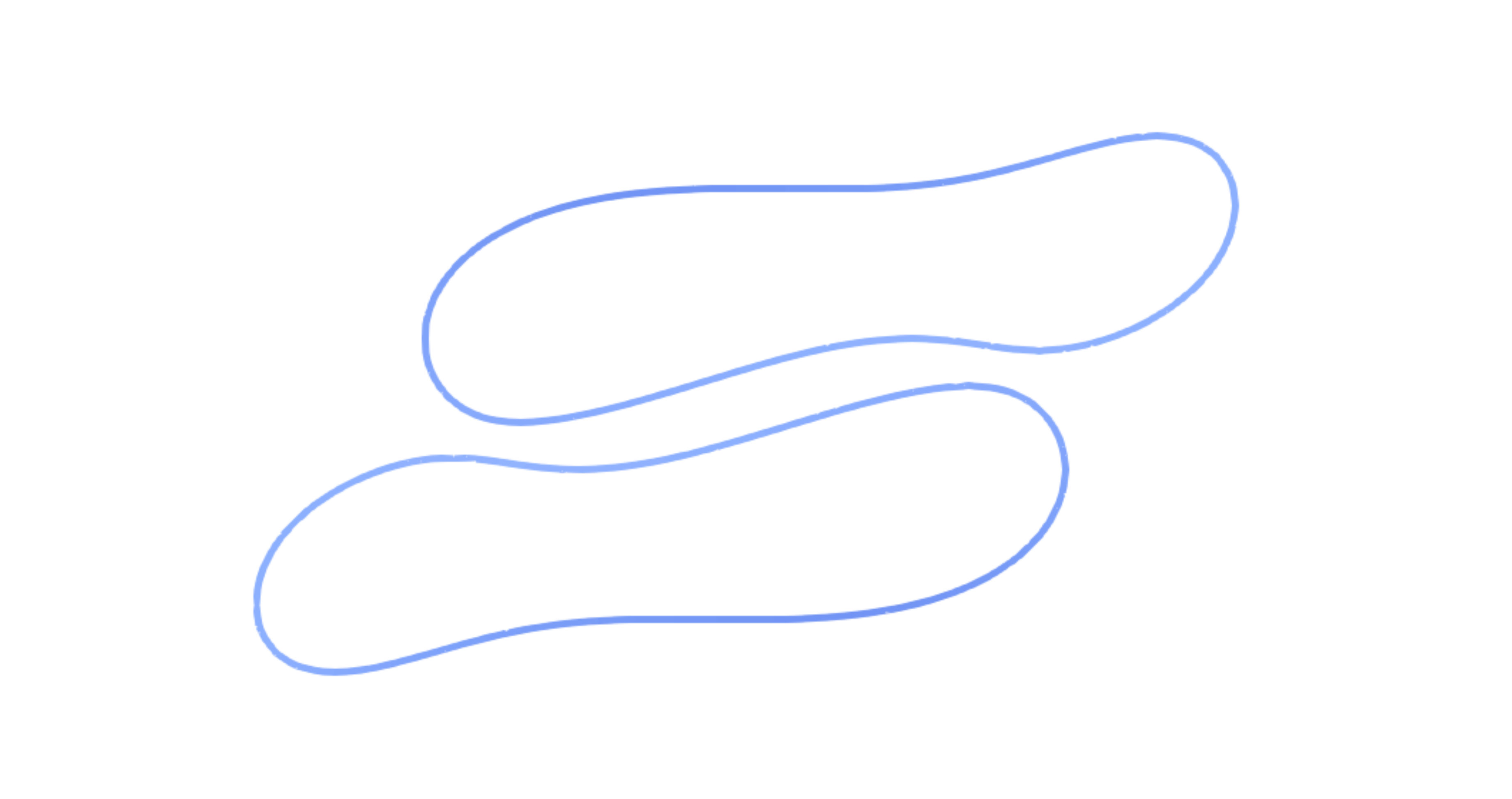}
        }    
         \hfill  
        \subfloat[$T=10$]
    	{
        \scalebox{0.3}    	
	\centering	
	\includegraphics[width=0.25\linewidth]{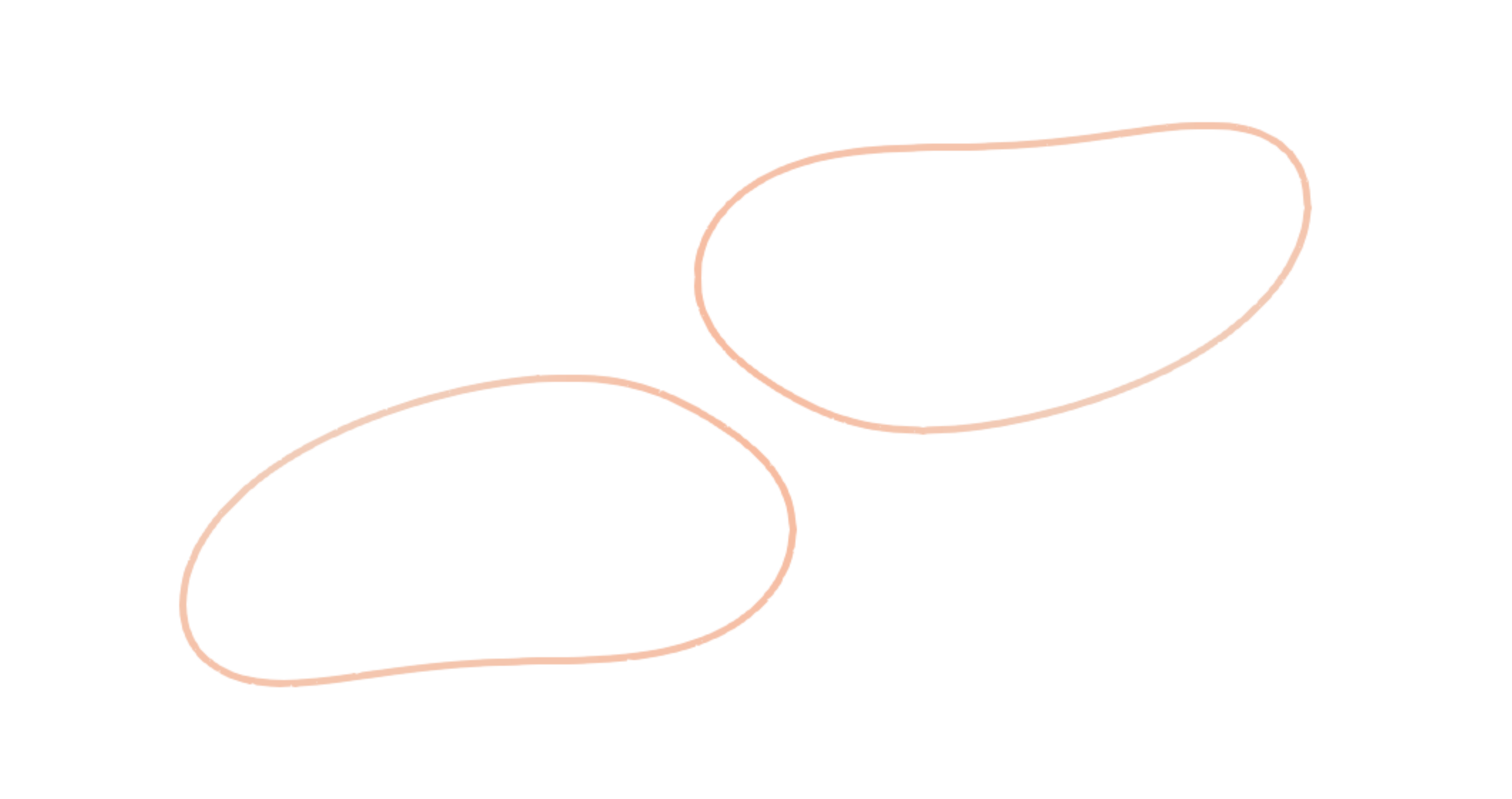}
        }           \\
        \vspace{-0.3cm}
        \subfloat
    	{
        \scalebox{0.75}    	
	\centering	
	\includegraphics[width=0.4\linewidth]{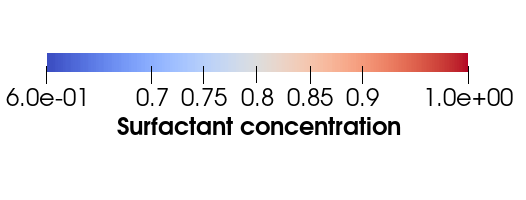}
        }         	
        \vspace{-0.5cm}
        \caption{Test case 2 in Example~\ref{sec:dropdrop}:  surfactant is added, no coalescence occur.}
        \label{fig:2drops test case 2}
\end{figure}

The background mesh is shown in Figure \ref{fig:meshDoubleDrop}, where close to the boundaries, $\partial \Omega$, the mesh size is $1/25$ and around the origin the finest elements are of size $1/100$. 
In Figure \ref{fig:2drops test case 1} and  \ref{fig:2drops test case 2}  we show  the drop shape in the two test cases. In the first case, with clean interfaces, the two drops coalesce during contact and form a single drop but in the second case, in the presence of surfactant ($\beta = 0.6$),  the two drops become more elongated and pass each other without merging.
Note that here topological changes such as merging is easily handeled by the level set method.

\newpage

\subsection{Rising drop in 3D}\label{sec:risingdrop3D}
In this section, we consider the first 3D example in \cite{BaGaNu15}. The interface is initialy a sphere with radius $0.25$ centered at $(0.5,0.5,0.5)$. The computational domain is $\Omega = [0, \ 1]\times [0, \ 1] \times [0, \ 2]$ and $I=[0,T]$ with $T=3$. The physical parameters are chosen as in Table \ref{table - set of parameters}.

\begin{figure}[h!]
\centering          
        \setcounter{subfigure}{0}   
         \subfloat[$\beta = 0$]
    	{
        \scalebox{0.25}    	
	\centering	
	\includegraphics[width=0.25\linewidth]{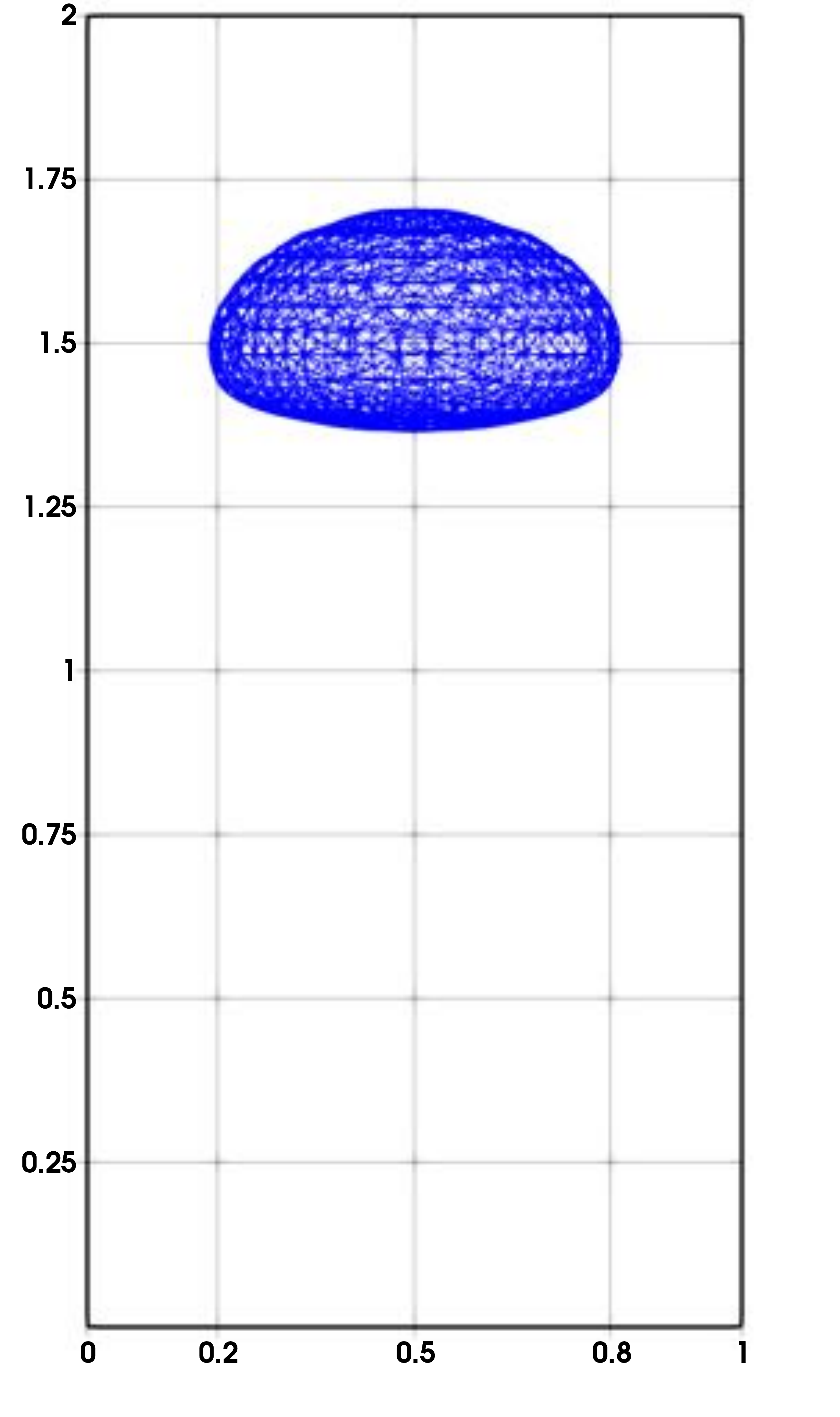}
        }     
         \subfloat[$\beta = 0.5$]
    	{
        \scalebox{0.25}    	
	\centering	
	\includegraphics[width=0.25\linewidth]{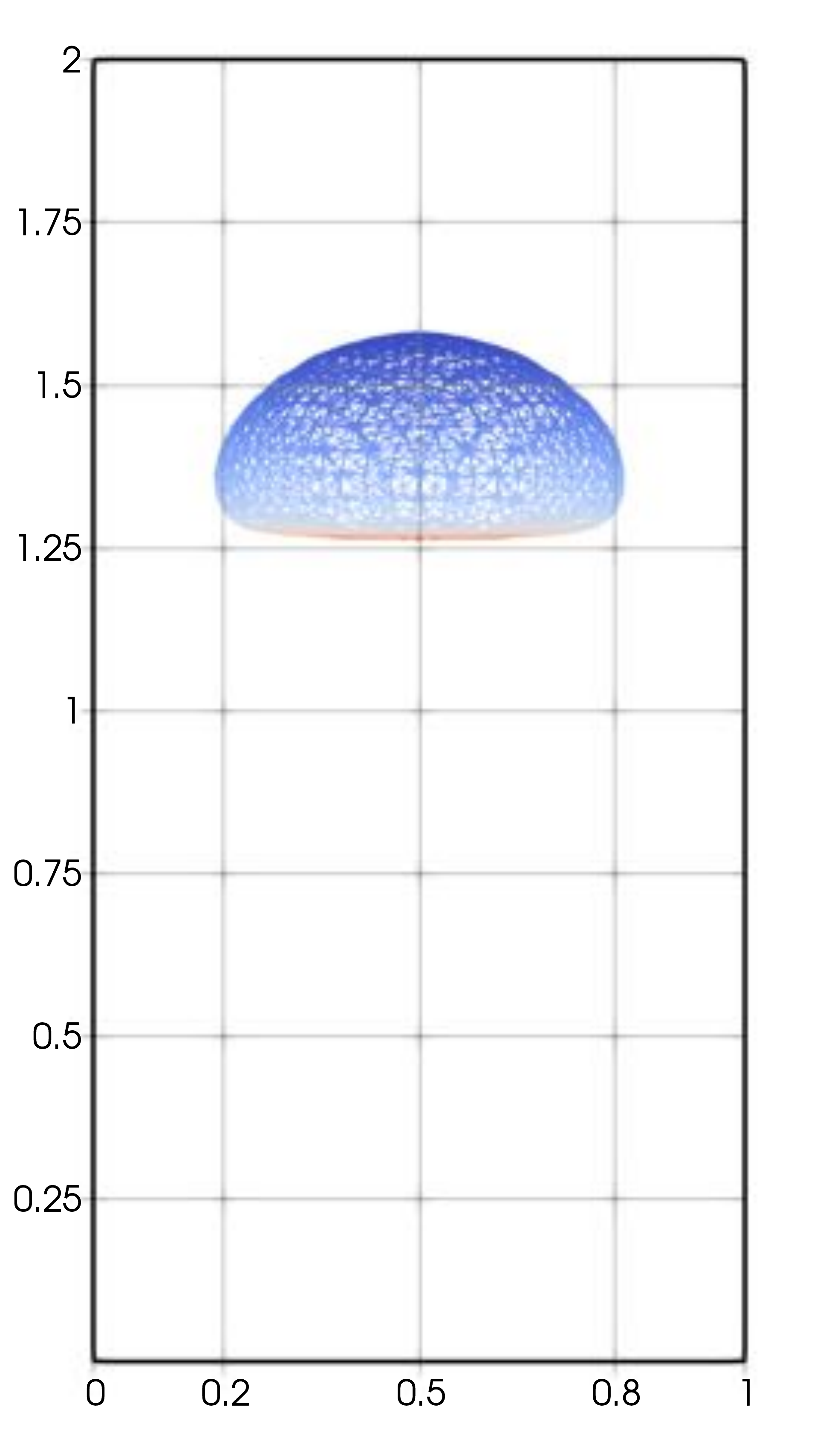}
        }     
         \subfloat
    	{
        \scalebox{0.25}    	
	\centering	
	\includegraphics[scale =0.25]{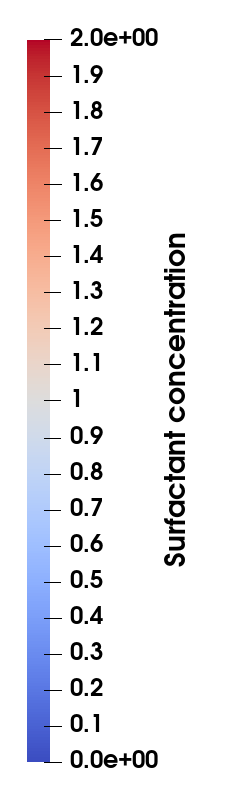}
        }    
        \hfill
        \caption{Shape of the drop at the final time $t = 3$ for Example \ref{sec:risingdrop3D}. } \label{fig:3D_finalShape}
\end{figure}

We compute the same quantities as in \ref{example rising drop 2D} but with their natural extensions to three space dimensions. We define the circularity as 
\begin{align*}
c = \frac{S_a}{S_b},
\end{align*}
where $S_a$ is the surface area of the sphere which has the same volume as the drop and $S_b$ is the surface area of the drop.\\
The computations are done on a uniform mesh of size $h=1/30$ with a time step $\Delta t = 0.02$. 
One can see the domain and the final shape of the drops in Figure \ref{fig:3D_finalShape}. As in two space dimensions the presence of surfactant increase the deformation of the drop and decrease the rise velocity. 
Table \ref{table - result rising drop 3D} shows some characteristic quantities obtained for the simulation with and without surfactant. The results agree with the results obtained in \cite{BaGaNu15}. 

\begin{table}[ht]
\centering
\begin{tabular}{ p{2cm} p{3cm} p{3cm}  }
\hline
 & without surfactant &  with surfactant\\
\hline
$c_\text{min}$	& 0.95130 & 0.9279	\\
$t_{c_\text{min}}$	& 3 & 2.8535	\\
$u^3_{c, \text{max}}$ & 0.3816	& 0.3147	\\
$t_{u^3_{c, \text{max}}}$ & 0.9438	& 0.85	\\
$z_{c}(t=3)$ & 1.53815& 1.3985 \\

\hline   
\end{tabular}
\captionof{table}{Some characteristic values for the rising drop in 3D, Example \ref{sec:risingdrop3D}.  \label{table - result rising drop 3D}}
\end{table}

\section{Conclusions}\label{sec:con}
In this work our focus was on developing an accurate discretization for the surfactant transport equation coupled to the incompressible Navier--Stokes equations. The proposed discretization of the surface PDE and the bulk PDEs utilize the same computational mesh which does not need to conform to the evolving interfaces. We have presented a new space-time weak formulation for the surfactant transport equation which results in discrete conservation of surfactant mass without enforcing the mass with a Lagrange multiplier. In this paper we considered insoluble surfactants but in future work we will also consider soluble surfactants. In the presented scheme the computation of the mean curvature vector is avoided and this simplifies the computations. The presented space-time scheme is based on using quadrature rule in time to approximate the space-time integrals and has a straightforward implementation. We used linear elements in both space and time and a piecewise linear approximation of the interface from an approximate level set function and observed second order convergence. However, if better approximations of the interface is used the proposed space-time CutFEM method can be used with high order elements and result in higher order discretization than second order.

Our scheme relies on an accurate representation and evolution of the interface. We used the level set method which simplified the representation and evolution of the interface in three space dimensions and in the simulations when coalescence occured. However, we note that with our approximate level set function we do not conserve the total area of the drop and in future work we would like to improve on the representation of the interface.

\section*{Acknowledgement}
This research was supported by the Swedish Research Council Grants No. 2018-05262 and the Wallenberg Academy Fellowship KAW 2019.0190.

\bibliographystyle{elsarticle-num}

\biboptions{sort&compress}
\bibliography{StokesSurfactantArxiv}

\end{document}